\magnification =\magstep1
\baselineskip = 14pt

\hfill Duke Math. J.  76  (1994),  no. 1, 203--272.\break

\vskip .5cm

\centerline {\bf KOSZUL DUALITY FOR OPERADS}

\vskip 1cm

\centerline {\bf Victor Ginzburg and  Mikhail Kapranov}

\vskip 1.5cm

\centerline{\bf Table of Contents}

 \item{\bf 0.} Introduction.
\item{\bf 1.} Operads in algebra and geometry.
\item{\bf 2.} Quadratic operads.
\item{\bf 3.} Duality for $dg$ - operads.
\item{\bf 4.} Koszul operads.

\vskip 1cm

\centerline {\bf 0. Introduction.}

\vskip 1cm

{\bf (0.1)} 
The purpose of this paper is to relate two seemingly disparate
developments. One is the  theory of graph cohomology of Kontsevich
[Kon 2 - 3] which arose out of earlier works of Penner [Pe] and Kontsevich 
[Kon 1]
on
the cell decomposition and intersection theory on the moduli
spaces of curves.
The other is the theory of Koszul duality for 
quadratic associative
algebras  which was introduced by Priddy [Pr]
and has found many applications in homological algebra,
algebraic geometry and
representation theory (see e.g., [Be] [BGG] [BGS] [Ka 1] [Man]).
 The unifying concept here is
that of an operad. 

This paper can be divided into two parts consisting of
chapters 1, 3 and 2, 4, respectively.  The purpose of the first
part is to establish a relationship between operads, moduli spaces
of stable curves and graph complexes. To each operad we associate a 
collection of sheaves on moduli spaces. We introduce, in a natural way,
the  cobar complex of an operad and show that it is nothing but
a (special case of the)
graph complex, and that both constructions  can be interpreted
as the Verdier duality functor on sheaves.

In the second part we introduce a class of operads, called
quadratic,  and introduce a distinguished subclass of Koszul operads.
The main reason for introducing
Koszul operads (and in fact for writing this paper) is that 
most of the operads "arising from nature" are Koszul, cf. (0.8)
below. 
 We define a natural duality on quadratic operads
(which is analogous to the duality of Priddy [Pr] for quadratic
associative algebras) and
show that it is intimately related  to the cobar - construction, i.e.,
to graph complexes.  

\vskip .3cm

  {\bf (0.2)} Before going further into discussion of the
results of the paper, let us make some comments for the 
reader not familiar with the
notion of an operad.
Operads were introduced by J.P. May [May] in 1972
for the needs of homotopy theory.
Since then it was  gradually realized that this concept has in fact 
a  fundamental significance for mathematics in general. 
 From an algebraic point
of view, an operad is a system of data that formalizes properties
of a collection of maps
$\,X^n \rightarrow X\,$, a certain set for each $n = 1,2, ...$,
which are closed under permutations of arguments of the maps 
and under
all possible superpositions. 
Such a collection may be
generated by iterated
compositions of some primary maps, called generators of
the operad, and the
whole structucture of the operad may be determined,
in principle, by giving
the the list of relations among the generators. 
This
is very similar to defining  a group by generators and
relations. As for  groups, the consideration of the whole
operad and not only generators and relations presents
several obvious advantages. For instance, we can develop 
``homological algebra'' for operations, i.e., study higher syzygies.
 Note, in particular, that all the  types of algebras 
encountered in practice
(associative, Lie, Poisson,  Jordan etc.) are governed by suitable operads
(we take the binary operations involved as generators and the
identities which they satisfy as relations).

The place of operads
among other structures may be illustrated (very roughly)
by the following table.

 \vskip .4cm
 \halign{ \vbox to .6cm{}\quad#\quad&\vrule width .8pt\vbox to
 .6cm{}\quad#\quad&\vrule\vbox to .6cm{}\quad#\quad&\vrule\vbox
 to .6cm{}\quad#\quad&\vrule\vbox to .6cm{}\quad#\quad\cr
& {\bf Algebra}& {\bf Geometry}& \quad\quad
$\matrix{{\bf Linear}\cr {\bf Physics}}$
\quad\quad& $\matrix {{\bf Non - linear}\cr {\bf Physics}}$\cr
\noalign{\hrule height .8pt}
{\bf Spin 1}& Modules & Vector bundles & $\matrix {{\rm Maxwell}\cr{\rm
  equations}}$ & $\matrix{{\rm Yang 
- Mills} \cr{\rm theory}}$
\cr
\noalign{\hrule}
{\bf Spin 2}& Algebras& Manifolds & $\matrix{
{\rm Linear \quad gravity}\cr{\rm  equations}}$ &
 $\matrix{ {\rm Einstein}\cr
 {\rm gravity}}$\cr
\noalign{\hrule}
{\bf Spin 3}& Operads &{\bf ?}(Moduli spaces) & 
$\matrix{{\rm Rarita - Schwinger}\cr
{\rm  equations}}$& 
{\bf ?}$\pmatrix{ {\rm Conformal}\cr {\rm field  \,\,\,\,theory}}$ \cr}

\vskip .6cm

Question marks indicate that the name put at the corresponding 
square of the table is just the first approximation to an unknown 
ultimate name. The relation of operads to algebras in the first
column is similar to the relation of algebras to modules. 
Given an algebra $A$,
one has a notion of an $A$-module. Similarly, given an operad
${\cal P}$, there is a notion of a ${\cal P}$-algebra. Further,
for any ${\cal P}$-algebra $A$, there is a well defined
abelian category of $A$-modules, see \S 1.6 below.
This explains the hierarchy of the first column of the above table.
In the geometric column, vector bundles on a fixed manifold
correspond to modules over the 
(commutative) algebra of functions
on the manifold. Similarly, giving an operad is analogous
to fixing the class of geometric objects we want to consider.
This is equivalent to studying
the moduli
space of these geometric objects.
The relevance of the notion of an
operad to conformal field theory can be justified  by the following
fact which plays a crucial role in the
theory of operads and in the present paper in particular:

\vskip .3cm

\item{}{\it
The collection $\,{\cal M}=\{{\overline {\cal M}_{0,n+1}} 
\,,\, n=2,3,\ldots\}\,$
(Grothendieck-Knudsen moduli spaces
 of stable genus 0 curves with $n+1$ punctures) has
a natural structure of an operad of smooth manifolds.}

\vskip .3cm

  {\bf (0.3)} The paper is organized as follows. 
Chapter 1 begins with introducing
a category of trees that plays the key role (cf. [BV])
throughout the paper. We then
recall the definition of
an operad and produce a few elementary examples of
operads. Next, the above mentioned operad ${\cal M}$ formed by
Grothendieck-Knudsen moduli spaces is described in some detail.
The significance of this operad for the whole theory of operads
is explained: any operad can be described as a collection of
sheaves on $\cal M$.  

\vskip .3cm

  {\bf (0.4)}
Chapter 2 is devoted to {\it quadratic operads}, the ones generated
by binary operations subject to relations involving three arguments
only. Most of the structures that one encounters in algebra, e.g.,
associative, commutative, Lie, Poisson, Jordan, etc. algebras correspond
to quadratic operads.

Given a quadratic operad ${\cal P}$ , we define the {\it 
quadratic
dual} operad ${\cal P}^!$ in the way
 analogous to the definition of quadratic duality (Priddy)
of quadratic associative algebras. 
In particular, the operads ${\cal C}om$
governing commutative and ${\cal L}ie$ (governing Lie algebras)
 are quadratic dual to each
other. In some informal sense, as a correspondence between 
the categories  of (differential
graded) commutative and Lie algebras, this relation goes back
at least to the work of Quillen [Q 1-2]
and Moore [Mo]. Our theory exhibits a very simple
and precise algebraic fact which is the reason for this relation. 
The operad ${\cal A}s$ 
describing associative (not necessary commutative) algebras
is self - dual in our sense. 

There is a natural concept of a quadratic algebra over a quadratic
operad and for the dual quadratic operads $\cal P$ and ${\cal Q}$ we
 construct a duality
between quadratic (super-) algebras over $\cal P$ and $\cal Q$.
For the case of associative algebras (whose operad is self - dual)
we recover the construction of Priddy. 

Quadratic operads have several nice features. In \S 2.2
we introduce on the category of such operads the internal $hom$
in the spirit of Manin [Man]. We show that the quadratic duality
 can be interpreted as $hom (-, {\cal L}ie)$ where ${\cal L}ie$ is the
Lie operad which  therefore plays the role of
 a dualizing object in our theory.

\vskip .3cm

  {\bf (0.5)} 
In chapter 3 we introduce a contravariant duality functor
${\bf D}$ on the category of
differential graded
($dg$ -) operads (as opposed to the quadratic duality functor
${\cal P} \mapsto {\cal P}^!\,$ studied in the previous chapter).
 We present various approaches to duality. From the
algebraic point of view, the duality ${\bf D}$ is an analogue of the
cobar constuction and a generalization of
the tree part of graph complexes.
 From the geometric point of view, the duality
is an analogue of the Verdier duality for sheaves. In more detail,
let ${\cal W}_{n+1}$  be the moduli space of $n$ - labelled trees
with metric, cf. [Pe] [Kon 2].
This is a contractible topological space with
a natural cell decomposition. This cell decomposition is dual, in 
a sense, to the canonical stratification of ${\overline {\cal M}_
{0,n+1}}$. Moreover, the collection  ${\cal W}= \{{\cal W}_{n+1}\}$
forms a {\it co-operad} which plays the role dual to that of the
operad ${\cal M}$. We show that any $dg$ - operad gives rise
to a compatible collection of constructible complexes on the
spaces ${\cal W}_{n+1}$, one for each $n$.
Furthermore the collection arising from
the ${\bf D}$-dual $dg$ - operad turns out to be
 formed by the Verdier duals of 
the complexes corresponding to the original operad.
The spaces ${\cal W}_{n+1}$ are similar in nature to
Bruhat-Tits buildings, and there is yet another description of
duality in terms of sheaf cohomology, which is reminicent
of the Deligne-Lusztig duality [DL 1-3] for representations of finite
Chevalley groups. 

Next, to any $dg$ -operad we associate its 
generating function
which is in fact (see Definition 3.1.8) a certain formal map
${\bf C}^r \rightarrow {\bf C}^r$. It turns out (Theorem 3.3.2) that for
$dg$ - operads dual in our sense their generating maps are, up to
signs, composition inverses to each other.
Recall (see, e.g., [BGS] [L\"o]) 
that for associative algebra $A$ over a field $k$ 
its 
 cobar - construction(i.e., a suitable $dg$ - model for the Yoneda
algebra Ext$^\bullet_A(k,k)$) has generating function which is
multiplicative inverse to the generating function of $A$. 
The role of composition (change of coordinates on a manifold) versus 
multiplication  (change of coordinates in a vector bundle) for generating 
functions of operads  also fits nicely into the table above.

\vskip .3cm

 {\bf (0.6)} Chapter 4 is devoted to 
{\it Koszul operads}, the quadratic  operads whose
quadratic dual is quasi-isomorphic (canonically)
to the ${\bf D}$-dual. We prove that
the operads ${\cal A}s, {\cal C}om$ and ${\cal L}ie$ are Koszul.
 Associated to any
qudratic operad is its Koszul complex. We show that a quadratic 
operad ${\cal P}$ is
Koszul if and only if its Koszul complex is exact, which is equivalent
also to  vanishing  of higher homology for free ${\cal P}$-algebras.
Given a Koszul operad ${\cal P}$ , we introduce the notion of a
{\it Homotopy ${\cal P}$-algebra} which reduces in the special
cases of Lie and Commutative algebras to that introduced earlier
 by  Schlessinger and Stasheff
[SS] and exploited in an essential way by
Kontsevich [Kon 2]. In fact, Koszul operads provide
 the most natural framework
for the "formal non-commutative geometry".

\vskip .3cm

  {\bf (0.7)} The concept of  an operad in its present form
 had several important precursors. 
One should mention the formalism of "theories" of Lawvere (see  [BV]) and
the  pioneering work of Stasheff [St] on homotopy associative $H$ - spaces.
 A little earlier, in the 1955 paper [L], Lazard considered what
we would now call formal groups in an operad. He used the notion
of "analyseur" which is essentially equivalent 
(though formally different) to the modern
notion of an operad. In (2.2.14) we give a natural interpretation of 
Lazard's  ``Lie theory" for
formal groups in  analyseurs
in terms of Koszul duality .
 We are grateful to Y.I.Manin for pointing 
out to us the reference [L].

In late 1950's Kolmogoroff and Arnold [Kol] [Ar 1] have studied, in 
connection with Hilbert's 13 - th problem, what in our present language
is the operad of continuous operations ${\bf R}^n \rightarrow {\bf R}$.
It was proved in these papers that 
 any continuous function in $n\geq 3$
variables can be represented as a superposition of continuous functions
in only one and two variables, i.e., the  above operad is generated by 
 unary and binary operations. 
From the point of view of the present paper (0.2),
 this result raises an interesting
question whether all the {\it relations} among the binary generators
follow from those provided by functions of three variables. It
seems that this question has never been addressed.
We would like to thank I.M. Gelfand for drawing our attention to
 Kolmogoroff's
work.

\vskip .3cm

  {\bf (0.8)} Our interest in this subject originated
 from an attempt
to explain
a striking similarity between combinatorics involved in the 
Graph complex, introduced by Kontsevich in his work
on Chern-Simons theory, and combinatorics of the 
Grothendieck-Knudsen moduli spaces used by A. Beilinson and
the first author in their work on local geometry of moduli spaces of
$G$-bundles [BG 1]. In particular, the main motivation for the
study of Koszul operads begun in this paper is the investigation
of the operads formed by Clebsch - Gordan spaces, see
(1.3.12) below, for representations
of quantum groups and affine Lie algebras. In a future publication
we plan to show that these operads are Koszul. 
In short, Koszul algebras should be replaced by Koszul operads whenever
the category under considerations has a tensor - type
(e.g., fusion) structure.

\vskip .3cm

  {\bf (0.9)} We are very much indebted to Maxim Kontsevich whose
ideas stimulated most of our constructions. We are also very
grateful to Ezra Getzler who informed us about
his work  in progress
with J.D.S. Jones [Ge J] and
 was the first to suggest  that
the constructions we were working with were related to operads,
the notion unknown to the first author at the time (June 1992).
His remarks helped to clarify several important points. 
We much benefited from conversations with Sasha Beilinson, whose
current joint work with the first author (see [BG 1], [BG 2])
is closely related to the subject. 
 We would  like to thank I.M. Gelfand, Y.I. Manin, J.P. May,
V.V. Schechtman and
J.D.  Stasheff for the discussions of the results of this paper.
Several people kindly responded to the call for comments to the
preliminary version. In particular, we are indebted to D. Wright for
 the references [MTWW] and [W] enabling us to
give a proper attribution of Theorem 3.3.9 and to A.A. Voronov for
correcting our use of det - spaces.
 Special thanks are due to J.D. Stasheff for pointing out
numerous inaccuracies in the earlier version.

\hfill\vfill\eject

\centerline{\bf 1. OPERADS IN ALGEBRA AND GEOMETRY}

\vskip 1.5 cm

\centerline  {\bf 1.1. Preliminaries on trees. }

\vskip 1cm

\noindent {\bf (1.1.1)} By a {\it tree} we   mean in this paper a 
non-empty connected oriented graph $T$ 
without loops (oriented or not) with the following property:
there is at least one incoming edge 
and exactly one  outgoing edge at each vertex of $T$, see Fig. 1a.

\vbox to 4cm{}

Trees are viewed as abstract graphs (1-dimensional topological spaces), a 
plane picture being irrelevant.  We allow some edges of
a tree to be  
bounded by a vertex at  one end only. Such edges will
be called {\it external}. All other edges (those bounded by
vertices at both ends) will be called {\it internal}. 
Any tree
has a unique outgoing external edge, called the {\it output} or the
{\it root} of the tree, and several ingoing external edges, called
{\it inputs} or {\it leaves} of the tree. Similarly, the edges
going in and out of a vertex $v$ of a tree will be referred to
as inputs and outputs at $v$.
A tree with possibly several inputs and   a single vertex
is called a {\it star}. There is also a tree (see Fig. 1b)
with a single input and without vertices called the {\it
degenerate} tree. 

\vskip .2cm

We   use the notation ${\rm In}(T)$ for the set of input edges
of a tree $T$; for any vertex $v\in T$ we denote by ${\rm In}(v)$
the set of input edges at $v$. Similarly, we denote by ${\rm Out}(T)$
the unique output edge (root) of $T$ and for every vertex
$v\in T$ we denote by ${\rm Out}(v)$ the output edge at $v$.

\vskip .2cm

 Let $I$ be a finite set. A tree $T$ equipped with a bijection
between $I$ and the set ${\rm In} (T)$ will be
referred to as $I$ - {\it labelled} tree or
an $I$-tree for short. Two $I$ - trees $T, T'$ are called
{\it isomorphic} if there exists an isomorphism of trees $T 
\rightarrow T'$ preserving orientations and 
the labellings of the inputs. 

We   denote by $[n]$ the finite set $\{1, 2, ..., n\}$
and call $[n]$ - trees simply $n$ - trees. 

\vskip .3cm

{\noindent \bf (1.1.2) Composition.}
 Let $I_{1}$ be a set
and let $I_{2}$ be another set with a marked element 
$i \in I_{2}\,$. Define the  composition 
of $I_{1}$ and  $I_{2}$ along $i$  as $I = I_{1} \circ_i I_{2} =
I_{1} \cup (I_{2} \,\setminus\,\{i\})$.
 Given an $I_{1}$-tree $T_{1}$ and an $I_{2}$-tree
$T_{2}$ , let $T=T_{1} \circ_i T_{2}\,$ be the $I$-tree obtained by
identifying the output of $T_{1}$ with the $i$ - th 
input of $T_{2}$  as depicted in Fig. 2. 

\vbox to 4cm{}

The tree $T$ is called the {\it composition} of $T_{_1}$ and $T_{_2}$
(along $i$).
Note that any tree can be obtained as an iterated composition of stars.

\vskip .2cm

Let $T$ be a tree and $v\buildrel e\over\longrightarrow w$  an
internal edge of $T$. Then we can form a new tree $T/e$ by contracting
$e$ into a point. This new point is a vertex of $T/e$ denoted by
$<e>$. Clearly, we have
$$ {\rm In} (<e>) = {\rm In} (w) \circ_e {\rm In}(v).\leqno 
{\bf (1.1.3)}$$
If $T$ is $I$ - labelled then so is $T/e$. If $T, T'$ are two
$I$ - labelled trees then we write
 $T\buildrel e\over\longrightarrow T'$ if the tree $T'$ is
isomorphic (as a labelled tree) to $T/e$. 

 \vskip .2cm

Write $T \geq T'$ if there is a sequence of edge contractions
 $\,T \rightarrow
T_{1} \rightarrow \ldots \rightarrow
T_{k} \rightarrow T'\,$. Thus, $\geq\,$ is 
a partial order on the set of all
$I$-trees. An $I$-star is the  unique minimal  element with respect
to that order.

\vskip .2cm

 A tree $T$ (with at least one vertex)
is called a {\it binary tree} if there are exactly 2 inputs
at each vertex of $T$. Binary trees are the maximal elements with respect
to the partial order $\leq$. It can be shown that 
the number of non-isomorphic binary
$n$ - trees is equal to $(2n-3)!! = 
1\cdot 3\cdot 5 \cdot ... \cdot (2n-3)$. A simple proof of this fact
using generating functions will be given in Chapter 3.

\vskip .3cm

\noindent {\bf (1.1.4) The category of trees.} Let 
$T, T'$ be trees (viewed as 1-dimensional topological spaces).
 By a {\it morphism} from $T$ to $T'$ we understand a
continuous surjective map $f: T \rightarrow T'$ 
with the following properties:

\item{(i)} $f$ takes  each vertex to
a vertex and each edge into an edge or a vertex. 

\item{(ii)} $f$ is monotone, i.e., preserves the orientation.

\item{(iii)} The inverse image of any point of $T'$ under $f$ is
a connected subtree in $T$.

 Thus any morphism is a composition of an isomorphism 
and several edge contractions.
In this way we get
a category which we denote $Trees$.

\vskip .3cm

\noindent {\bf (1.1.5)} Let $\Sigma_n$ denote the symmetric group
of order $n$. For any two sets $I, J$  of the same cardinality
we denote by Iso$(I,J)$
the set of all bijections $I \rightarrow J$. We write $\Sigma_I$ for
Iso$(I,I)$, so  that $\Sigma_n =
{\rm Iso}([n], [n])$. Clearly, Iso$(I,J)$ is a principal homogeneous left
$\Sigma_I$ - space and a principal homogeneous right $\Sigma_J$ -
space.

Let $W$ be a vector space (over some field $k$) with an action of
$\Sigma_n$. There is a canonical way to construct (out of
$W$) a functor $I \mapsto W(I)$ from the category of $n$ - element
sets and  bijections to the category of $k$ vector
spaces. Namely, put
$$W(I) = \left( \bigoplus_{f\in {\rm Iso}([n], I)} W\right)_{\Sigma_n},
\leqno {\bf (1.1.6)}$$
the coinvariants with respect to the simultaneous action
of $\Sigma_n$ on Iso$([n], I)$ and $W$.
The original space   $W$ is recovered
 as the value of this functor on the set $[n]$.
Similarly, if $W$ is a set, or topological space with $\Sigma_n$ -
action, we can construct a functor $I \mapsto W(I)$
from the category of $n$ - element sets and their bijections
 to the category of
sets, or topological spaces as above, by the following analog of (1.1.6):
$$W(I) = {\rm Iso}([n], I) \times_{\Sigma_n} W.$$

\hfill\vfill\eject

\centerline {\bf 1.2. $k$ - linear operads.}

\vskip 1cm

\noindent {\bf (1.2.1)} Let $k$ be a field of characteristic 0.
A $k$ - linear operad ${\cal P}$ is a collection $\{{\cal P}(n),
n\geq 1\}$ of $k$ - vector spaces equipped with the following set
of data:

\item{(i)} An action of the symmetric group $\Sigma_n$ on ${\cal P}
(n)$  for each $n\geq 1$.

\item{(ii)} Linear maps (called compositions)
$$\gamma_{m_1,...,m_l}: {\cal P}(l) \otimes {\cal P}(m_1) \otimes ...
\otimes {\cal P}(m_l) \longrightarrow {\cal P}(m_1+...+m_l)$$
 for all $m_1,...,m_l \geq 1$. We write $\mu(\nu_1, ..., \nu_l)$
instead of 
$\gamma_{m_1, ..., m_l}(\mu \otimes \nu_1 \otimes ... \otimes \nu_l)$.

\item{(iii)} An element $1\in {\cal P}(1)$, called the unit, such that
$\mu(1,...,1) = \mu$ for any $l$ and any $\mu\in {\cal P}(l)$.

\vskip .2cm

It is required that these data satisfy the conditions (associativity and
equivariance with respect to symmetric group actions) specified by
May ([May], \S 1). These conditions are best expressed in terms of
trees.
Observe first that the datum (i) allows to
assign to any finite set $I$ a vector space ${\cal P}(I)$ as in (1.1.5).
Next, we associate to any tree $T$ the vector space
$$\tilde{\cal P}(T) = \bigotimes_{v\in T} {\cal P}({\rm In}(v)).
\leqno {\bf (1.2.2)}$$
To the degenerate tree without vertices we associate, by
definition, the field $k$. 

Note in particular that to the trees $T(n)$ and $T(m_1,...,m_l)$
depicted in Fig.3 we associate the spaces ${\cal P}(n)$
and ${\cal P}(l) \otimes {\cal P}(m_1) \otimes ... \otimes {\cal P}(m_l)$.

\vbox to 4cm{}

Thus the datum (ii) gives a map
$$\tilde{\cal P}(T(m_1,...,m_l)) \longrightarrow
\tilde{\cal P}(T(m_1+...+m_l)) = {\cal P}(m_1 + ...+ m_l).
\leqno {\bf (1.2.3)}$$
For any $n$ -  tree $T$ there exists a sequence of trees
$$T = T_0 \rightarrow T_1 \rightarrow ... \rightarrow T_r = T(n)
\leqno {\bf (1.2.4)}$$ 
where each $T_i$ is obtained from $T_{i-1}$ by replacing a fragment of type
$T(m_1,...,m_l)$ by $T(m_1 + ... + m_l)$. So maps (1.2.3) give rise
to a sequence of maps
$$\tilde{\cal P}(T) = \tilde{\cal P}(T_0) \rightarrow\tilde {\cal P}(T_1) 
\rightarrow ... \rightarrow \tilde{\cal P}(T_r) = {\cal P}(n).$$
The associativity condition is equivalent to the requirement that
the composite map $\tilde {\cal P}(T) \rightarrow P(n)$
does not depend on the choice of a sequence (1.2.4). 

\vskip .3cm

\noindent {\bf (1.2.5)} Observe that the tree $T(m_1+...+m_l)$ is
obtained from $T(m_1, ...,   m_l)$ by contracting all the
$l$ internal edges. The existence of unit $1\in {\cal P}(1)$
makes in possible to decompose the map (1.2.3) corresponding to this
contraction, into more elementary ones, each consisting of
contracting a single edge.

Namely, let $T$ be a tree, $v \buildrel e\over\longrightarrow w$
be an internal edge of $T$
and $T/e$ be the tree obtained by contracting $e$. Let 
$<e>$ be the vertex of $T/e$ obtained from the contracted edge.
We define a map
$${\cal P}({\rm In}(v)) \otimes {\cal P}({\rm In}(w)) \longrightarrow
{\cal P}\left( {\rm In}(v) \circ_e {\rm In}(w)\right) 
\buildrel (1.1.3)\over = {\cal P}({\rm In}(<e>))$$
by the formula
$\mu \otimes \nu \mapsto \mu(1,..., \nu, ..., 1)$ where $\nu$ 
is placed at the entry corresponding to the edge $e$.
By tensoring this map with the identity elsewhere on $T$, we obtain
a map
$$\tilde\gamma_{_{T,e}}: \tilde {\cal P}(T) \rightarrow \tilde
{\cal P}(T/e).\leqno {\bf (1.2.6)}$$
More generally, if $I$ is a finite set and $T, T'$ are two $I$ - 
trees such that $T\geq T'$ then by composing maps of the 
type $\tilde\gamma_{_{T,e}}$ we get a map
$$\tilde\gamma_{_{T,T'}}: \tilde {\cal P}(T) \rightarrow \tilde {\cal P}(T')
\leqno {\bf (1.2.7)}$$
which is well - defined due to  the associativity condition.

Thus a $k$ - linear operad ${\cal P}$ gives rise to a functor
$$\tilde {\cal P}: Trees \longrightarrow {\rm Vect},\quad T \mapsto 
\tilde {\cal P}(T)$$
equipped with the following additional structures:

\item{(i)} For any trees $T_1$ and $T_2$ and any $j\in {\rm In}(T_2)$
one has a functorial isomorphism
$$\Phi_{T_1, T_2}: \tilde{\cal P}(T_1) \otimes \tilde{\cal P}(T_2) \rightarrow
\tilde{\cal P}(T_1 \circ_j T_2). \leqno {\bf (1.2.8)}$$

\item{(ii)} The isomorphisms in (i) satisfy the  associativity constraint
saying that for any trees $T_1, T_2, T_3$ and any $i\in {\rm In}(T_2),
j\in {\rm In}(T_3)$ the following diagram commutes:

$$\matrix{
&\tilde{\cal P}(T_1) \otimes\tilde {\cal P}(T_2) \otimes 
\tilde{\cal P}(T_3) &
\buildrel {\rm Id} \otimes \Phi_{T_2, T_3} \over
\longrightarrow &\tilde{\cal P}(T_1) \otimes
\tilde{\cal P}(T_2\circ_j T_3) &\cr
\Phi_{T_1, T_2}\otimes {\rm Id} & \big\downarrow &&\big\downarrow &\Phi_{T_1, 
T_2\circ T_3}\cr
&\tilde{\cal P}(T_1 \circ_i T_2) \otimes\tilde {\cal P}(T_3) &
\buildrel \Phi_{T_1\circ T_2, T_3} \over\longrightarrow &
\tilde P(T_1\circ_i T_2\circ_j T_3)&}$$

\vskip .2cm

\noindent {\bf (1.2.9)} Let $V$ be a $k$ - vector space. Its {\it
operad of endomorphisms}, ${\cal E}_V$, consists of vector spaces
$${\cal E}_V(n) = {\rm Hom}(V^{\otimes n}, V).$$
with compositions
and the $\Sigma_n$ - action on ${\cal E}_V(n)$ being defined in an
obvious way. We have  ${\cal E}_V(1) = {\rm End}(V)$.

\vskip .3cm

\noindent {\bf (1.2.10)} Observe that for any $k$ - linear operad ${\cal P}$
the space $K = {\cal P}(1)$ has a natural structure of an
associative $k$ - algebra with unit. Conversely, if $K$ is a $k$ -
algebra then the collection $\{{\cal P}(1) = K, {\cal P}(n) = \{0\},
n>1 \}$ forms an operad. Furthermore, for an
arbitrary $k$ - linear operad ${\cal P}$
the space ${\cal P}(n)$ has several ${\cal P}(1)$ - module structures.
These structures are summarized in the following definition.

\proclaim (1.2.11) Definition. Let $K$ be an associative $k$ - algebra
with unit. A $K$ - collection is a collection $E = \{E(n), n\geq 2\}$ of $k$
- vector spaces equipped with the following structures:

\item{(i)} A left $\Sigma_n$ - action on $E(n)$,  for each $n\geq 2$.

\item{(ii)} A structure of a left $K$ - module and a right $K^{\otimes n}$ -
module on $E(n)$, $n\geq 2$.

\noindent These structures are required to satisfy the following compatibility
condition: for any $s\in \Sigma_n$ and any $\lambda_1,...,\lambda_n \in K,
a\in E(n)$ we have
$$s\left(a\cdot (\lambda_1\otimes ... \otimes \lambda_n)\right) =
a\cdot (\lambda_{s(1)} \otimes ... \otimes \lambda_{s(n)}).$$

If ${\cal P}$ is a $k$ - linear operad and $K = 
{\cal P}(1)$ then $\{{\cal P}(n), n\geq 2\}$ is, clearly, a $K$ - collection.

\vskip .3cm

\noindent {\bf (1.2.12)} It will be convenient for the future purposes to
give a reduced, in a certain sense, version of the tree formalism above.

We call a tree $T$ {\it reduced} if there are at least two inputs at
each vertex $v\in T$ (the degenerate tree without vertices is also
assumed to be reduced). 

Let $K$ be an associative $k$ - algebra and $E$  a $K$ - collection.
As in (1.1.5) we extend $E$ to a functor on finite sets and bijections.
To each reduced tree $T$ we associate a vector space $E(T)$ as follows.
For the degenerate tree $\rightarrow$ we set $E(\rightarrow) = K$
 and for a tree $T$ with non - empty set of vertices we set
$$E(T) = \bigotimes_{v\in T}{}_{_K} E({\rm In}(v)). \leqno{\bf (1.2.13)}$$
This tensor product is taken over the ring $K$ by using the 
$(K, K^{\otimes {\rm In}(v)})$ - bimodule  structure  on
 each $E({\rm In}(v))$,
see Fig.4.

\vbox to 4cm{}

Explicitly, this means that $E(T)$ is the quotient of the $k$ - tensor
product $\bigotimes_v E({\rm In}(v))$ by associativity conditions
described as follows.
Suppose that $e\mapsto \lambda_e$ is any $K$ - valued function
on the set of all edges such that $\lambda_e = 1$ for all external edges.
Then for any collection of $a_v \in  E({\rm In}(v))$ we impose the
relation
$$\bigotimes_{v\in T} a_v \cdot \biggl(\bigotimes_{e\in {\rm In}(v)}
\lambda_e\biggl) = \bigotimes_{v\in T} \lambda_{{\rm Out}(v)}\cdot a_v.$$

\vskip .2cm

\noindent {\bf (1.2.14)} Now let ${\cal P} = \{{\cal P}(n)\}$ be a 
$k$ - linear operad and $K = {\cal P}(1)$. Since $\{{\cal P}(n), n\geq 2\}$
form a $K$ - collection, we can assign to any reduced tree $T$ a vector space
${\cal P}(T)$ by formula (1.2.13). This assignment has the following two
fundamental structures:

\item{(i)} A linear map
$$\gamma_{T,T'}: {\cal P}(T) \rightarrow {\cal P}(T') \leqno {\bf (1.2.15)}$$
defined whenever $T\geq T'$.

\item{(ii)} An isomorphism
$${\cal P}(T_1\circ_i T_2) \rightarrow {\cal P}(T_1) \otimes_K
{\cal P}(T_2)$$
 given for any $i\in {\rm In}(T_2)$.

\noindent The maps in (i) and (ii) are  induced by (1.2.7)
and (1.2.8), respectively.

Observe further that the assignment $T \rightarrow {\cal P}(T)$
extends to a functor on the (sub) category of reduced trees, see
(1.1.4).

\vskip .3cm

\noindent {\bf (1.2.16) Convention.} In the rest of this paper we
shall consider only reduced trees, unless specified otherwise. 

\hfill\vfill\eject

\centerline {\bf 1.3. Algebraic operads.}

\vskip 1cm

\noindent {\bf (1.3.1)} Let ${\cal P}, {\cal Q}$ be two $k$ - linear
operads. A {\it morphism of operads} $f: {\cal P} \rightarrow
{\cal Q}$ is a collection of linear maps which are equivariant with respect
to $\Sigma_n$ - action, commute with compositions $\gamma_{m_1,...,m_l}$
in ${\cal P}$ and ${\cal Q}$ and take the unit of ${\cal P}$ to
the unit of ${\cal Q}$, see [May].

\vskip .2cm

\proclaim (1.3.2) Definition. Let ${\cal P}$ be a $k$ - linear operad.
A ${\cal P}$ - algebra is a $k$ - vector space $A$ equipped with a 
morphism of operads $f: {\cal P} \rightarrow {\cal E}_A$ where 
${\cal E}_A$ is the operad of endomorphisms of $A$, see (1.2.9).

Clearly, giving a structure of a ${\cal P}$ - algebra on $A$ is the same
as giving a collection of linear maps
$$f_n: {\cal P}(n) \otimes A^{\otimes n} \longrightarrow A\leqno {\bf
(1.3.3)}$$
satisfying natural associativity, equivariance and unit conditions. We
  write \hfill\break
$\mu (a_1, ... , a_n)$ for $f_n(\mu \otimes (a_1 \otimes ...
\otimes a_n))$, $\mu\in {\cal P}(n), a_i\in A$.

\vskip .3cm

\noindent {\bf (1.3.4) Free algebras.} Let ${\cal P}$ be a $k$ - linear
operad and $K = {\cal P}(1)$. As noted in (1.2.10), $K$ is an associative
$k$ - algebra and every ${\cal P}(n)$ is a left $K$ - module and a right
$K^{\otimes n}$ - module. 

Let $V$ be any left $K$ - module. Form the following graded vector space
$$F_{\cal P}(V) = \bigoplus_{n\geq 1} \biggl(
 {\cal P}(n) \otimes_{K^{\otimes n}}
 V^{\otimes n}
\biggl)_{\Sigma_n} \leqno {\bf (1.3.5)}$$
where $V^{\otimes n}$ is the $n$ -th tensor power of $V$ over $k$. 

\proclaim (1.3.6) Lemma. Compositions in ${\cal P}$ induce natural maps
${\cal P}(n) \otimes F_{\cal P}(V)^{\otimes n} \rightarrow F_
{\cal P}(V)$. There maps
make $F_{\cal P}(V)$ into a $\cal P$ - algebra.

We call $F_{\cal P}(V)$ the {\it free ${\cal P}$ - algebra} generated by $V$.
Note that $F_{\cal P}(V)$ has a natural grading given by the decomposition
(1.3.5).

Here are some other examples of operads and algebras:

\vskip .3cm

\noindent {\bf (1.3.7) Associative algebras.} For any $n$ let ${\cal A}s
(x_1,...,x_n)$ denote the free associative $k$ - algebra on generators $x_1,
..., x_n$ (i.e., the algebra of non - commutative polynomials). Let
${\cal A}s(n) \i {\cal A}s(x_1,...,x_n)$ be the subspace spanned by
monomials containing each $x_i$ exactly once. There are exactly $n!$
such monomials namely $x_{s(1)} \cdot ... \cdot x_{s(n)}$, $s\in
\Sigma_n$. Clearly, ${\cal A}s(n)$ has a natural $\Sigma_n$ - action and as a 
$\Sigma_n$ - module it is isomorphic to the regular representation of
$\Sigma_n$. The collection As$ = \{{\cal A}s(n)\}$ forms an operad
called the {\it associative operad}. The composition maps (see (1.2.1))
$${\cal A}s(l) \otimes {\cal A}s(m_1) \otimes ... \otimes {\cal A}s(m_l)
\longrightarrow {\cal A}s(m_1 + ... + m_l)$$
are given by substituting the monomials $\phi_1 \in {\cal A}s(m_1),
..., \phi_l\in {\cal A}s(m_l)$ in place of generators $x_1,...,x_l$ into
a monomial $\psi \in {\cal A}s(l)$. 

We leave to the reader to verify that an ${\cal A}s$ - algebra is nothing
but an associative algebra in the usual sense (possibly without unit). 

\vskip .3cm

\noindent {\bf (1.3.8) Commutative (associative) algebras.} For any $n$
let ${\cal C}om (x_1,...,x_n) = $\hfill\break $= k[x_1, ..., x_n]$ be the 
free commutative algebra on generators $x_1, ..., x_n$ 
(i.e., the algebra of polynomials). There exists precisely one monomial
containing each $x_i$ exactly once, namely $x_1 \cdot ... \cdot x_n$.
Let ${\cal C}om (n) \i {\cal C}om (x_1, ..., x_n)$ be the 1-dimensional
subspace spanned by this monomial. The collection 
${\cal C}om = \{{\cal C}om (n)\}$ forms an operad with respect to the
trivial actions of $\Sigma_n$ on ${\cal C}om (n)$ and compositions
defined similarly to (1.3.7). We call ${\cal C}om$ the 
{\it commutative operad}.

Again, it is straightforward to see that a ${\cal C}om$ - algebra
is nothing but a commutative associative algebra in the usual
sense (possibly without unit).

\vskip .3cm

\noindent {\bf (1.3.9) Lie algebras.} For any $n$ let ${\cal L}ie
(x_1, ..., x_n)$ be the free Lie algebra over $k$ generated by 
$x_1, ..., x_n$. Let ${\cal L}ie (n) \i {\cal L}ie
(x_1, ..., x_n)$ be the subspaces spanned by all bracket monomials 
containing  each $x_i$ exactly once. Note that such monomials are
not all linearly independent due to the Jacobi identity. 
The subspace ${\cal L}ie (n)$ is invariant under the action
of $\Sigma_n$ on ${\cal L}ie (x_1,...,x_n)$ by permutations of $x_i$.
It is known that
$$ {\rm dim} \,\,\, {\cal L}ie (n) = (n-1)! \leqno {\bf (1.3.10)}$$
Moreover, if $k$ is algebraically closed, then A. Klyachko [Kl]
an isomorphism of $\Sigma_n$ - modules
$${\cal L}ie (n) \cong {\rm Ind}_{{\bf Z}/n}^{\Sigma_n} (\chi)
\leqno {\bf (1.3.11)}$$
where $\chi$ is the 1-dimensional representation of the
cyclic group ${\bf Z}/n$ sending the generator into a
primitive $n$ - th root of 1 (the induced module does not
depend on the choice of such a $\chi$). 

The collection ${\cal L}ie = \{{\cal L}ie (n)\}$ forms an operad with
respect to the composition operations defined similarly to the
ones described in (1.3.4). An algebra over the operad ${\cal L}ie$
is the same as a Lie algebra in the usual sense. 

\vskip .2cm

It should be clear to the reader at this point how to construct operads
governing other types of algebras encountered in practice:
Poisson algebras, Jordan algebras etc.

\vskip .3cm

\noindent {\bf (1.3.12)} A collection of finite - dimensional
$k$ - vector spaces
$$\left\{ E^j(a_1, ..., a_r), \quad a_1, ..., a_r \in {\bf Z}_+, \quad
\sum a_i \geq 1, \quad j=1, ...,r\right\}$$
is called an $r$ - fold collection if the following holds:

\item{(i)} For any $a_1, ...,a_r \in {\bf Z}_+$ the space
$E^j(a_1, ..., a_r)$ is equipped with an action of the group
$\Sigma_{a_1} \times ... \times \Sigma_{a_r}$.

\item{(ii)} $$E^j(0, ..., 0, 1, 0, ...,0) \,\,\, (1 \,\,\,
{\rm on\quad the\quad}
j - {\rm th\quad place}) = \cases {$0$ & if $i\neq j$\cr
$k$ & if $i=j$\cr}.$$

Now let $K$ be a semisimple $k$ - algebra and $\{M_1, ..., M_r\}$
be a complete collection of (the isomorphism classes of) simple
left $K$ - modules. Given a $K$ - collection $E = \{E(n), n\geq 2\}$
(1.2.11), we regard $E(n)$ as a left module over the algebra 
$\Pi_n = K^{\otimes n} \otimes K^{op}$. For any $a_1, ..., a_r
\in {\bf Z}_+$ such that $\sum a_i = n$, we define
$$E^i(a_1, ..., a_r) = {\rm Hom}_{\Pi_n} \bigl(
M_1^{\otimes a_1} \otimes ... \otimes M_r^{\otimes a_r}
 \otimes M_i^*, E(n)\bigl).
\leqno {\bf (1.3.13)}$$
The $k$ - vector spaces $E^j(a_1, ..., a_r)$ thus defined clearly form
an $r$ - fold collection. It will be called the $r$ - fold
collection {\it associated} to the $K$ - collection $E$.

\vskip .3cm

\noindent {\bf (1.3.14) Clebsch - Gordan spaces and operads.}
 Let ${\cal A}$
be a semisimple  Abelian $k$ - linear category equipped with
 a symmetric monoidal
structure $\otimes$ (for example, the tensor category of 
finite - dimensional representations of a
finite, or an algebraic group).

Fix any integer $r\geq 1$ and any set $X_1, ..., X_r$ of pairwise
non - isomorphic simple objects  of ${\cal A}$.
Let $X = \bigoplus_{i=1}^r X_i$. Associated to $X$ is the operad
${\cal P}_X$ defined by collection of vector spaces
$${\cal P}_X (N) = {\rm Hom}_{\cal A} (X^{\otimes n}, X).\leqno
{\bf (1.3.15)}$$
Observe that ${\cal P}_X(1) = \bigoplus {\rm End}(X_i)$ is a semisimple
$k$ - algebra. We put $K = {\cal P}_X(1)$ and view the collection
(1.3.15) as a $K$ - collection. Clearly, the $r$ - fold collection
associated to that $K$ - collection via (1.3.13) is formed by the
Clebsch - Gordan spaces
$${\cal P}^i(a_1, ..., a_r) = {\rm Hom}_{\cal A} (X_1^{\otimes a_1} \otimes
... \otimes X_r^{\otimes a_r}, \,\, X_i).\leqno {\bf (1.3.16)}$$

\vskip .3cm

\noindent {\bf (1.3.17) Operads in monoidal categories.} Let $Vect$ be the
category of $k$ - vector spaces. Note that the concept of a $k$ -
linear operad  appeals only to the category $Vect$ with its
symmetric monoidal structure given by the tensor product. Clearly, one
can define operads in any symmetric monoidal category, i.e., a category 
${\cal A}$ (not necessarily additive or abelian) equipped with a bifunctor
$\otimes: {\cal A} \times {\cal A} \rightarrow {\cal A}$ and 
natural associativity and commutativity constraints for this functor [Mac].

All the preceding constructions (e.g., the notion of a ${\cal P}$ -algebra)
can be carried over to the setup of operads in any symmetric monoidal
category $(\cal A, \otimes)$. Given such an operad, one defines, as
in (1.1.5), for any finite set $I$, an object ${\cal P}(I) \in {\cal A}$and for  tree $T$  an object
$\tilde {\cal P}(T) \in {\cal A}$.

In this section we concentrate on algebraic examples of categories
${\cal A}$.

\vskip .3cm

\noindent {\bf (1.3.18) Two categories of graded vector spaces.} Let
$gVect^+$ be the category whose objects are graded vector spaces
$V^\bullet = \bigoplus_{i\in {\bf Z}} V^i$ and morphisms are
linear maps preserving the grading. For $v\in V^i$ we write deg$ (v) = i$.
We introduce on $gVect^+$ a symmetric monoidal structure defining 
with the
tensor product 
$$ (V^\bullet \otimes W^\bullet)^m = \bigoplus_{i+j=m} V^i \otimes W^j.
\leqno {\bf (1.3.19)}$$
The associativity morphism $V^\bullet \otimes (W^\bullet \otimes X^\bullet)
\rightarrow (V^\bullet \otimes W^\bullet) \otimes X^\bullet$ takes
$v\otimes (w\otimes x) \mapsto (v\otimes w) \otimes x$. The commutativity
isomorphism $V^\bullet \otimes W^\bullet \rightarrow W^\bullet \otimes
V^\bullet$ takes $v\otimes w \mapsto w\otimes v$.

\vskip .2cm

The symmetric monoidal category $gVect^-$ has the same objects, morphisms,
tensor product and associativity isomorphism as $gVect^+$ but the
commutativity isomorphism is changed to be
$$v \otimes w \mapsto (-1)^{ {\rm deg}(v) \cdot {\rm deg}(w)} w\otimes v.
\leqno {\bf (1.3.20)}$$

Any $k$ - linear operad ${\cal P}$ can be regarded as an operad
in either of the categories $gVect^\pm$ and so we can speak about
${\cal P}$ - algebras in these categories. For example, if ${\cal P} = 
{\cal C}om$ is the commutative operad, then a ${\cal C}om$ - algebra
in $gVect^+$ is a commutative associative algebra equipped with a 
grading compatible with the algebra structure. A ${\cal C}om$ -
algebra in $gVect^-$ is an associative graded algebra which is
graded (or super-) commutative.

\vskip .3cm

\noindent {\bf (1.3.21) The determinant operad $\Lambda$.} 
In the operad ${\cal C}om$ each space ${\cal C}om(n)$ is 1-dimensional
and equipped with the trivial action of $\Sigma_n$. We  
now introduce an operad $\Lambda$ in the category
$gVect^-$ which is an ``odd" analog of ${\cal C}om$.

\vskip .2cm

We define $\Lambda(n)$ to be the 1-dimensional vector space
$\bigwedge^n(k^n)$ (the sign representation of $\Sigma_n$)
placed in degree $(1-n)$. In order to describe compositions in
$\Lambda$ we first describe what is to  be  a $\Lambda$ - algebra in
$gVect^-$.

We consider graded vector spaces $A = \bigoplus A_i$ with one binary
operation $(a,b) \mapsto ab$ satisfying the following identities:

\item{(i)} ${\rm deg}(ab) = {\rm deg}(a) + {\rm deg}(b) -1$,

\item{(ii)} $ab = (-1)^{{\rm deg}(a) + {\rm deg}(b) + 1} ba$,

\item{(iii)} $a(bc) = (-1)^{{\rm deg}(a) + 1} (ab)c.$

\noindent The operation $ab$ corresponds to the generator $\mu \in
\Lambda(2)$. The condition (i) means that  ${\rm deg}(\mu) = -1$
and the condition (ii) means that $\Sigma_2$ actis on $\mu$ by the
sign representation.

We fix integers $d_1,...,d_n$ and
fet $\Lambda(x_1, ..., x_n)$ be the free algebra with the above
identities generated by symbols $x_i$, of degree $d_i$.

 \proclaim (1.3.22) Lemma. The subspace 
$E_{d_1,...,d_n}$ in $\Lambda(x_1,...,x_n)$ spanned
by non-associative monomials containing each $x_i$ exactly once,
has dimension 1. Moreover, any two such monomials are proportional
with coefficients $\pm 1$.

The meaning of this lemma is that the
 identities (i) - (iii) above are
consistent, i.e., do not lead to the contradiction 1=0. In other
words, any two ways of comparing the signs of the monomials
lead to the same result.

\vskip .2cm

\noindent {\sl Proof of the lemma:}
As in Mac Lane's axiomatics of symmetric 
monoidal categories [Mac], it is enough to check the three
elementary ambiguities, i.e., the two ways of comparing
$a(bc)$ and $(bc)a$, of $(ab)c$ and $c(ab)$ and of
$a(b(cd))$ and $((ab)c)d$. We leave this to the reader.

\vskip .2cm

Another way to see Lemma 1.3.22 can be based on the following
remark (made to us by E. Getzler). If $A$ is a graded
commutative algebra (i.e., a ${\cal C}om$ - algebra in the category
$gVect^-$) then the same algebra with the grading shifted by one
and new multiplication $a*b = (-1)^{{\rm deg}(a)} ab$ will
satisfy the identities (i) - (iii) above.

\vskip .3cm

\noindent {\bf (1.3.23)} To finish the construction of the operad 
$\Lambda$, we take, in the situation of the
 Lemma 1.3.22, $n$ generators $x_1,...,x_n$ of degree 0.
We denote $\Lambda'(n) = E_{0,...,0}$ the 1-dimensional
subspace from this lemma.
 The space
$\Lambda'(n)$ is $\Sigma_n$ - invariant and the action of
$\Sigma_n$ on $\Lambda'(n)$ is given by the sign representation.
To see the last assertion it is enough to show that any transposition
$(ij)$ acts by $(-1)$. But we can take the basis vector of
$\Lambda'(n)$ given by any product $...(x_ix_j)...$ in which
$x_i$ and $x_j$ are bracketed together. By (ii) we have $x_ix_j = -
x_jx_i$ whence the assertion. 

So we can identify $\Lambda'(n)$ with $\Lambda(n) = \bigwedge^n(k^n)$
and the substitution of monomials in place of generators, as in (1.3.7),
defines the operad structure on $\Lambda$.

\hfill\vfill\eject

\centerline {\bf 1.4. Geometric operads.}

\vskip 1cm

\noindent {\bf (1.4.1)} The category of topological space
has an obvious symmetric monoidal structure given by the Cartesian
product. Operads in this category will be called {\it
topological operads}. (This was the original context of  [May].)

Given a topological operad ${\cal P}$, the total homology spaces
$H_\bullet({\cal P}(n), k)$ form an operad in the category
$gVect^-$ (by K\"unneth formula). Furthermore, for any $q\geq 0$, the
subspaces $H_{q(n-1)}({\cal P}(n), k)$ form a sub - operad.

\vskip .3cm

\noindent {\bf (1.4.2)} An important example of a topological operad is 
given by the little $m$ - cubes operad ${\cal C}_m$ of Boardman - Vogt -
 May [BV] [May]. By definition, ${\cal C}_m(n)$ is the space of numbered $n$
- tuples of non-intersecting $m$ - dimensional cubes inside the standard
cube $I^n$, with faces parallel to those of $I^m$. The operad
${\cal C}_2$ of little squares has an algebro - geometric analog
which will be particularly interesting for us and which we proceed to
describe.

\vskip .3cm

\noindent {\bf (1.4.3) The moduli space ${\cal M}(n)$.}
Let $M_{0, n+1}$ be the moduli space of $(n+1)$ - tuples 
$(x_0, ..., x_n)$ of distinct points on the complex projective line
${\bf C}P^1$ modulo projective authomorphisms. Choose a point
$\infty \in {\bf C}P^1$ so ${\bf C}P^1 = {\bf C} \cup \{\infty\}$.
Letting $x_0 = \infty$, one gets an isomorphism of $M_{0, n+1}$
with the moduli space of $n$ - tuples of distinct points on {\bf C}
modulo affine automorphisms.

\vskip .2cm

The space $M_{0,n+1}$ has a canonical compactification ${\cal M}(n) 
\supset M_{0,n+1}$ introduced by Grothendieck and Knudsen [De 2] [Kn].
The space ${\cal M}(n)$ is the moduli space of stable $(n+1)$ -
pointed curves of genus 0, i.e., systems $(C, x_0, ..., x_n)$ where $C$
is a possibly reducible curve (with at most nodal singularities)
and $x_0, ..., x_n \in C$ are distinct smooth points such that:

\vskip .1cm

\item{(i)} Each component of $C$ is isomorphic to ${\bf C}P^1$.

\item{(ii)} The graph of intersections of components of $C$ is a tree.

\item{(iii)} Each component of $C$ has at least 3 special points
where by a special point we mean either one of the $x_i$ or a singular 
point of $C$.

\vskip .1cm

The space $M_{0,n+1}$ forms an open dense part of ${\cal M}(n)$ consisting of
$(C, x_0, ..., x_n)$ such that $C$ is isomorphic to ${\bf C}P^1$. The space
${\cal M}(n)$ is a smooth complex projective variety of dimension
$n-2$. It has the following elementary construction, see  [BG 1] [FM]
[Ka 3].

\vskip .2cm

Let Aff $= \{ x\mapsto ax+b\}$ be the group of affine transformations
of {\bf C}. The group Aff acts diagonally on ${\bf C}^n$ preserving the
open part ${\bf C}^n_*$ formed by points with pairwise distinct coordinates.
As we noted before, we have an isomorphism $M_{0,n+1} \cong {\bf C}^n_*/{\rm
Aff}$.  Denote by $\Delta \i {\bf C}^n$ the principal diagonal, i.e.,
the space of points $(x,x,...,x)$. Then we have an embedding
$$M_{0,n+1} = {\bf C}^n_*/{\rm Aff} \i ({\bf C}^n - \Delta)/{\rm Aff} =
{\bf C}P^{n-2}.$$
The coordinate axes in ${\bf C}^n$ give $n$ distinguished points $p_1,
..., p_n \in {\bf C}P^{n-2}$. Let us blow  up all the points $p_i$,
then blow up the proper transforms
(= closures of the preimages of some open parts)
 of the lines $<p_i, p_j>$, then
blow up the proper transforms of the planes $<p_i, p_j, p_k>$, and so on.
It can be shown that the resulting space is isomorphic to
${\cal M}(n)$.

\vskip .3cm

\noindent {\bf (1.4.4) The configuration operad ${\cal M}$.}
The family of spaces ${\cal M} = \{{\cal M}(n), n\geq 1\}$, forms a 
topological operad. The symmetric group action on ${\cal M}(n)$ is
given by
$$(C, x_0, ..., x_n) \mapsto (C, x_0, x_{s(1)}, ..., x_{s(n)}), \quad s\in 
\Sigma_n.$$
The composition map
$${\cal M}(l) \times {\cal M}(m_1) \times ... \times {\cal M}(m_l)
\longrightarrow {\cal M}(m_1 + ... + m_l)$$
is defined by
$$ (C, y_0, ..., y_l), (C^{(1)}, x^{(1)}_0, ..., x^{(1)}_{m_1}), ..., 
(C^{(l)}, x^{(l)}_0, ..., x^{(l)}_{m_l}) \longmapsto$$
$$ \longmapsto (C', y_0, x{(1)}_1, ..., x^{(1)}_{m_1}, ..., 
x^{(l)}_1, ..., x^{(l)}_{m_l})$$
where $C'$ is the curve  obtained from the disjoint union 
$C \sqcup C^{(1)} \sqcup ... \sqcup C^{(l)}$ by identifying
$x_0^{(i)}$ with $y_i$, $i=1,...,l$ (Fig.5)

\vbox to 5cm{}

We   call ${\cal M}$ the {\it configuration operad}, since
${\cal M}(n)$ can be regarded as (compactified) configuration
spaces of points on ${\bf C}P^1$. 

\vskip .3cm

\noindent {\bf (1.4.5) The stratification of the space ${\cal M}(n)$.}
Given a point $(C, x_0, ..., x_n) \in {\cal M}(n)$, we associate to it
an $n$ - tree $T = T(C, x_0, ..., x_n)$ as follows. The vertices of $T$
correspond to the irreducible components of $C$. The vertices corresponding to
two components $C_1, C_2$ are joined by an (internal) edge if
$C_1\cap C_2 \neq \emptyset$. An external edge is assigned to each of the
marked points $x_0, ..., x_n$. The input edge labelled by $i \neq 0$ is
attached to the vertex corresponding to the component containing $x_i$.
The output edge is attached to the vertex corresponding to the component
containing $x_0$. The property (iii) of stable curves ensures that
$T(C, x_0, ..., x_n)$ is a reduced tree. 

\vskip .2cm

For any reduced $n$ - tree $T$ let ${\cal M}(T) \i {\cal M}(n)$
be the subset consisting of points $(C, x_0, ..., x_n)$ such that
$T(C, x_0, ..., x_n) = T$. In this way we obtain an algebraic stratification
${\cal M}(n) = \bigcup {\cal M}(T)$. This stratification has the following
properties (cf. [BG 1]):

$${\rm codim}\,\,\, {\cal M}(T) = \# {\rm \quad internal 
\quad edges \quad of
\quad} T.\leqno {\bf (1.4.6)}$$

$${\cal M}(T) \i \overline {{\cal M}(T')} \quad \Longleftrightarrow
\quad 
 T\geq T'.\leqno {\bf (1.4.7)}$$

In particular, 0 -  dimensional strata are labelled by binary trees.
Codimension 1 strata correspond to trees with two vertices. Their closures
are precisely the irreducible components of ${\cal M}(n) - M_{0,n+1}$
 which is a normal crossing divisor. Moreover, the entire stratification
 above can be recovered by intersecting these components in all
possible ways. In addition, we have the following result. 

\proclaim (1.4.8) Proposition. There are canonical direct product
decompositions
$${\cal M}(T) = \prod_{v\in T} M_{0, |{\rm In}(v)| + 1},$$
$$\overline {{\cal M}(T)} = \prod_{v\in T} {\cal M}({\rm In}(v)).$$
In particular, the closure of each stratum is smooth. 

\vskip .2cm

Let $T(m_1, ..., m_l)$ be the tree in Fig.3. 

\proclaim (1.4.9) Corollary. The structure maps (1.2.1)(ii) of the
operad ${\cal M}$ can be identified with the embedding of the stratum
$${\cal M}(l) \times {\cal M}(m_1) \times ... \times {\cal M}(m_l) 
= {\cal M}(T(m_1,...,m_l) \hookrightarrow {\cal M}(n).$$

\hfill\vfill\eject

\centerline {\bf 1.5. Operads and sheaves.}

\vskip 1cm

\noindent {\bf (1.5.1)} Let ${\cal P}$ be a $k$ - linear operad.
The compositions in ${\cal P}$ make it possible to construct, for
any $n$, a sheaf ${\cal F}_{\cal P}(n)$ on the  moduli space
${\cal M}(n)$, as we now proceed to explain. 

\vskip .2cm

Let $X$ be any CW - complex and $S = \{X_\alpha\}$ be a Whitney
stratification [GM] of $X$ into connected strata $X_\alpha$. We say that
a sheaf ${\cal F}$ of $k$ - vector spaces on $X$ is $S$ - 
{\it combinatorial}
if the restriction of ${\cal F}$ to each stratum $X_\alpha$ is a 
constant sheaf. Thus ``$S$ - combinatorial" is more restrictive that
``$S$ - constructible" [KS] [GM] which means that ${\cal F}|_{X_\alpha}$ are
only locally constant. 

It is well known that giving an $S$ - combinatorial sheaf ${\cal F}$ is
equivalent to giving the following linear algebra data:

\item{(i)} Vector spaces $F_\alpha = H^0(X_\alpha, {\cal F})$,
one for each stratum $X_\alpha$;

\item{(ii)} {\it Generalization maps} $g_{\alpha\beta}: F_\alpha
\rightarrow F_\beta$ defined whenever $X_\alpha \i \overline{X_\beta}$
and satisfying the transitivity condition:
$$g_{\alpha\gamma} = g_{\beta\gamma} \circ g_{\alpha\beta} \quad
{\rm if} \quad X_\alpha \i \overline{X_\beta} \i \overline{X_\gamma}.$$

\vskip .3cm

\noindent {\bf (1.5.2)} We now consider the moduli space ${\cal M}(n)$
defined in \S 1.4 together with the stratification ${\cal M}(n) = 
\bigcup _T {\cal M}(T)$ labelled by the set of $n$ - trees. 
A sheaf on ${\cal M}(n)$ combinatorial 
 with respect to this stratification
will be referred to as a
``combinatorial sheaf of ${\cal M}(n)$".

\vskip .2cm

Let ${\cal P}$ be a $k$ - linear operad. To any  $n$ - tree $T$ we have 
associated in (1.2.13) a vector space ${\cal P}(T)$ and to any
pair of $n$ - trees $T'\leq T$ we have associated a linear map
$\gamma_{T,T'}: {\cal P}(T) \rightarrow {\cal P}(T')$. Note that we
have an inclusion ${\cal M}(T) \i \overline {{\cal M}(T')}$
precisely when $T'\leq T$. Thus associating to a stratum ${\cal M}(T)$
the space ${\cal P}(T)$ and using the $\gamma_{T,T'}$ as generalization 
maps, we get a combinatorial sheaf on ${\cal M}(n)$ which we denote by
${\cal F}_{\cal P}(n)$. 

The sheaves ${\cal F}_{\cal P}(n)$ for different $n$ enjoy certain
compatibility with the operad structure on $\{{\cal M}(n)\}$. 
Those compatibility properties are formalized in the next subsection.
Before proceeding to do this, let us agree on the following
shorthand notation. If $X, Y$ are topological spaces,
 $\cal F$ is a sheaf on $X$ and $\cal G$ is a sheaf on
$Y$, then  the sheaf
$p_X^*{\cal F} \otimes p_Y^*{\cal G}$ on $X \times Y$, where
$p_X: X\times Y \rightarrow X, p_Y: X\times Y \rightarrow Y$ are
the natural projections, will be shortly denoted by ${\cal F}
\otimes {\cal G}$. Similarly for more spaces.  

\vskip .3cm

\noindent {\bf (1.5.3) Sheaves on an operad.} Let $Q$ be a topological
operad and
$$\gamma_{m_1,...,m_l}: Q(l) \times Q(m_1) \times ... \times Q(m_l)
\rightarrow Q(m_1 + ... + m_l)$$
be its structure maps. A sheaf on $Q$ is, by definition, a collection 
${\cal F} = \{{\cal F}(n)\}$ where ${\cal F}(n)$ is a sheaf of
$k$- vector spaces on $Q(n)$ together with the following data:

\item{(i)} A structure of  $\Sigma_n$ - equivariant sheaf on ${\cal F}(n)$.

\item{(ii)} For each $m_1, ..., m_l$, a homomorphism of sheaves on
$Q(l) \times Q(m_1) \times ... \times Q(m_l)$:
$$r_{m_1, ..., m_l}: \gamma_{m_1...m_l}^* {\cal F}(m_1 + ... + m_l)
\rightarrow {\cal F}(l) \otimes {\cal F}(m_1) \otimes ... \otimes
{\cal F}(m_l).$$

\item{(iii)} A homomorphism $\epsilon: {\cal F}(1)_1 \rightarrow k$ where
${\cal F}(1)_1$ is the stalk of the sheaf ${\cal F}(1)$ at the point $1\in
Q(1)$.

\vskip .2cm

These data should satisfy three conditions of associativity, equivariance
and counit which we now explain.

\vskip .3cm

\noindent {\bf (1.5.4) The coassociativity condition.} Observe that the
associativity condition for the topological operad $Q$ 
amounts to the commutativity, for all $l, m_1, ..., m_l$, 
$m_{ij}, i=1, ..., l, j=1, ..., \nu_i$, of the diagrams
$$\matrix{& Q(l) \times \prod Q(m_i) \times\prod_{i,j}Q(m_{ij}) &\buildrel
1\times\gamma\over\longrightarrow & Q(l) \times \prod_i Q\biggl(
\sum_j m_{ij}\biggl)&\cr
\gamma\times 1&\big\downarrow&&\big\downarrow&\gamma\cr
&Q\biggl(\sum_im_{ij}\biggl) \times \prod_{i,j}Q(m_{ij})
 &\buildrel \gamma\over
\longrightarrow &Q\biggl(\sum_{i,j}m_{ij}\biggl)&}.$$

Correspondingly,  the structure data (ii) for a sheaf ${\cal F}$ 
give rise to the
following diagram of sheaves on $Q(l) \times \prod_i Q(m_i) \times
\prod_{i,j} Q(m_{ij})$:
$$\overfullrule =0pt\matrix{& (\gamma\circ(1\times\gamma))^* 
{\cal F}\biggl(\sum_{i,j}m_{ij}
\biggl)& \buildrel r\over\longrightarrow & (1\times\gamma)^* \left(
{\cal F}(l) \otimes \bigotimes_i {\cal F}\biggl(\sum_j m_{ij}\biggl)
\right)&\cr
r&\big\downarrow&&\big\downarrow& 1\times r\cr
&(\gamma\times 1)^*\left({\cal F}\biggl(\sum m_i\biggl) \otimes
\bigotimes_{i,j} {\cal F}(m_{ij})\right)& \buildrel r\otimes 1 \over
\longrightarrow &{\cal F}(l) \otimes \bigotimes_i {\cal F}(m_i)
\otimes \bigotimes_{i,j} {\cal F}(m_{ij})&}.$$

It is required that all such diagrams  commute. 

\vskip .3cm

\noindent {\bf (1.5.5) The equivariance condition.} It is required, first
of all, that the map $r_{m_1, ..., m_l}$ commutes with the natural
action of $\Sigma_{m_1} \times ... \times \Sigma_{m_l}$ on
${\cal F}(m_1 + ... + m_l)$ and ${\cal F}(m_1) \otimes ... \otimes
{\cal F}(m_l)$. 

Secondly, for any permutation $s\in \Sigma_l$ we denote by $\tilde s \in
\Sigma_{m_1 + ... + m_l}$ the block permutation which moves segments of
lengths $m_1, ..., m_l$ according to $s$. Note that the equivariance
condition for the topological operad $Q$ implies the commutativity
of the diagram
$$\matrix{& Q(l) \times Q(m_1) \times ... \times Q(m_l) &
\buildrel \gamma\over\longrightarrow & Q(m_1 + ... + m_l)&\cr
{\rm Id} \times s& \big\downarrow&&\big\downarrow &\tilde s^*\cr
&Q(l) \times Q(m_{s(1)}) \times ... \times Q(m_{s(l)}) &
\buildrel \gamma\over\longrightarrow & Q(m_{s(1)} + ... + m_{s(l)})&}$$
It is now required that the following diagram of sheaves on 
$Q(l) \times Q(m_1) \times ... \times Q(m_l)$  commutes:

$$\overfullrule = 0pt \hskip -1cm
\matrix{ & \pmatrix{(\gamma\circ(\tilde s\times {\rm Id}))^* {\cal F}
(m_{s(1)} + ... + m_{s(l)})\cr
= (\tilde s^*\circ\gamma)^*  {\cal F}(m_{s(1)} + ... + m_{s(l)})}&
\buildrel r\over\longrightarrow& (s\times {\rm Id})^*\left[
{\cal F}(l) \otimes {\cal F}(m_{s(1)}) \otimes ... \otimes 
{\cal F}(m_{s(l)}\right]&\cr
{\rm eq.}&\big\downarrow&&\big\downarrow&{\rm eq.}\cr
&\gamma^*{\cal F}(m_1 + ... + m_l) &\buildrel r\over\longrightarrow&
{\cal F}(l) \otimes {\cal F}(m_1) \otimes ... \otimes {\cal F}(m_l)&}$$
where ``eq." denotes morphisms of equivariance.

\vskip .3cm

\noindent {\bf (1.5.6) The counit condition.}  Recall that the element $1\in Q(1)$ is such
that
the composition
$$Q(l) \buildrel j\over \rightarrow Q(l) \times Q(1) \times ... \times Q(1)
\buildrel \gamma\over\rightarrow Q(l)$$
where $j(q) = (q,1,...,1)$, is the identity map. It is required,
in addition, that the
composition
$${\cal F}(l) = (\gamma\circ j)^* {\cal F}(l)\buildrel r\over \rightarrow 
j^*\left({\cal F}(l)\otimes {\cal F}(1) \otimes ... 
\otimes {\cal F}(1)\right)\buildrel {\rm Id} \otimes\epsilon \otimes ...
\otimes \epsilon \over\longrightarrow {\cal F}(l)$$
is the identity homomorphism. 

This completes the definition of a sheaf on a topological operad. 

\vskip .3cm

\noindent {\bf (1.5.7) Example.} Let $Pt$ be the topological operad
having $Pt(n) = \{pt\}$ (one - point space) for any $n$ and all
the structure maps  being the identities . We   call $Pt$
the {\it trivial operad} (it defines associative and 
commutative $H$ - spaces, see [May]). Let ${\cal F}$ be a sheaf on $Pt$.
For any $n$ the  sheaf ${\cal F}$ on the space $Pt(n)$ is
just a vector space $F(n)$. The structure data (i) - (iii) of (1.5.3)
amount to $\Sigma_n$ - action on the space $F(n)$, linear maps
$$F(m_1 + ... + m_l) \rightarrow F(l) \otimes F(m_1) \otimes
... \otimes F(m_l)$$
and a linear functional $\epsilon: F(1) \rightarrow k$.

It is immediate to see that these data make the collection of the dual
vector spaces $F(n)^*$ into a $k$ - linear operad. Conversely,
 every $k$ - linear operad ${\cal P}$ with finite - dimensional spaces
${\cal P}(n)$ defines a sheaf on the operad $Pt$.

\vskip .3cm

\noindent {\bf (1.5.8) Remark.} One might say that the $F(n)$ in the above
example form a {\it cooperad}, a structure dual to that of an operad. We
prefer to postpone the discussion of this structure until Chapter 3. 

\vskip .2cm

The considerations of (1.5.7) can be generalized as follows.

\proclaim (1.5.9) Proposition. Let $Q$ be a topological operad and
${\cal F}$ be a sheaf on $Q$. Then the graded spaces $H^\bullet(
Q(n), {\cal F}(n))^*$ (the duals to the total cohomology spaces)
form an operad in the category $gVect^-$ of graded vector spaces.
Moreover, for any $q\geq 0$ the subspaces $H^{q(n-1)}(Q(n),
{\cal F}(n))^*$ form a sub - operad.

\noindent {\sl Proof:} Obvious from the axioms. Left to the reader. 

\vskip .2cm

The fact (1.4.1) that the homology spaces $H_\bullet(Q(n),k)$ form an 
operad is a particular case of this proposition. Indeed, for any
topological operad $ Q$ the constant sheaves $\underline{k}_{
{ Q}(n)}$
form a sheaf on $ Q$.

\vskip .3cm

\noindent {\bf (1.5.10)} Let $ Q$ be a topological operad and ${\cal F}$
be a sheaf on $ Q$. We   say that ${\cal F}$ is an {\it iso - 
sheaf} if all the maps $r_{m_1,...,.m_l}$ (data (ii) of a sheaf) are 
isomorphisms. Now we can formulate the precise relation
between $k$ - linear operads and sheaves on the configuration
operad ${\cal M}$.

\proclaim (1.5.11) Theorem. a) If ${\cal P}$ be a $k$ - linear
operad, then the sheaves ${\cal F}_{\cal P}(n)$ on the spaces
${\cal M}(n)$ introduced in (1.5.2) form  a sheaf ${\cal F}_{\cal P}$
on the operad ${\cal M}$.\hfill\break
b) If ${\cal P}(1) = k$ then ${\cal F}_P$ is an iso - sheaf.\hfill
\break
c) Any  combinatorial
iso - sheaf ${\cal F}$ on ${\cal M}$ 
has the form ${\cal F}_P$ for some $k$ - linear operad ${\cal P}$
with ${\cal P}(1) = k$.

The proof is straightforward and left to the reader.

\vskip .3cm

\noindent {\bf (1.5.12) Remark.} If one replaces in (1.5.11)(c) the
adjective ``combinatorial" by ``constructible" then one obtains a
 notion of a {\it braided operad} introduced by 
Fiedorowicz [F]. 

\vfill\eject

\centerline {\bf 1.6. Modules over an algebra over an operad.}

\vskip 1cm

\noindent {\bf (1.6.1)} Let ${\cal P}$ be a $k$ - linear operad and
$A$ be a ${\cal P}$ - algebra. An $A$ - {\it module} (or a 
$({\cal P}, A)$ - module, if ${\cal P}$ is  to be specified explicitly)
is a $k$ - vector space $M$ together with a collection of linear maps
$$q_n: {\cal P}(n) \otimes A^{\otimes (n-1)} \otimes M \longrightarrow
M, \quad n\geq 1$$
satisfying the following conditions:

\item{(i)} (Associativity)  For any natural numbers $n, r_1, ..., r_n$,
any elements $\lambda \in {\cal P}(n)$, $\mu_i \in {\cal P}(r_i)$,
$a_{ij}\in A$, $j=1,...,r_i$, $i=1,...,n-1$ and also $a_{n1}, ..., 
a_{n,r_n-1}\in A$, $m\in M$ we have the equality
$$q_{m_1+...+m_n}(\lambda(\mu_1, ..., \mu_n) \otimes a_{11} \otimes
... \otimes a_{1,r_1} \otimes a_{21} \otimes  ... ... \otimes a_{n, 
r_n-1} \otimes m) =$$
$$= q_n\biggl(\lambda \otimes \mu_1(a_{11}, ..., a_{1,r_1}) \otimes ...
\otimes \mu_{n-1}(a_{n-1,1} , ..., a_{n-1, r_{n-1}}) \otimes$$
$$\otimes q_{r_n}(\mu_n \otimes a_{n1} \otimes ... \otimes a_{n, r_n-1}
\otimes m)\biggl).$$

\item{(ii)} (Equivariance) The map $q_n$ is equivariant with respect
to the action of $\Sigma_{n-1} \i \Sigma_n$ on ${\cal P}(n) \otimes
A^{\otimes (n-1)}\otimes M$ given by
$$s(\lambda \otimes (a_1\otimes ... \otimes a_{n-1}) \otimes m) =
s^*(\lambda) \otimes (a_{s(1)} \otimes ... \otimes a_{s(n-1)})
\otimes m.$$

\item{(iii)} (Unit) We have $q_1({\bf 1} \otimes 1 \otimes m) = m$ for any
$m\in M$ where ${\bf 1} \in {\cal P}(1)$ is the unit of ${\cal P}$
and 1 is the unit element of $k = A^{\otimes 0}$.

\vskip .2cm

We   often write $\lambda(a_1, ..., a_{n-1}, m)$ for
$q_n(\lambda\otimes (a_1 \otimes ... \otimes a_{n-1}) \otimes m)$.

\vskip .3cm

\noindent {\bf (1.6.2) Examples.} a) Every ${\cal P}$ - algebra is
a module over itself.

\vskip .2cm

b) Let ${\cal P} = {\cal A}s$  and $A$ be a $\cal P$ - algebra, i.e.,
 an
associative algebra. An $({\cal A}s, A)$ - module is the same as an $A$ -
bimodule in the usual sense. Indeed, let us realize ${\cal A}s (n)$
as the space of non-commutative polynomials in $x_1,...,x_n$ spanned
by the monomials $x_{s(1)} \cdot ... \cdot x_{s(n)}$, $s\in \Sigma_n$, 
see
(1.3.6). Given an $A$ - bimodule $M$, we define the map 
$q_n: {\cal A}s (n)
 \otimes A^{\otimes (n-1)} \otimes M \rightarrow M$ by the rule
$$q_n( x_{s(1)} \cdot ... \cdot x_{s(n)} \otimes
(a_1\otimes ... \otimes a_{n-1}) \otimes m) = a_{s(1)}\cdot ... \cdot 
a_{s(n)}$$
where we set $a_n = m$. 

\vskip .2cm

c) Let ${\cal P} = {\cal C}om$ and $A$ be a $\cal P$ - algebra, i.e.,
a commutative algebra. A $({\cal C}om, A)$ - module is the same
as an $A$ - module in the usual sense. 

\vskip .2cm

d) Let ${\cal P} = {\cal L}ie$ and $A$ be a Lie algebra. A $({\cal L}ie, A)$ -
module is the same as an $A$ - module (representation of the Lie algebra
$A$) in the usual sense.

\vskip .3cm

\noindent {\bf (1.6.3)} Let $M, M'$ be two $({\cal P}, A)$ - modules.
A morphism of modules $f: M\rightarrow M'$ is a $k$ - linear map
such that $f(\lambda(a_1, ..., a_{n-1},m)) = \lambda(a_1, ..., a_{n-1},
f(m))$ for any $n$, any $\lambda \in  {\cal P}(n), a_i \in A, m\in M$.
Clearly, all $({\cal P}, A)$ - modules form an Abelian category. 

\vskip .3cm

\noindent {\bf (1.6.4) The universal enveloping algebra.}
As for every Abelian category, it is natural to expect that the
category of $({\cal P}, A)$ - modules can be described in therms of  left
modules over a certain associative algebra. Such an algebra indeed
 exists and is called the {\sl universal enveloping algebra} of A, to be
denoted $U(A)$ or $U_{\cal P}(A)$. Here is the construction. 

\vskip .2cm

By definition, $U(A)$ is generated by symbols $X(\lambda; a_{1},\ldots
, a_{n-1})$, for every $n$, every $\lambda\in {\cal P} (n)$ and every 
 $a_{1}, \ldots
, a_{n-1}\in A$. (In particular,
when $n=1$ then no $a_{i}$ should be specified and we have generators
$X(\lambda), \lambda\in {\cal P} (1)$). These symbols are
 subject to conditions of
 polylinearity with respect to each argument and to the following 
identifications:
$$X(\lambda ; \mu_{1}(a_{11},\ldots
,a_{1,r_{1}}), \ldots
, \mu_{n-1}(a_{n-1,1}, \ldots
, a_{n-1,r_{n-1}})) = \leqno {\bf (1.6.5)}$$
$$= X( \lambda (\mu_{1},\ldots
,\mu_{n-1}, {\bf 1}); a_{11},\ldots
, a_{n-1,r_{n-1}})$$
for any $\lambda\in {\cal P} (n), \mu_{i}\in {\cal P}(m_{i}), 
a_{ij}\in A$.

The multiplication in the algebra $U(A)$ is given by the formula:
\par
$$X (\lambda; a_{1},\ldots
, a_{n-1}) X(\mu, b_{1},\ldots
, b_{l-1}) =
X( \lambda ({\bf 1},\ldots
, {\bf 1},\mu); a_{1},\ldots
, a_{n-1}, b_{1},\ldots
, b_{l-1}).$$

\noindent It is immediate to verify that we get
 in this way  an associative
algebra $U(A)$ and the image of $X({\bf 1})$ is the unit of this algebra.

\vskip .2cm

Heuristically, $X(\lambda ; a_{1},\ldots
, a_{n-1})$ corresponds to the operator
$$m \mapsto \lambda (a_1, ..., a_{n-1}, m)$$
\noindent which acts on every $({\cal P}, A)$ - module M.

The following fact is obvious and its proof is left to the
reader.

\proclaim (1.6.6) Proposition. Let A be a ${\cal P}$-algebra. The category of
$A$ - modules is equivalent to the category of left modules over the
associative algebra $U_{\cal P}(A)$.

 \noindent {\bf  (1.6.7) Examples.} a) If ${\cal G}$ is a Lie algebra 
then $U_{{\cal L}ie}({\cal G})$ defined above is
the ordinary universal enveloping algebra of ${\cal G}$.

\vskip .2cm

b) If $A$ is an associative algebra and we regard 
it as an ${\cal A}s$-algebra
 then $U_{{\cal A}s}(A) = A \otimes 
 A^{op}$.

\vskip .2cm

 c) If $A$ is an accosiative commutative algebra and 
we regard it as
${\cal C}om$ - algebra  then
 $U_{{\cal C}om}(A) = A$.

\vskip .3cm

\noindent {\bf (1.6.8)} 
 By construction, the universal enveloping 
algebra $U_{\cal P}(A)$ has the form
$$U_{P}(A) = \bigoplus {\cal P} (n)\otimes A^{\otimes 
(n-1)}/\equiv $$ 
  where $\equiv$ is the equivalence relation
described in (1.6.5). The $n$ - th summand  above is just the space of
 generators $X(\lambda ; a_{1},\ldots
, a_{n-1})$.

Observe that $U(A) = U_{\cal P}(A)$ has a natural
(multiplicative) filtration $F_\bullet$ where
$F_{n}U(A)$ is the subspace consisting of images of generators
$X(\lambda ;a_{1},\ldots
,a_{r-1})$ for $r\le n$.

\hfill\vfill\eject

\centerline {\bf 2. QUADRATIC OPERADS.}

\vskip 1.5cm

\centerline {\bf 2.1. Description of operads by generators
and relations.}

\vskip 1cm

\noindent {\bf (2.1.1)} Let $K$ be an associative $k$ - algebra.
Given an arbitrary $K$ - collection $E = \{E(n), n\geq 2\}$
(1.2.11) we define an operad $F(E)$ called the {\it free operad}
generated by $E$. By definition
$$F(E)(n) = \bigoplus_{ n - {\rm trees}\,\,\, T} E(T)$$
where $T$ runs over isomorphism classes of $n$ - trees and
$E(T)$ was defined in (1.2.13). The composition maps (1.2.1)(ii)
for $F(E)$
$$F(E)(l) \otimes F(E)(m_1) \otimes ... \otimes F(E)(m_l) \rightarrow
F(E)(m_1 + ... + m_l)$$
are defined by means of maps
$$E(T) \otimes E(T_1) \otimes ... \otimes E(T_l) \rightarrow
E(T(T_1, ..., T_l))\leqno{\bf (2.1.2)}$$
where $T$ is an $l$ - tree, $T_i, i=1,...,l$, is a $m_i$ - tree and
$T(T_1, ..., T_l)$ is their composition, see Fig.6.

\vbox to 5cm{}

The  definition of the map (2.1.2) is obvious,  keeping
 in mind that both the
LHS and RHS are tensor products of the same spaces but over different
rings (some over $k$ and some over $K$).

\vskip .3cm

\noindent {\bf (2.1.3) Ideals.} 
 Let ${\cal P} = \{{\cal P} (n)\}$ be  a $k$-linear operad.
 A (two-sided) ideal in $\cal P$
is a collection $\cal I$ of  vector subspaces
 ${\cal I} (n)\subset {\cal P} (n)$ 
satisfying the following three conditions:

\item{(i)} For each $n$ the space ${\cal I} (n)$ is preserved by the action 
of $\Sigma_{n}$ on
${\cal P} (n)$.

\item{(ii)} 
If $\lambda \in {\cal P}(n), \mu_{1}\in {\cal P} (m_{1}),\ldots
, \mu_{n}\in {\cal P} (m_{n})$, and for at least one $j$ we have
 $\mu_{j}\in {\cal I} (m_{j})$, then the composition $\lambda
(\mu_{1},\ldots, \mu_{n})$ belongs to
${\cal I}(m_{1}+\ldots  +m_{n})$.

\item{(iii)}
If $\lambda\in {\cal I}(n)$ and $\mu_{i}\in {\cal P} (m_{i}), i=1,\ldots
,n$, then $\lambda (\mu_{1},\ldots
, \mu_{n})\in{\cal I}(m_{1} + \ldots + m_{n})$.

\vskip .2cm

If ${\cal I}$ is an ideal in an operad ${\cal P}$ then we can construct
 the quotient
operad ${\cal P}/{\cal I}$ with components $({\cal P}/{\cal I})(n) = 
{\cal P}(n)/{\cal I}(n)$ (quotient linear
spaces). The conditions (ii) and (iii) above imply that compositions in
$\cal P$ induce well defined compositions in ${\cal P}/{\cal I}$.

It is straightforward to see that the kernel of a morphism of $k$ -
linear operads $f: {\cal P} \rightarrow {\cal Q}$ is an ideal in ${\cal P}$. 

\vskip .3cm

\noindent {\bf (2.1.4)} Let $V$ be a finite-dimensional $k$ - vector space
and let ${\cal L}ie (V)$, ${\cal A}s (V)$, ${\cal C}om (V)$ be respectively
the free Lie algebra, free associative (tensor) algebra and free
commutative (polynomial) algebra generated by $V$. There are canonical
linear maps
$${\cal L}ie (V) \buildrel \epsilon \over\hookrightarrow {\cal A}s(V)
\buildrel \pi \over\longrightarrow {\cal C}om (V).\leqno {\bf (2.1.5)}$$
Both maps are the identity on $V$ are are uniquely determined by the
requirement that $\pi$ is a homomorphism of associative algebras
and $\epsilon$ is a homomorphism of Lie algebras (with the structure
of a Lie algebra on ${\cal A}s(V)$ given by $[a,b] = ab-ba$). The maps
(2.1.5) give rise to a collection of linear maps
$${\cal L}ie (n) \buildrel \epsilon_n \over\hookrightarrow {\cal A}s(n)
\buildrel \pi_n \over\longrightarrow {\cal C}om (n), \quad n=1,2,...$$
such that $\pi_n \circ \epsilon_n = 0$ for  $n\geq 2$. 

These maps give morphisms of operads ${\cal L}ie \buildrel \epsilon
 \over\hookrightarrow {\cal A}
\buildrel \pi \over\longrightarrow {\cal C}om$. Let ${\cal L}ie^+$ be
the collection of spaces $\{{\cal L}ie(n), n\geq 2\}$ and $({\cal L}ie^+)$
be the minimal ideal in ${\cal A}s$ containing ${\cal L}ie^+$. The
proof of the following result is left to the reader.

\proclaim (2.1.6) Proposition. $({\cal L}ie^+) = {\rm Ker} (\pi)$  that is,
${\cal C}om \cong {\cal A}s /({\cal L}ie^+)$. 

\vskip .2cm

\noindent {\bf (2.1.7) Quadratic operads.} Let $K$ be a semisimple $k$ 
- algebra. Let $E$ be a $(K, K^{\otimes 2})$ - bimodule with an involution
$\sigma: E \rightarrow E$ such that 
$$\sigma(\lambda e) = \lambda \sigma(e), \quad \sigma (e\cdot (\lambda_1
\otimes \lambda_2)) = \sigma(e) \cdot (\lambda_2\otimes\lambda_1),\quad
\forall \lambda, \lambda_1, \lambda_2 \in K, e\in E.$$
We form the space $E\otimes_K E$, the tensor product with respect to
the right $K$ - module structure on the first factor
 given by $e\cdot\lambda = e\cdot (\lambda
\otimes 1)$. This space has two structures:

\item{(i)} A $\Sigma_2$ - action given by the action of $\sigma$ on the
second factor $E$.

\item{(iii)} A structure of  $(K, K^{\otimes 3})$ - bimodule.

\vskip .2cm

Therefore the induced $\Sigma_2$ - module ${\rm Ind}_{\Sigma_2}^{\Sigma_3}
(E\otimes_K E)$ inherits the $(K, K^{\otimes 3})$ - bimodule structure.
Let $R \i {\rm Ind}_{\Sigma_2}^{\Sigma_3}
(E\otimes_K E)$ be a $\Sigma_3$ - stable $(K, K^{\otimes 3})$ - sub - bimodule.
To any such data $(E,R)$ we associate an operad ${\cal P}(K, E, R)$ in
the following way. 

\vskip .2cm

We form the $K$ - collection  $\{E(2) = E, E(n) = 0, n > 2\}$
(denoted also by $E$) and the corresponding free operad $F = F(E)$.
Observe that $F(E)(3) = {\rm Ind}_{\Sigma_2}^{\Sigma_3}
(E\otimes_K E)$, see Fig.7.

\vbox to 5cm{}

More generally,
$$F(E)(n) = \bigoplus_{ {\rm binary}\atop n- {\rm trees} \,\,\,T}
 E(T).\leqno
{\bf (2.1.8)}$$
Let $(R)$ be the ideal in $F(E)$ generated by the subspace $R\i F(E)(3)$.
We put $ {\cal P} = {\cal P}(K,E,R) = F(E)/(R)$. An operad of
 type ${\cal P}(K,E,R)$ is called
a {\it quadratic operad} with the space of generators $E$ and the space
of relations $R$. Note that $K, E, R$ can be recovered from $\cal P$ as
$K = {\cal P}(1)$, $E = {\cal P}(2)$, and $R = {\rm Ker} \{ F(E)(3)
\rightarrow {\cal P}(3)\}$. 

\vskip .3cm

\noindent {\bf (2.1.9) The quadratic duality.} Given a semisimple $k$ -
algebra $K$ and a finite - dimensional left $K$ - module $V$ with a
$\Sigma_n$ - action, we define $V^\vee = {\rm Hom}_K (V, K)$.
This is a right $K$ - module, i.e., a left module over the opposite
algebra $K^{op}$. We   always equip $V^\vee$ with the transposed action
of $\Sigma_n$ twisted by the sign representation. 

\vskip .2cm

Let ${\cal P} = {\cal P}(K, E, R)$ be a quadratic operad. The space
$E^\vee$ has a natural structure of $(K^{op}, K^{op} \otimes K^{op})$ -
bimodule. Observe that $F(E^\vee)(3) = F(E(3))^\vee$. Let $R^\bot \i
F(E^\vee)(3)$ be the orthogonal complement of $R$. It is stable under the
$\Sigma_3$ - action and the three $K^{op}$ - actions on $F(E^\vee)(3)$.
We define the {\it dual quadratic operad} ${\cal P}^!$ to be 
$${\cal P}^! = {\cal P}( K^{op}, E^\vee, R^\bot).$$

\vskip .3cm

\noindent {\bf (2.1.10) Examples.} Suppose  that $K = k$ and ${\cal P} = 
{\cal P}(k, E, R)$ is a quadratic operad. A ${\cal P}$ - algebra is
a vector space $A$ with several {\it binary} operations (parametrized
by $E$) which are subject to certain identities (parametrized
by $R$) each involving three arguments. This is precisely the way
of defining the types of algebras most commonly encountered in practice.
In particular, the operads ${\cal A}s, {\cal C}om, {\cal L}ie$
governing, respectively, associative, commutative and Lie algebras,
are quadratic.

\vskip .2cm

Note that the structure constants of an algebra $A$ over a quadratic 
operad satisfy quadratic relations. For example, if $\{e_i\}$ is a
basis of $A$, define $\mu: A \times A \rightarrow A$ by
$\mu(e_i\otimes e_j) = \sum c_{ij}^k e_k$. The condition that $\mu$
is associative means that 
$$\sum_k c_{ij}^k c_{kl}^m = \sum_k c_{jk}^k c_{ik}^m \quad \forall i,j, m.$$
Similarly in other examples. This explains the name ``quadratic operad."

\vskip .2cm

\proclaim (2.1.11) Theorem. We have isomorphisms of operads
$${\cal C}om^! = {\cal L}ie, \quad {\cal L}ie^! = {\cal C}om, \quad
{\cal A}s^! = {\cal A}s.$$

The idea that the ``worlds" of commutative, Lie and associative
algebras are, in some sense, dual to each other as described above,
was promoted by Drinfeld [Dr] and Kontsevich [Kon 2]
(the ${\cal C}om - {\cal L}ie$  duality was implicit already in Quillen's
paper 
[Q 1] and in Moore's Nice talk [Mo]).
The concept of
quadratic operads and their Koszul duality allows us to make this
idea into a theorem. Observe also that the isomorphism
${\cal C}om \cong {\cal A}s /({\cal L}ie^+)$ of Proposition
2.1.6 is, in a sense, ``self - dual".

\vskip .3cm

\noindent {\sl Proof of the theorem:} The group $\Sigma_3$ has three
irreducible representations which we denote {\bf 1} (the identity
representation), Sgn (the sign representation) and $V_2$ (the 2-dimensional
representation in the hyperplane $\sum x_i = 0$ in the 3-dimensional
 space of $(x_1, x_2, x_3)$ ).
Both ${\cal C}om$ and ${\cal L}ie$ have 1-dimensional spaces of generators
with $\Sigma_2$ acting trivially on ${\cal C}om(2)$ and by sign
on ${\cal L}ie(2)$. An elementary calculation of  group characters
shows that we have isomorphisms of $\Sigma_3$ - modules
$$F({\cal C}om (2))(3) = {\bf 1} \oplus V_2, \quad
F({\cal L}ie(2))(3) = {\rm Sgn} \oplus V_2.$$
This implies that we have
$${\cal C}om (3) = {\bf 1}, \quad R_{{\cal C}om} = V_2, 
\quad {\cal L}ie (3) = V_2, \quad R_{{\cal L}ie} = {\rm Sgn}.$$
The duality between ${\cal C}om$ and ${\cal L}ie$ follows. 

To prove that ${\cal A}s^! = {\cal A}s$, we consider the space
$F({\cal A}s(2))(3)$. This space has dimension 12 and is spanned by 
the 12 expressions of the form $x_{s(1)} (x_(2)x_{s(3)}), 
x_{s(1)}(x_{s(2)}x_{s(3)}), s\in \Sigma_3$.
 
We introduce on this space a scalar product $<\,\, , \,\, >$
by setting all these
products orthogonal to each other and putting
$$<x_i(x_jx_k), x_i(x_jx_k)> = {\rm sgn}\pmatrix{1&2&3\cr i&j&k}, \quad
<(x_ix_j)x_k, (x_ix_j)x_k> = - {\rm sgn}\pmatrix{1&2&3\cr i&j&k}.$$
This product is sign - invariant with respect to the $\Sigma_3$ - action,
i.e., $<s(\mu), s(\nu)> = {\rm sgn} (s) <\mu, \nu>$. The (6-dimensional)
space of
relations $R_{{\cal A}s}$ is spanned by all the associators
$(x_i(x_jx_k) - (x_ix_j)x_k$. This space coincides with its own
annihilator
with respect to the described scalar product. This shows that
${\cal A}s^! = {\cal A}s$.

\hfill\vfill\eject

\centerline {\bf 2.2. The analogs of Manin's tensor products. The Lie operad
as a dualizing object.}

\vskip 1cm

\noindent {\bf (2.2.1)} In this section we   consider only quadratic 
operads ${\cal P}$ with ${\cal P}(1) = k$. Such an operad 
 is defined by a vector space 
of generators $E = {\cal P}(2)$ with an
involution $\sigma$ (action of $\Sigma_2$) and an $\Sigma_3$ -
invariant subspace
$R$ of relations inside $F(E)(3)$ where $F(E)$ is the free operad generated
by $E$. Observe that $F(E)(3)$ is the direct sum of three copies of
$E\otimes E$, see Fig.6. So we   refer to this space as $3(E\otimes E)$.

Our aim is to present an operad - theoretic version of [Man].

\vskip .3cm

\noindent {\bf (2.2.2)} Let $({\cal A}, \otimes)$ be any symmetric monoidal
category and ${\cal P}, {\cal Q}$ be two operads in ${\cal A}$.
The collection of objects ${\cal P}(n) \otimes {\cal Q}(n)$ forms
a new operad in ${\cal A}$ denoted by ${\cal P} \otimes {\cal Q}$.
The main property of this operad is that if $A$ is a ${\cal P}$ -
algebra and $B$ is a ${\cal Q}$ - algebra (in ${\cal A})$ then
$A\otimes B$ is a ${\cal P} \otimes {\cal Q}$ - algebra. 

We are interested in the case of $k$ - linear operads, i.e., ${\cal A} = 
Vect$.

\vskip .3cm

\noindent {\bf (2.2.3)}  Suppose that ${\cal P, Q}$
are quadratic  operads with ${\cal P}(1) = {\cal Q}(1) = k$
so ${\cal P}(2)$ and ${\cal Q}(2)$ are their spaces of generators.
We denote by ${\cal P}\circ {\cal Q}$ the sub - operad in
${\cal P} \otimes {\cal Q}$ generated by ${\cal P}(2) \otimes
{\cal Q}(2)$. 
 Let $R_{\cal P}\i 3({\cal P}(2)\otimes {\cal P}(2))$ and 
$R_{\cal Q}\i 3({\cal Q}(2)\otimes {\cal Q}(2))$
be the spaces of relations of $\cal P$ and $\cal Q$. 
Then the operad
${\cal P}\circ {\cal Q}$ is  described as follows.

\vskip .2cm

 The
space of generators of ${\cal P} \circ {\cal Q}$ is, by definition, 
${\cal P}(2)\otimes {\cal Q}(2)$. The third
component of the free operad generated by ${\cal P}(2)\otimes {\cal Q}(2)$ is 
$3\bigl ({\cal P}(2)\otimes {\cal Q}(2) \otimes {\cal P}(2) \otimes
 {\cal Q}(2)\bigl)$ 
which  can be regarded in two
ways as: 

\item{1.}  $({\cal P}(2)\otimes {\cal P}(2)) \otimes 3({\cal Q}(2)
\otimes {\cal Q}(2))$ and regarded as
such, it contains the subspace $({\cal P}(2)\otimes {\cal P}(2)) \otimes
 R_{\cal Q}$.

\item{2.}  $3({\cal P}(2) \otimes {\cal P}(2))\otimes ({\cal Q}(2)
\otimes {\cal Q}(2))$ and regarded 
as such, it contains the subspace
$R_{\cal P}\otimes ({\cal Q}(2)\otimes {\cal Q}(2))$.

\vskip .2cm

\proclaim (2.2.4) Proposition. If ${\cal P,Q}$ are quadratic then
${\cal P} \circ {\cal Q}$ is also quadratic with the space of generators
${\cal P}(2)\otimes {\cal Q}(2)$ and the space of relations
$$\bigl(({\cal P}(2)\otimes {\cal P}(2)) \otimes
 R_{\cal Q}\bigl) \quad +  \quad \bigl(R_{\cal P}\otimes 
({\cal Q}(2)\otimes {\cal Q}(2))\bigl).$$
(The sum is not necessarily direct.)

The proof is straightforward.

\vskip .2cm

\noindent Thus ${\cal P} \circ {\cal Q}$ is the analog of the white circle
product $A\circ B$ for associative algebras considered by Manin [Man].

\vskip .3cm

\noindent {\bf (2.2.5)} Given two quadratic operads ${\cal P, Q}$ as before
we define the quadratic operad ${\cal P} \bullet {\cal Q}$ (the black
circle product, cf. [Man]) to have the same space of generators
${\cal P}(2) \circ {\cal Q}(2)$ as ${\cal P} \otimes {\cal Q}$ but the
space of relations
$$\bigl(({\cal P}(2)\otimes {\cal P}(2)) \otimes
 R_{\cal Q}\bigl) \quad \cap \quad \bigl(R_{\cal P}\otimes 
({\cal Q}(2)\otimes {\cal Q}(2))\bigl ).$$

\vskip .2cm

\proclaim (2.2.6)  Theorem.  Each of the products 
$\circ, \bullet$ defines on the
category of quadratic operads (and morphisms defined in (1.3.1))
a symmetric monoidal structure. Moreover, we have:
\item{(a)} $({\cal P} \circ {\cal Q})^! = {\cal P}^! \bullet
{\cal Q}^!$.
\item{(b)} $Hom (P\bullet Q, R) = Hom (P, Q^!\circ R)$.
In particular,  the commutative operad is
a unit object with respect to $\circ$ and the Lie operad is
a unit object with respect to $\bullet$.

\noindent {\sl Proof:} a) It is similar to the argument of Manin
[Man] for quadratic algebras by using the obvious relations
between sums, intersections and orthogonal complements of subspaces
in a vector space.

\vskip .2cm

b) The fact that ${\cal C}om$ is a unit object with respect to $\circ$
follows because ${\cal C}om (n) = k$ for any $n$. The fact that
${\cal L}ie$ is a unit object  with respect to $\bullet$, follows from
part a) and the fact that ${\cal L}ie^! = {\cal C}om$.

\vskip .3cm

\noindent {\bf (2.2.7)} Given quadratic operads ${\cal P}$ and ${\cal Q}$,
we define the quadratic operad $hom ({\cal P}, {\cal Q})$ as follows.
Its space of generators is set to be
$$hom ({\cal P}, {\cal Q})(2) = {\rm Hom}_k ({\cal P}(2), {\cal Q}(2)).$$
The space of relations $R_{hom ({\cal P}, {\cal Q})}$ is defined to be
the minimal subspace in 
\hfill\break $F({\rm Hom}_k ({\cal P}(2), {\cal Q}(2)))(3)$
such that the canonical map
$${\rm can}: {\cal Q}(2) \rightarrow {\rm Hom}_k ({\cal P}(2), {\cal Q}(2))
\otimes {\cal P}(2)) = {\cal P}(2)^* \otimes {\cal Q}(2) \otimes
{\cal P}(2)$$
extends to a morphism of operads ${\cal Q} \mapsto 
hom ({\cal P}, {\cal Q}) \circ {\cal P}$. More precisely, we define 
$R_{hom ({\cal P}, {\cal Q})}$ to be minimal among subspaces $J$
with the property that
$$J \otimes ({\cal P}(2) \otimes {\cal P}(2)) \i 
F({\rm Hom}_k ({\cal P}(2), {\cal Q}(2)) \otimes {\cal P}(2))(3)$$
contains the image of the embedding
$${\rm can}_*: F({\cal Q}(2))(3) \hookrightarrow
F({\rm Hom}_k ({\cal P}(2), {\cal Q}(2)) \otimes {\cal P}(2))(3)$$
induced by can. 

\vskip .2cm

If $\mu_1, ..., \mu_m$ is a basis of ${\cal P}(2)$, $\nu_1, ..., \nu_n$
is a basis of ${\cal Q}(2)$ and $a_{ij}$, $i=1,...,m$, $j=1, ..., n$ is
the corresponding basis of ${\rm Hom} ({\cal P}(2), {\cal Q}(2))$
then the elements
$$\tilde \nu_i = \sum_j a_{ij} \otimes \mu_j \in \left( hom({\cal P}, {\cal Q})
\circ {\cal P}\right)(2)$$
satisfy all the quadratic relations holding for actual $\nu_i$. 

\proclaim (2.2.8) Theorem. We have a natural isomorphism of operads
$$hom ({\cal P}, {\cal Q}) \cong {\cal P}^! \bullet {\cal Q}.$$

\noindent {\sl Proof:} Simple linear algebra.

\proclaim (2.2.9) Corollary. a) For any quadratic operad ${\cal P}$
we have
$${\cal P}^! = hom ({\cal P}, {\cal L}ie).$$
b) 
There exists a morphism of operads ${\cal L}ie \rightarrow
{\cal P} \otimes {\cal P}^!$ which takes the generator of the
1-dimensional space ${\cal L}ie (2)$ into the identity operator
in ${\cal P}(2) \otimes {\cal P}^!(2) = {\cal P}(2) \otimes {\cal P}(2)^*$.
In particular, for any $\cal P$ - algebra $A$ and ${\cal P}^!$
- algebra $B$ the vector space $A\otimes_k B$ has a natural Lie
algebra structure. 

\vskip .2cm

\noindent {\bf (2.2.10)} Let us give another interpretation of the
operad $hom({\cal P}, {\cal Q})$ in terms of algebras.  

Suppose that $V_1, ..., V_m, W$ are $k$ - vector spaces. By a
${\cal P}$ - {\it multilinear map} $V_1 \times ... \times V_m 
\rightarrow W$ we   understand an element
$$\Phi \in F_{\cal P}(V_1^* \oplus ... \oplus V_m^*) \otimes W$$
which is homogeneous of degree 1 with respect to dilatations of any of
$V_i$. (Recall (1.3.6) that $F_{\cal P}$ means the free $\cal P$ - algebra
generated by a vector space). Such maps form a vector space which
we denote by ${\rm Mult}_{\cal P}(V_1, ..., V_m| W)$. For example
${\cal P}(n) = {\rm Mult}_{\cal P}(k, ..., k|k)$ ($m$ copies of $k$
before the bar), cf. (1.3.7) - (1.3.9). The space of ordinary 
multilinear maps, i.e., ${\rm Hom}_k (V_1 \otimes ... \otimes V_m, W)$
is nothing but ${\rm Mult}_{{\cal C}om}(V_1, ..., V_m| W)$ where
${\cal C}om$ is the commutative operad. 

\vskip .2cm

Similarly to usual multilinear maps, ${\cal P}$ - multilinear maps can
be composed. In particular, for any vector space $V$ the spaces
$${\cal E}_{{\cal P}, V}(n) = {\rm Mult}_{\cal P}(V, ..., V| V)
\quad \quad (n \quad {\rm copies\quad of \quad} V)\leqno
{\bf (2.2.11)}$$ form a $k$ - linear operad which we call the 
{\it operad of endomorphisms of $V$ in ${\cal P}$}. 

\vskip .3cm

\noindent {\bf (2.2.12)} Let $\cal Q$ be another $k$ -
 linear operad. By a $\cal Q$ -
algebra in ${\cal P}$ we mean a $k$ - vector space $A$ together with
a morphism of operads $f: {\cal Q} \rightarrow {\cal E}_{{\cal P}, A}$.
When ${\cal P} = {\cal C}om$ is the commutative operad, we get the usual
notion of a ${\cal Q}$ - algebra. 

\proclaim (2.2.13) Theorem. Let ${\cal P}$ and ${\cal Q}$ be quadratic
operads with ${\cal P}(1) = {\cal Q}(1) = k$. Then ${\cal Q}$ - 
algebras in ${\cal P}$ are the same as $hom({\cal P}, {\cal Q})$ - 
algebras (in the usual sense).

\noindent {\sl Proof:} Let $A$ be a ${\cal Q}$ - algebra in ${\cal P}$.
The corresponding morphism $f: {\cal Q} \rightarrow {\cal E}_{{\cal P},
A}$ is defined by its second component
$$f_2: {\cal Q}(2) \rightarrow {\cal E}_{{\cal P}, A}(2) =
\left({\cal P}(2) \otimes (A^*)^{\otimes 2}\right)_{\Sigma_2} \otimes A.$$
By taking a partial transpose of $f_2$ we get a $\Sigma_2$ - 
equivariant map
$$f_2^\dagger : {\rm Hom}_k ({\cal P}(2), {\cal Q}(2)) \longrightarrow
{\rm Hom}_k (A \otimes A, A) = {\cal E}_A(2).$$
In order that a given linear map $f_2$ come from a morphism of operads
$f: Q \rightarrow {\cal E}_{{\cal P}, A}$, the quadratic relations among
the generators of $Q$ should be satisfied. By definition of relations
in $hom({\cal P}, {\cal Q})$ this is equivalent to the condition that
$f_2^\dagger$ extends to a morphism $hom({\cal P}, {\cal Q}) \rightarrow
{\cal E}_A$, i.e., that we have on $A$ a structure of a 
$hom({\cal P}, {\cal Q})$
- algebra. 

\vskip .3cm

\noindent {\bf (2.2.14) Lazard - Lie theory for formal groups in operads
and Koszul duality.} The idea of considering formal groups in operads 
goes back to an important paper of Lazard [L]. He used the concept
of ``analyseur" which is essentially equivalent to the modern notion
of operad. The main result of Lazard is a version of Lie theory
(= correspondence between (formal) Lie groups and Lie algebras) for
his generalized formal groups. It turns out that Lazard - Lie theory
has a very transparent interpretation in terms of Koszul duality for
operads. 

\vskip .2cm

We start with formulating basic definitions in the modern language. Let
${\cal P}$ be a $k$ - linear operad. We assume  ${\cal P}(1) = k$.
Let $W$ be a $k$ - vector space and
$$\hat F_{\cal P}(W) = \prod_{n\geq 1} \left( {\cal P}(n) 
\otimes W^{\otimes n}\right) _{\Sigma_n}, \leqno {\bf (2.2.15)}$$
the completed free ${\cal P}$ - algebra on $W$. If $z_1, ..., z_r$ form
a basis of $W$ then elements of $\hat F_{\cal P}(W)$ can be regarded
as formal series in $z_1, ..., z_r$ whose terms are products of
$z_i$ with respect to operations in ${\cal P}$. Thus, for example,
if ${\cal P} = {\cal C}om$, we get the usual power series algebra,
if ${\cal P} = {\cal A}s$, we get the algebra of non - commutative
power series etc. 
Note that our series do not have constant terms since all the free
algebras are without unit. 

There is a natural projection (on the first factor)
$$\hat F_{\cal P}(W) \rightarrow W. \leqno {\bf (2.2.16)}$$
It can be thought of as the differential at zero. 

\vskip .2cm

Let $V, W$ be two $k$ - vector spaces. We define the space of
{\it formal ${\cal P}$ - maps} from $V$ to $W$ as
$${\rm FHom}_{\cal P}(V, W) = V \otimes \hat F_{\cal P}(W^*).$$
For $\Phi \in \hat F_{\cal P}(W)$ we denote by $d_0\Phi \in
{\rm Hom}_k (V, W)$ its differential at 0, i.e., the image of $\Phi$
under the natural projection fo $V\otimes W^*$ induced by (2.2.16).
There are obvious composition maps
$${\rm FHom}_{\cal P}(V, W) \times {\rm FHom}_{\cal P}(W, X) 
\longrightarrow {\rm FHom}_{\cal P}(V, X)$$
(given by inserting  of power series into arguments of other
power series)
which make the collection of vector spaces and formal ${\cal P}$ -
maps into a category. 

A simple generalization of the classical inverse function theorem 
shows that $\Phi \in {\rm FHom}_{\cal P}(V, W)$ is invertible if
and only if $d_0\Phi \in {\rm Hom}_k (V,W)$ is invertible.

\vskip .2cm

\proclaim (2.2.17) Definition. Let $\cal P$ be a $k$ - linear operad. A
$\cal P$ - formal group is a pair $(V, \Phi)$ where $V$ is a $k$ -
vector space and $\Phi \in {\rm FHom}_{\cal P}(V \oplus V, V)$
is a formal ${\cal P}$ - map satisfying the two conditions:
{\item{(i)} $d_0\Phi: V \otimes V \rightarrow V$ is the addition map
$(v, v') \mapsto v + v'$.
\item{(ii)} $\Phi$ is associative i.e., we have the equality
$\Phi(x, \Phi(y,z)) = \Phi(\Phi(x,y), z)$ of formal maps
$V\oplus V \oplus V \rightarrow V$.
\hfill\break}
\noindent A formal homomorphism of formal ${\cal P}$ - groups
$(V,\Phi)$ and $(W, \Psi)$ is a formal ${\cal P}$ - map
$f\in {\rm FHom}_{\cal P}(V,W)$ such that $f(\Phi(x,y)) = 
\Psi( f(x), f(y))$. 

\vskip .2cm

\noindent {\bf (2.2.18)} If $(V, \Phi)$ is a ${\cal P}$ - formal group,
then we define, following Lazard, its Lie bracket by
$$[x,y] = \Phi(x,y) - \Phi(y,x) \quad ({\rm mod. \quad cubic\quad terms}).$$
This construction makes $V$ into a Lie algebra in ${\cal P}$ in the
sense of (2.2.12). 
Furthermore, we obtain, from Theorem 2.2.13 and the results of Lazard [L],
  the
following theorem.

\proclaim (2.2.19) Theorem. Let $\cal P$ be a quadratic operad and
${\cal P}^!$ be its quadratic dual. The category of finite - dimensional
${\cal P}$ - formal groups (and their formal homomorphisms) is equivalent
to the category of finite - dimensional ${\cal P}^!$ - algebras. 

A special case of this theorem corresponding to ${\cal P} = {\cal A}s$,
the associative operad, was pointed out to us earlier by M. Kontsevich.
Since ${\cal A}s^! = {\cal A}s$, the theorem in this case says that
for
formal groups defined by means of power series with non - commuting
variables the role of Lie algebras is played by associative algebras
(possibly without unit). 

\hfill\vfill\eject

\centerline {\bf 2.3. Quadratic algebras over a quadratic operad.}

\vskip 1cm

\noindent {\bf (2.3.1)} Let ${\cal P}$ be a $k$ - linear operad and
$A$  a ${\cal P}$ - algebra. An {\it ideal} in $A$ is a linear subspace
$I \i A$ such that for any $n$, any $\mu\in {\cal P}(n)$ and any
$a_1, ..., a_{n-1} \in A, i\in I$, we have $\mu (a_1, ..., a_{n-1}, i)\in I$. 
Given any ideal $I\i A$ the quotient vector space $A/I$ has a natural
structure of a ${\cal P}$ - algebra.

\vskip .3cm

\noindent {\bf (2.3.2) Quadratic algebras.} Let ${\cal P} =
{\cal P}(K, E, R)$ be a quadratic operad (2.1.7) so that $K = {\cal P}(1)$
is a semisimple $k$ - algebra and $E = {\cal P}(2)$ is a $(K, K^{\otimes 2})$
- bimodule. Let also $V$ be a $K$ - bimodule. The tensor product
$E \otimes_{K^{\otimes 2}}V^{\otimes 2}$ has a natural structure of a
$(K, K^{\otimes 2})$ - bimodule (the left $K$ - action comes from that
on $E$ and the right $K^{\otimes 2}$ - action comes from that on 
$V^{\otimes 2}$). Moreover, there is a $\Sigma_2$ - action on
$E \otimes_{K^{\otimes 2}}V^{\otimes 2}$ given by
$$\sigma(e \otimes (v_1 \otimes v_2)) = \sigma(e) \otimes (v_2\otimes v_1).
\leqno {\bf (2.3.3)}$$
Observe that the space of coinvariants $\left( E\otimes_{K^{\otimes 2}}
V^{\otimes 2}\right)_{\Sigma_2}$ inherits a $K$ - bimodule structure. Let
$$S \i \left( E\otimes_{K^{\otimes 2}}
V^{\otimes 2}\right)_{\Sigma_2} \leqno {\bf (2.3.4)}$$
be a $K$ - sub -  bimodule. Given $V$ and $S$ we   construct a ${\cal P}$
- algebra $A = A(V,S)$ as follows.

Let $F_{\cal P}(V)$ be the free ${\cal P}$ - algebra generated by $V$,
see (1.3.5). Observe that $F_{\cal P}(V)_2$, the degree 2 graded component
of $F_{\cal P}(V)$, is $\left( E\otimes_{K^{\otimes 2}}
V^{\otimes 2}\right)_{\Sigma_2}$. Let $(S) \i F_{\cal P}(V)$ be the
ideal generated by $S$ (= the minimal ideal containing $S$). Put
$$A(V, S) = F_{\cal P}(V)/(S).\leqno {\bf (2.3.5)}$$
An algebra of this type will be called a {\it quadratic ${\cal P}$ - 
algebra} (with the space of generators $V$ and the space of relations
$S$). 

\vskip .2cm

Observe that $A(V,S)$ has a natural grading $A(V,S) = \bigoplus_{i\geq 1}
A_i (V, S)$. Furthermore, ${\cal P}$, as any operad in the category
$Vect$, can be regarded as an operad in the category $gVect^+$ of graded
vector spaces (1.4.2)(c) and $A(V,S)$ is a ${\cal P}$ - 
algebra in this category.

\vskip .3cm

\noindent {\bf (2.3.6) Quadratic superalgebras.} Let ${\cal P} =
{\cal P}(K, E, R)$ be a quadratic operad as before. We now view ${\cal P}$
as an operad in the other category of graded vector spaces (1.4.2) namely
$gVect^-$. A ${\cal P}$ - algebra in this category will be referred to 
as a (graded) ${\cal P}$ - {\it superalgebra}. 

\vskip .2cm

Let $V$ be a $K$ - bimodule. We replace the $\Sigma_2$ - action on
$E\otimes_{K^{\otimes 2}}V^{\otimes 2}$ given by (2.3.3) by the
following one:
$$\sigma( e\otimes (v_1\otimes v_2)) = - \sigma(e) \otimes (v_2\otimes v_1).
\leqno {\bf (2.3.7)}$$
Let $S \i \left( E\otimes_{K^{\otimes 2}}
V^{\otimes 2}\right)_{\Sigma_2}$ be a $K$ - sub   - bimodule 
where the coinvariants are taken with respect to the new action
(2.3.7). Given such data $(V,S)$, we define a ${\cal P}$ - 
superalgebra $A(V,S)^{-}$ as the quotient as the free
${\cal P}$ - superalgebra generated by $V$ (placed in degree 1)
by the ideal generated by $S$.

\vskip .3cm

\noindent {\bf (2.3.8) Quadratic duality.} Let ${\cal P} = 
{\cal P}(K,  E, R)$ be a quadratic operad and ${\cal P}^! = {\cal P}
(K^{op}, E^\vee, R^\bot)$  the dual operad (2.1.9). Given a quadratic
${\cal P}$ - algebra $A = A(V,S)$ we define the quadratic ${\cal P}^!$ -
superalgebra $A^! = A(V^\vee, S^\bot)^{-}$. Here $V^\vee = 
{\rm Hom}_K (V, K)$ and
$$S^\bot \i \left( \left( E \otimes_{K^{\otimes 2}}
 V^{\otimes 2}\right)_{\Sigma_2}\right)^\vee = \left( (V^\vee)^{\otimes 2}
\otimes_{K^{\otimes 2}} E^\vee\right)^{\Sigma_2} = 
\left( E^\vee \otimes_{(K^{op})^{\otimes 2}} (V^\vee)^{\otimes 2}\right)
_{\Sigma_2}$$ 
is the annihilator of $S$. 

In a similar way, given a quadratic $\cal P$ - superalgebra $B$, we define the
dual ${\cal P}^!$ - algebra (in the category $gVect^+$) $B^!$. The
assignment $A \mapsto A^!$ gives a 1-1 correspondence between quadratic
$\cal P$ - algebras (resp. superalgebras) and quadratic
${\cal P}^!$ - superalgebras (resp. algebras).

\vskip .3cm

\noindent {\bf (2.3.9) Examples.} a) Let ${\cal P} = {\cal A}s$ be the
associative
operad. Then ${\cal A}s^! = {\cal A}s$. Note that an ${\cal A}s$ - super
algebra (i.e., an ${\cal A}s$ - algebra in $gVect^-$) is the same as an
${\cal A}s$ - algebra (in $gVect^+$). Both concepts give the usual
notion of a graded associative algebra. In this case the quadratic duality
(2.3.8) reduces to the well known Koszul duality for quadratic associative
algebras introduced by Priddy [Pr]. 

\vskip .2cm

b) Let ${\cal P} = {\cal C}om$ be the commutative algebra and $A$ be a 
quadratic ${\cal C}om$ - algebra i.e., a quadratic commutative
(associative) algebra. Let $A^!_{as}$ be the Koszul dual of $A$ as of
an associative algebra. One can show, see [Q 2] 
that $A^!_{as}$ has a structure of
a graded - commutative Hopf algebra. Hence $A^!_{as}$ is the
enveloping algebra of a certain graded Lie superalgebra ( i.e., 
a ${\cal L}ie$ - algebra in the category $gVect^-$) which we
denote ${\cal G} = 
{\rm Prim} (A^!_{as})$. We have  $A^! = {\cal G}$. 

\vskip .3cm

\noindent {\bf (2.3.10) Quadratic duality and enveloping algebras.}
Let ${\cal P} = {\cal P}(K, E, R)$ be a quadratic operad. Let $A$
be a quadratic $\cal P$ - algebra. Then  we have (1.6.4)
the universal enveloping
algebra $U_{\cal P}(A)$, which is an associative algebra
in ordinary sende.  If $A$ is a quadratic $\cal P$ -
superalgebra, the same construction as
in (1.6.4) define its universal enveloping superalgebra 
$U_{\cal P}^-(A)$ which is an ${\cal A}s$ - algebra in the
category $gVect^-$. As noted in (2.3.9)(a), 
we can regard $U^-_{\cal P}(A)$ as an ordinary graded associative
algebra. 

\proclaim (2.3.11) Theorem. If $A$ is a quadratic ${\cal P}$
- algebra 
and $A^!$ the dual quadratic ${\cal P}^!$ - superalgebra then the
universal enveloping algebras $U_{\cal P}(A), U_{{\cal P}^!}^-(A^!)$
are quadratic associative algebras in the ordinary sense and
$$\left(U_{\cal P}(A)\right)^! = U_{{\cal P}^!}^-(A^!).$$

\noindent {\sl Proof:}
Let ${\cal P} = {\cal P}(K, E, R)$ and $A = A(V, S)$.
The algebra $U_{\cal P}(A)$ has an obvious grading in which
the generator $X(\lambda; a_1, ..., a_{n-1})$, $\lambda \in {\cal P}(n),
a_i \in A_{n_i}$, has degree $i_1 + ... + i_{n-1}$. The degree 1
component of $U_{\cal P}(A)$ is linearly spanned by $X(\lambda, a)$,
$\lambda \in {\cal P}(2) = E, a\in A_1 = V$ and is isomorphic to
$E\otimes_K V$. Obviously this component generates $U_{\cal P}(A)$
as an algebra. 

To describe the relations among these generators, consider the space
$Y = F(E)(3) \otimes V^{\otimes 2}$ (here and in the remainder of this
section all tensor products are taken over $K$). This space splits
into the direct sum of three components $Y = L_1 \oplus L_2 \oplus L_3$
depicted in Fig. 8.

\vbox to 5cm{}

The letters ``$K$'' on some of the edges mean that the tensor product over
$K$ is taken with respect to the structures of left/right $K$ - 
module represented by the ends of
this edge. Note that every $L_i$ and hence $X$ is
a $K$ - bimodule with respect to the actions corresponding to the edges
not marked ``$K$''.  Elements of $L_1$ can be viewed as formal
expressions $\sum X(\lambda_i; a_i) X(\mu_i; b_i)$, as can elements
of $L_3$. The group $\Sigma_2$ permuting the left two inputs of the trees in
Fig. 8, maps $L_1$ isomorphically to $L_3$ and preserves $L_2$. Thus,
denoting $X = Y_{\Sigma_2}$ the space of coinvariants, we have
$$X \quad  = \quad  E\otimes V \otimes E\otimes V \quad\oplus\quad
E\otimes \bigl(E\otimes V^{\otimes 2}\bigl)_{\Sigma_2}.$$
Let $W\i X$ be the image, under the canonical projection $Y \rightarrow X$,
of the subspace $R\otimes V^{\otimes 2}$, where $R \i F(E)(3)$
is the space of relations of $\cal P$.

\proclaim (2.3.12) Proposition. The algebra $U_{\cal P}(A)$ is defined
by the space of generators $E\otimes V$ and the space of
quadratic relations
$$E\otimes V \otimes E\otimes V \quad\cap\quad (W \,\, + \,\,
E\otimes S)\leqno {\bf (2.3.13)}$$

where $S \i (E\otimes V^{\otimes 2})_{\Sigma_2}$ is the space
of relations in $A$ and the intersection is taken inside $X$.

\noindent {\sl Proof:} Straightforward. Left to the reader.

\vskip .2cm

Similarly, the algebra $U_{{\cal P}^!}(A^!)$ has the space of generators
$E^\vee \otimes V^\vee$ and the space of relations
$$E^\vee\otimes V^\vee\otimes E^\vee\otimes V^\vee \quad \cap\quad
(W^\bot \,\, + \,\, E^\vee \otimes S^\bot). \leqno {\bf (2.3.14)}$$
We want to show that (2.3.13) and (2.3.14) are the orthogonal complement
to each other. This is a particular case of the following
general lemma.

\proclaim (2.3.15) Lemma. Let $X$ be any $K$ - bimodule decomposed into
a direct sum of bimodules $X = M \oplus N$. Let $X^\vee = M^\vee 
\oplus N^\vee$ be the corresponding  decomposition of $X^\vee
= {\rm Hom}_K (X, K)$. Let $W\i X$ and $L\i N$ be any sub - bimodules.
Then the orthogonal complement (in $M^\vee)$ of
$M\cap (W+L)$ coincides with $M^\vee \cap (W^\bot_X + L^\bot_N)$,
where $W^\bot_X$ and $L^\bot_N$ are the orthogonal complements in 
$X$ and $N$, respectively.

\noindent {\sl Proof of the lemma:} Since $K$ is semisimple,
we can write $N = L\oplus P$ where $P$ is another sub-bimodule, so
$X = M\oplus L \oplus P$ and $X^\vee = M^\vee \oplus L^\vee \oplus P^\vee$.
If $Z \i X$ is any sub-bimodule then
$$ (M\cap Z)^\bot_M \quad = \quad {\rm Im}
 \left\{ (M\cap Z)^\bot_X \buildrel \pi
\over\rightarrow M^\vee\right\} \quad = \quad
M^\vee \cap\left( (M\cap Z)^\bot_X + L^\vee + P^\vee\right),
\leqno {\bf (2.3.16)}$$
where $\pi: X^\vee \rightarrow M^\vee$ is the projection dual to
the embedding $M\hookrightarrow X$, i.e., the projection along
$L^\vee \oplus P^\vee$. Let us apply this to $Z = W+L$ and note that
$$\left( M\cap (W+L)\right)^\bot_X \quad = \quad
M^\bot_Z + \left(W^\bot_X \cap L^\bot_X\right) \quad = \quad
L^\vee + P^\vee + \left( W^\bot_X \cap (M^\vee + P^\vee)\right).$$
We get that the RHS of (2.3.16) is equal to
$$M^\vee \cap \left( W^\bot_X \cap (M^\vee + P^\vee) + L^\vee + P^\vee\right).
\leqno {\bf (2.3.17)}$$
To finish the proof of the lemma, it remains to show that (2.3.17)
coincides with $M^\vee \cap (W^\bot_X + P^\vee)$. (Note that $P^\vee$
is the same as $L^\bot_N$.)

To show this, suppose that $m = w + p \in M^\vee \cap (W^\bot_X + P^\vee)$,
so that $w \in W^\bot_X, p\in P^\vee$. Then $w = m-p \in M^\vee + P^\vee$,
so writing $m = w + 0 + p$, we get that $m$ belongs to
(2.3.17). Conversely, let $m$ belong to (2..3.17), so
$m = w+l+p$ with $l\in L^\vee, p\in P^\vee$ and $w\in W^\bot_X\cap
(M^\vee + P^\vee)$, so $w = m' + p',
m'\in M^\vee, p' \in P^\vee$. Then we have $m = m' + p' + l + p$
and, since $X^\vee = M^\vee \oplus L^\vee\oplus P^\vee$,
we get $m=m', p+p' = 0, l=0$. This $m = w + p$, wo
$m \in M^\vee \cap(W^\bot_X + P^\vee)$. Lemma 2.3.15 and
 hence Theorem 2.3.11
are proven.

\hfill\vfill\eject

\centerline {\bf 3.  DUALITY FOR  $dg$ - OPERADS.}

\vskip 1.5cm

\centerline {\bf  3.1. $dg$ - operads. Generating maps.}

\vskip 1cm

\noindent {\bf (3.1.1)} Let $dgVect$ be the symmetric monoidal category
 of differential graded ($dg$ - ) vector spaces, i.e., of complexes
over the base field $k$. By definition, an object of $dgVect$ is a graded
vector space $V^\bullet$ together with a linear map
(differential) $d: V^\bullet \rightarrow V^\bullet$ of degree 1 such that
$d^2 = 0$. Morphisms are linear maps preserving gradings and differentials.
The tensor product of two complexes $V^\bullet$ and $W^\bullet$ is
defined by (1.4.3) with the differential  given by the Leibnitz rule
$$d(v\otimes w) = d(v) \otimes w + (-1)^i v \otimes d(w), \quad
v\in V^i, w\in W^j.$$
The symmetry isomorphism $V^\bullet \otimes W^\bullet \rightarrow
W^\bullet \otimes V^\bullet$ is the same as in the category $gVect^-$,
see (1.4.2 (c)). As usual, for a complex $V^\bullet \in dgVect$ we
define
the dual complex $V^*$ by
$$(V^*)^i = (V^{-i})^*, \quad d_{V^*} = (d_V)^* \leqno {\bf (3.1.2)}$$
and the shifted complex
$V^\bullet [i]$, $\i \in {\bf Z}$, by
$$(V^\bullet [i])^j = V^{i+j}, \quad d_{V^\bullet [i]} = (-1)^i d_V.
\leqno{\bf (3.1.3)}$$

An operad in the category $dgVect$ is called a $dg$ - {\it operad}
(over $k$). An algebra over a $dg$ - operad in the category
$dgVect$ will be called a $dg$ - algebra. Note that any $k$ - linear
operad ${\cal P}$ can be regarded as a $dg$ - operad  (each
${\cal P}(n)$ is placed in degree 0). 

\vskip .3cm

\noindent {\bf (3.1.4)} For any $dg$ - operad ${\cal P}$ the
collection of cohomology vector spaces $H^\bullet {\cal P}(n)$
form an operad $H^\bullet ({\cal P})$ in the category $gVect^-$.
A morphism of $dg$ - operads $f: {\cal P} \rightarrow {\cal Q}$ is
called a {\it quasi - isomorphism} if the induced morphism
$H^\bullet (f): H^\bullet ({\cal P}) \rightarrow H^\bullet ({\cal Q})$
is an isomorphism. Similarly for  $dg$ - algebras over $dg$ - operads.

\vskip .3cm

\noindent {\bf (3.1.5)} A $dg$ - operad $\cal P$ will be called
{\it admissible} if the following conditions hold:

\vskip .1cm

\item{(i)} Each ${\cal P}(n)$ is a finite - dimensional $dg$ - vector 
space. 

\item{(ii)} The space ${\cal P}(1)$ is concentrated
 in degree 0, and is a semisimple $k$ - algebra.

\vskip .2cm

\noindent Given a semisimple $k$ - algebra $K$, we denote by $dgOP(K)$
the category of admissible $dg$ - operads $\cal P$ with ${\cal P}(1) = K$
 and 
with morphisms identical on first components.
\vskip .3cm

\noindent {\bf (3.1.6)} Our next aim is to define a kind of
``generating function" for an admissible $dg$ - operad $\cal P$.
It will be convenient for the future to work in a slightly greater
generality.

\vskip .2cm

Let $K$ be a semisimple $k$ - algebra. By a $K - dg$ - collection
we mean a collection $E = \{E(n), n\geq 2\}$ of finite - dimensional
$dg$ - vector spaces $E(n)$ together with a left $\Sigma_n$ - action
and a structure of $(K, K^{\otimes n})$ - bimodule on each $E(n)$
which satisfy the compatibility condition identical to the 
one given in
(1.2.11). 

\vskip .2cm

Clearly if $\cal P$ is an admissible $dg$ - operad and $K = {\cal P}(1)$
then ${\cal P}(n), n\geq 2$, form a $K - dg$ - collection. Given any
$K - dg$ - collection $E$, one defines the {\it free $dg$ - operad}
$F(E)$ as in (2.1.1).

\vskip .2cm

By an $r$ {\it fold $dg$ - collection} we mean, similarly to (1.3.15),
a collection of finite - dimensional complexes
$$E^i(a_1, ..., a_r), \quad a_i \in {\bf Z}_+, 
i=1,...,r, \quad\sum a_i \geq 1$$
of $\Sigma_{a_1} \times ... \times \Sigma_{a_r}$ - modules
such that 
$$ E^i(0, ..., 0, \underbrace{1}_{j}, 0, ..., 0) = \cases{ k & for $i=j$;\cr
0 & for $i\neq j$\cr}$$
Given a semisimple $k$ - algebra $K$ with $r$ simple summands and a
$K - dg$ - collection $E$, we define the $r$ - fold $dg$ -
collection associated to $E$ by the formula identical to (1.3.16).

\vskip .2cm

\proclaim (3.1.8) Definition. Let $K$ be a semisimple $k$ - algebra,
$r$ be the number of simple summands in $K$ and $E$ be a $K - dg$ -
collection. The generating map of $E$ is the $r$ - tuple of
formal power series
$$g^{(i)}_E(x_1, ..., x_r) = \sum_{a_1, ..., a_r = 0}^\infty \chi \left[
E^i(a_1, ..., a_r)\right] {x_1^{a_1}\over a_1!} \cdot ...
\cdot {x_r^{a_r}\over a_r!}, \quad i=1, ..., r \leqno {\bf (3.1.9)}$$
where $\{E^i(a_1, ..., a_r)\}$ is the $r$ - fold $dg$ -
collection associated to $E$ and $\chi$ stands for
 the Euler characteristic.

The $r$ - tuple $g_{_E}(x) = (g^{(1)}_E(x), ..., g^{(r)}_E(x))$ will
be regarded as a formal map
$$g_{_E}: {\bf C}^r \rightarrow {\bf C}^r, \quad x = (x_1, ..., x_r) 
\mapsto (g^{(1)}_E(x), ..., g^{(r)}_E(x)). \leqno {\bf (3.1.10)}$$
Note that by (3.1.7) we have
$$g_E^{(i)}(x) = x_i + ({\rm higher \quad order \quad terms}).
\leqno {\bf (3.1.11)}$$

In particular, for any admissible $dg$ - operad $\cal P$ we have
its generating map $g_{\cal P}$.

\vskip .2cm

\noindent {\bf Special case:} If $\cal P$ is a $k$ - linear operad
with ${\cal P}(1) = k$ then its generating map is a single powr series
$$g_{\cal P}(x) = \sum_{n=1}^\infty {\rm dim} \,\, {\cal P}(n)
{x^n\over n!}.$$

\vskip .3cm

\noindent {\bf (3.1.12)} Examples. a) The operads ${\cal A}s, {\cal C}om$
and ${\cal L}ie$, see (1.3.7) - (1.3.9), are admissible operads
with $K = k$ and trivial $dg$ - structure. Therefore $r=1$ and the
generating maps of these operads are the following power series in
one variable:
$$g_{{\cal A}s}(x) = {1\over 1-x} - 1, \quad g_{{\cal C}om}(x) = 
e^x - 1, \quad g_{{\cal L}ie}(x) = - \log (1-x).$$

\vskip .2cm

b) Let $\cal A$ be the symmetric monoidal category of representations of
the group ${\bf Z}/2$. It has two (1-dimensional) irreducible objects:
ther trivial representation $I$ and the sign representation $J$. Taking
$X = I\oplus J$, we define an admissible $k$ - linear operad $\cal P$
with ${\cal P}(n) = {\rm Hom}_{\cal A}(X^{\otimes n}, X)$, see (1.3.12).
We have ${\cal P}(1) = k\oplus k$. As explained in (1.3.16),
 the 2-fold collection associated to $\cal P$ consists of Clebsch -
Gordan spaces
$${\rm Hom}_{\cal A}(I^{\otimes a} \otimes J^{\otimes b}, I) = \cases{ 0 &
for $b$ odd\cr
$k$& for $b$ even\cr}$$
$${\rm Hom}_{\cal A}(I^{\otimes a} \otimes J^{\otimes b}, J) = \cases{ 0 &
for $b$ even\cr
$k$& for $b$ odd\cr}.$$
Therefore the generating map of $\cal P$ is given by the formulas
$$g^{(1)}_{\cal P}(x_1, x_2) = \sum_{a,b} {\rm dim}\,\,\,
{\rm Hom}_{\cal A}(I^{\otimes a} \otimes J^{\otimes b}, I)
 {x_1^a\over a!} {x_2^b\over b!} = e^{x_1} \cosh (x_2) - 1,$$
$$g^{(2)}_{\cal P}(x_1, x_2) = \sum_{a,b} {\rm dim}\,\,\,
{\rm Hom}_{\cal A}(I^{\otimes a} \otimes J^{\otimes b}, J)
 {x_1^a\over a!} {x_2^b\over b!} = e^{x_1} \sinh (x_2).$$

\vskip .2cm

\noindent {\bf (3.1.13)} It is possible to refine the generating map of
a $K - dg$ - collection so as to take into account not only the
dimensions of the graded components of the complexes $E^i(a_1, ..., a_r)$
but also the  symmetric groups actions.

\vskip .2cm

Let ${\cal R}_n$ be the Grothendieck group of the category of finite -
dimensional representations of the symmetric group $\Sigma_n$ over
$k$. For any such representation $V$ we denote by $[V]$ its class in
${\cal R}_n$. The maps
$${\cal R}_m \otimes {\cal R}_n \rightarrow {\cal R}_{m+n}, \quad
[V] \otimes [W] \mapsto \left[ {\rm Ind}_{\Sigma_m \times \Sigma_n}
^{\Sigma_{m+n}} (V \otimes W) \right] \leqno {\bf (3.1.14)}$$
make ${\cal R} = \bigoplus_{n\geq 0} {\cal R}_n$ into a commutative
graded ring. It is well known [Macd] that $\cal R$ is isomorphic to
the ring ${\bf Z}[e_1, e_2, ...]$ of symmetric functions in infinitely
many variables. There is a natural ring homomorphism
$h: {\cal R} \rightarrow {\bf Z}$ defined by
$$h([V]) = {{\rm dim}\,\,\, V\over n!}, \quad [V] \in {\cal R}_n.
\leqno {\bf (3.1.15)}$$
The tensor product ${\cal R}_{a_1} \otimes ... \otimes {\cal R}_{a_r}$
is naturally identified with the Grothendieck group of representations of
 the Cartesian product
$\Sigma_{a_1} \times ... \times \Sigma_{a_r}$. It is a direct summand
of the ring ${\cal R}^{\otimes r} = \bigoplus
{\cal R}_{a_1} \otimes ... \otimes {\cal R}_{a_r}$. If $V^\bullet$ is any
finite - dimensional complex of $\Sigma_{a_1} \times ... \times \Sigma_{a_r}$
- modules then by $[V^\bullet] \in {\cal R}^{\otimes r}$ 
we denote the alternating sum of classes of $V^i$ in 
${\cal R}_{a_1} \otimes ... \otimes {\cal R}_{a_r} \i {\cal R}^{\otimes r}$. 

\vskip .3cm

\noindent {\bf (3.1.16)} Let now $K$ be a semisimple $k$ - algebra with
$r$ simple summands and $E$ be a $K - dg$ - collection. We define the
{\it refined generating map} of $E$ to be the $r$ - tuple of formal
power series
$$G^{(i)}_E(x_1, ..., x_r) = \sum \left[ E^i(a_1, ..., a_r)\right]
x_1^{a_1} ... x_r^{a_r} \quad \in \quad {\cal R}^{\otimes r}[[ x_1, ...,,
x_r]] \leqno {\bf (3.1.17)}$$
where, as in (3.1.8), $\{E^i(a_1, ..., a_r)\}$ is the $r$ - fold
$dg$ - collection associated to $E$ (note that there are no denominators). 

Let  $h^{\otimes r}: {\cal R}^{\otimes r} \rightarrow {\bf Z}$ be
the tensor power of the homomorphism (3.1.15) 
 $$h_r: {\cal R}^{\otimes r} [[x_1, ..., x_r]] \rightarrow {\bf Z}
[[x_1, ..., x_r]]\leqno {\bf (3.1.18)}$$
be the associated homomorphism of power series rings. Then the
numerical generating map of $E$ can be recovered as
$$g^{(i)}_E(x_1, ..., x_r) = h_r (G^{(i)}_E(x_1, ..., x_r)).
\leqno {\bf (3.1.19)}$$

\hfill\vfill\eject

\centerline {\bf 3.2. The cobar - duality.}

\vskip 1cm

\noindent {\bf (3.2.0)} For a finite - dimensional $k$ - vector space
$V$ we denote by ${\rm Det}(V)$ the top exterior power of $V$.

Let $T$ be a tree (1.1.1) We denote by ${\rm Ed}(T)$ the set of all
edges of $T$ except the output edge Out$(T)$. We denote by Det$(T)$
the 1-dimensional vector space Det$(k^{{\rm Ed}(T)})$. 
Similarly, let  ed$(T)$ be the set of internal edges of $T$
and let $\det (T) = {\rm Det}(k^{{\rm ed}(T)})$.
  The number of internal edges
of $T$ will be denoted $|T|$. 

\vskip .3cm

\noindent {\bf (3.2.1)} Let $\cal P$ be an admissible $dg$ - operad
(3.1.5), so  that $K = {\cal P}(1)$ is a semisimple $k$ - algebra.
For any $n\geq 2$ we construct a complex
$${\cal P}(n)^* \otimes {\rm det}(k^n) \buildrel \delta\over\rightarrow
\bigoplus_{ n - {\rm trees}\,\,\,T \atop |T| = 1} {\cal P}(T)^* \otimes
{\rm det}(T) \buildrel \delta\over\rightarrow ... \buildrel 
\delta\over\rightarrow
\bigoplus_{ n - {\rm trees} \,\,\,T \atop |T| = n-2} {\cal P}(T)^* \otimes
{\rm det}(T) \leqno {\bf (3.2.2)}$$
where the sums are over isomorphism classes of 
(reduced) $n$ - trees and ${\cal P}
(T)$ was defined in (1.2.13). The differential $\delta$ is defined
by its matrix elements
$$\delta_{T', T}: {\cal P}(T')^* \otimes {\rm det}(T') \longrightarrow
{\cal P}(T)^* \otimes {\rm det}(T) \leqno {\bf (3.2.3)}$$
where $T, T'$ are $n$ - trees, $|T| = i, |T'| = i-1$.
By definition, $\delta_{T', T} = 0$ unless $T' = T/e$ is obtained from
$T$ be contracting an internal edge $e$. If this is the case, we set
$$\delta_{T', T} = (\gamma_{T, T'})^* \otimes l_e\leqno {\bf (3.2.4)}$$
where $\gamma_{T, T'}$ is the composition map from (1.2.15) and the map
$l_e: {\rm det}(T') \rightarrow {\rm det}(T)$ is defined by the formula
$$l_e(f_1\wedge ...\wedge f_m) = e\wedge f_1 \wedge ... \wedge f_m.$$
In this formula we use the natural identification ed$(T) = {\rm ed}(T') 
\cup \{e\}$ and regard $e$ as a basis vector of $k^{{\rm ed}(T)}$.

\vskip .3cm

\noindent {\bf (3.2.5)} It is straightforward to verify that
 $\delta^2 = 0$,
i.e.,(3.2.2) is a complex. We normalize the grading of this complex
by placing the sum over $T$ with $|T| = i$ in  degree $i+1$.

Observe that, for a $dg$ - operad $\cal P$, each term of (3.2.2) is a $dg$ -
vector space whose  differential we denote by $d$. Clearly $d$ commutes
with $\delta$ so (3.2.2) is a complex of $dg$ - vector spaces i.e.,
a double complex. 

\vskip .3cm

\noindent {\bf (3.2.6)} We now define a collection of $dg$ - vector
spaces $C({\cal P})(n), n\geq 1$. For $n=1$ we put
$C({\cal P})(1) = K^{op}$, placed in degree 0 (with trivial differential).
For $n > 1$ we define $C({\cal P})(n)$ to be the total complex
(= $dg$ - vector space) associated to the double complex (3.2.2).

\proclaim (3.2.7) Theorem. a) The collection $C({\cal P}) = 
\{ C({\cal P})(n), n\geq 1\}$ has a natural structure of an admissible
$dg$ - operad. \hfill\break
b) The correspondence ${\cal P} \mapsto C({\cal P})$ extends to 
contravariant functor $C: dgOP(K) \rightarrow dgOP(K^{op})$ (notation
of (3.1.5)). This functor takes quasi - isomorphisms to quasi - isomorphisms.

\noindent {\sl Proof:} a) The complexes ${\cal P}(n),  n\geq 2$, 
obviously form a $K - dg$ - collection (3.1.6). Therefore the shifted
dual complexes ${\cal P}(n)^*[-1]$ form a $ K^{op} - dg$ - collection
(shift and duality are defined by (3.1.2-3)). We denote this collection
by ${\cal P}^*[-1]$. Thus we can form the free $dg$ - operad
$F({\cal P}^*[-1])$.

\proclaim (3.2.8) Lemma. For any $n$ the complex $F({\cal P}^*[-1])(n)$
is  isomorphic to the total complex of the double complex
(3.2.2), the latter taken  with the internal differential $d$ (induced by
that on $\cal P$) only.

The lemma is proved by immediate inspection. The only point that needs
explanation is the appearance of vector spaces det$(T)$ in (3.2.2).
The reason for this is part b) of the following lemma.

\proclaim (3.2.9) Lemma. \item{a)} Let $I$ be a finite set of $m$
elements and $W_i^\bullet, i\in I$ be $dg$ - vector spaces.
Then there is a canonical isomorphism
$$\phi: \bigotimes_{i\in I} \left( W_i^\bullet [-1]\right) 
\buildrel \approx\over\longrightarrow \left( \bigotimes_{i\in I}
W_i^\bullet\right) [-m] \otimes {\rm Det} (k^I).$$
\item{b)} Let $E = \{ E(n), n\geq 2\}$, be a $K - dg$ - collection,
$E[-1]$ be the collection of shifted $dg$ - vector spaces and
$T$ be a tree with $m$ vertices. Then there is a canonical
isomorphism
$$E[-1](T) \quad \cong \quad E(T)[-m] \otimes \det (T).$$

\noindent {\sl Proof:} a) To define $\bigotimes_{i\in I} W_i^\bullet$,
we should choose some ordering $(i_1, ..., i_m)$ on $I$ and
consider the product $W^\bullet_{i_1} \otimes ... \otimes W^\bullet_{i_m}$.
Any other ordering will give the product related to this one by
a uniquely defined isomorphism.  The same for $\bigotimes W_i^\bullet [-1]$.
Let now $w_i \in W^\bullet_i [-1]$ be some elements and let $w'_i$
denote the same elements but regarded as elements of $W_i^\bullet$. We
define
$$\phi( w_{i_1} \otimes ... \otimes w_{i_m}) =
w'_{i_1} \otimes ... \otimes w'_{i_m} \otimes (i_1 \wedge ... \wedge i_m).$$
The proof that $\phi$ is independent on the choice of ordering 
$(i_1, ..., i_m)$ follows from the definition of the symmetric monoidal
structure in $dgVect$ and is left to the reader.

\vskip .2cm

b) Let Vert$(T)$ be the set of vertices of $T$. Let $\Gamma \i T$ be the
subtree consisting of all vertices and all internal edges. This is
a contractible 1-dimensional simplicial complex so taking its chain
complex we get the exact sequence
$$0\rightarrow k^{{\rm ed}(T)} \rightarrow k^{{\rm Vert}(T)}
\rightarrow k\rightarrow 0$$
whence the space Det$(k^{{\rm Vert}(T)})$ is canonically identified
with $\det (T) = {\rm Det} (k^{{\rm ed}(T)})$. Now part (b)
of the lemma follows from part (a), since $E(T) = \bigotimes_{v\in {\rm Vert}
(T)} E({\rm In}(v))$. 

Lemmas 3.2.9 and 3.2.8 are proven.

\vskip .3cm

\noindent {\bf (3.2.10)}  By Lemma 3.2.8, the compositions in the
free operad $F({\cal P}^*[-1])$ give rise to
 maps of graded vector spaces
$$C({\cal P})(l) \otimes C({\cal P})(m_1) \otimes ... \otimes C({\cal P})
(m_l) 
\rightarrow C({\cal P})(m_1 + ... + m_l) \leqno {\bf (3.2.11)}$$
which satisfy the Leibnitz rule with respect to
internal differentials $d$ in $C({\cal P})(n)$. To complete the
proof of Theorem 3.2.7 (a), it remains to check that the maps
(3.2.11) satisfy  also the Leibnitz rule  with respect to the differentials
$\delta$ in $C({\cal P})(n)$ defined in (3.2.3). We leave this straightforward
checking to the reader. This completes the proof of
part (a) of Theorem 3.2.7. Part (b) is  straightforward. 

\vskip .3cm

\noindent {\bf (3.2.12)} We call the $dg$ - operad $C({\cal P})$ the
{\it cobar - construction} of $\cal P$. We define the {\it dual dg -
operad} ${\bf D}({\cal P})$ by
$${\bf D}({\cal P}) = C({\cal P}) \otimes \Lambda = \{ C({\cal P}) (n)
\otimes \Lambda (n)\}, \leqno {\bf (3.2.13)}$$
where $\Lambda$ is the determinant operad (1.3.21) and the product 
$\otimes$ is defined in (2.2.2). Recall that $\Lambda (n)$ is a
1-dimensional vector space in degree $(1-n)$ with  the sign action of
$\Sigma_n$. Hence ${\bf D}({\cal P})(n)$ is the $dg$ - vector space
associated to the following complex which differs from (3.2.2) by 
shifting the grading by $(1-n)$ and by replacing det by Det:
$$ ... \rightarrow \bigoplus_{ n - {\rm trees}\,\,\,T\atop |T| = n-3}
{\cal P}(T)^* \otimes {\rm Det} (T) \rightarrow 
\bigoplus_{ n - {\rm trees} \,\,\,T\atop |T| = n-2}
{\cal P}(T)^* \otimes {\rm Det} (T). \leqno {\bf (3.2.14)}$$
The grading in (3.2.14) is  arranged  so
that the rightmost term is placed in degree
0. 

It follows from the definitions that the correspondence ${\cal P}
\mapsto {\bf D}({\cal P})$  extends to a contravariant functor
${\bf D}: dgOP(K) \rightarrow dgOP(K^{op})$ taking quasi - isomorphisms
to quasi - isomorphisms. 

\vskip .3cm

\noindent {\bf (3.2.15) Example.} 
 Let ${\cal }P = {\cal C}om$ be the commutative operad
 so  that ${\cal C}om (n) = k$ for
every $n$, see (1.3.8). The $n$-th component of ${\bf D}(P)$ is the 
{\it tree
complex}
\par
$$\bigoplus_{|T| = 0} {\rm Det}(T) \longrightarrow 
\bigoplus_{|T|=1}
{\rm Det} (T)\rightarrow \ldots
\rightarrow \bigoplus _{|T| = n-2}{\rm Det}(T)  $$
 
 This is a special case of more general graph
 complexes
considered by Kontsevich [Kon 2].
 In the paper [BG 1] it was
proven by using Hodge theory that  for $k={\bf C}$ this complex  is exact 
everywhere except the
rightmost term and the cokernel of the rightmost differential is
naturally isomorphic to ${\cal L}ie(n)$, the $n$-th space of the 
Lie operad. Thus we have an isomorphism ${\bf D}({\cal C}om) \cong
{\cal L}ie$. We   give a purely algebraic proof of this
later  in Chapter 4.

\vskip .3cm

\proclaim (3.2.16) Theorem. For any admissible $dg$ - operad $\cal P$
there is a canonical quasi - isomorphism ${\bf D}({\bf D}({\cal P})) 
\rightarrow {\cal P}$.  

\noindent {\sl Proof:} 
Let us write down ${\bf D}({\bf D}({\cal P}))(n)$ explicitly. By definition
(3.2.14),
$${\bf D}({\bf D}({\cal P}))(n)  \quad = \quad \bigoplus_{ n -\,{\rm trees}
\,\,\,S}
\biggl[ \bigotimes _{v\in S}{\bf D}( {\cal P})({\rm In}(v))^* \otimes 
{\rm det}(S)\biggl].$$
Substituting here the definitions (1.1.6)
of  the value  on a set of the functor
corresponding to ${\bf D}({\cal P})$,
   and keeping track of cancellation of some Det - factors, we get
$${\bf D}({\bf D}(P)) (n) = \bigoplus_{T \geq T'}
{\bigotimes_{v\in T}{\cal P}({\rm In}(v))\over
\bigotimes_{w\in T'} {\rm Det}(T_w)}.\leqno {\bf (3.2.17)}$$
Here the summation is over the  isomorphism classes of pairs $T, T'$ of
$n$ -  trees such that $T \geq T'$ i.e., $T'$ can be obtained
from $T$ by contracting some (possibly empty) set of edges. 
For $w\in T'$  we denote by $T_w$  the subtree of $T$ contracted into
$w$. Division by a 1-dimensional vector space is understood as 
tensoring with the dual space.

The construction may be understood better using
 Fig. 9a where vertices $w$ of the tree
$T'$ are indicated  as big circles ("regions") containing inside
them the corresponding trees $T_w$. Let $T_{n}$ denote
 the unique $n$ - tree
 with 
 a single vertex.
The summand  in (3.2.17) corresponding to the
pair $T_n \geq T_n$ is nothing but ${\cal P}(n)$ (Fig.  9b):

\vbox to 6cm{}

 It is
straightforward to  verify that the summand ${\cal P}(n)$ is in fact a
quotient complex of the whole  complex
${\bf D}({\bf D}({\cal P}))(n)$,
 i.e., the projection to
this summand along all  other summands is a (surjective) morphism of
complexes which we denote by $f_n$.
 It is also easily verified that we get in this way  a
morphism of $dg$ - operads $f: {\bf D}({\bf D}({\cal P}))
\rightarrow {\cal P}$. Let us show 
that $f$ is a
quasi-isomorphism i.e., that each subcomplex Ker $(f_n) \i
{\bf D}({\bf D}({\cal P}))(n)$, is acyclic. Note that Ker $(f_n)$,
as  part of ${\bf D}({\bf D}({\cal P}))(n)$, is actually a triple
complex so its differential is a sum of three partial differentials
$d_1 + d_2 + d_3$. The differential $d_1$ is induced by the
differential in $\cal P$; the differential $d_2$ is induced by the
composition in $\cal P$ (which induces the second differential in ${\bf
D}({\cal P})$). Finally, $d_3$ is induced by the composition in 
${\bf
D}({\cal P})$, i.e., by the grafting of trees. 

\vskip .2cm

It is enough to show that Ker $f_n$ is acyclic with respect to $d_3$.
If $T$ is an  $n$ -  tree with more than one vertex then
the summands in (3.2.17) with  all $T'\leq T$ form a subcomplex
$K^\bullet_T \i ({\rm Ker}(f_n), d_3)$ and Ker $(f_n)$ is the
direct sum of such $K^\bullet_T$.   We  shall prove that each
$K^\bullet_T$ is acyclic.

\vskip .2cm

The differential $d_3$ in $K^\bullet_T$ consists purely in redrawing the
boundaries among regions in Fig. 8 (a). More precisely, $K^\bullet_T$
is the tensor product of the vector space ${\cal P}(T)$ and a purely
combinatorial complex $C^\bullet_T$ where
$$C^i_T = \bigoplus_{\matrix{T'\leq T\cr |T| - |T'| = i}}
\bigotimes_{w\in T'} {\rm Det}(T_w)^*.$$
Observe that  specification of a tree $T'\leq T$ is equivalent to a
specification of a subset of internal edges of $T$ which are
contracted in $T'$. We see that $C^\bullet_T$ is isomorphic to the
augmented chain complex of a simplex whose vertices correspond
 to internal edges of $T$. Thus,
 $C^\bullet_T$ and $K^\bullet_T$ are acyclic.
Theorem 3.2.16 follows. 

\vskip .3cm

\proclaim (3.2.18) Proposition. For any admissible $dg$ - operad $\cal P$
with ${\cal P}(1) = k$ there exists a natural morphism of $dg$ -
operads $\lambda: {\bf D}({\cal C}om) \rightarrow {\cal P} \otimes {\bf D}
({\cal P})$. In particular, for any $dg$ - algebra $A$ over $\cal P$
and any $dg$ - algebra $B$ over 
${\bf D}({\cal P})$ the tensor product $A \otimes_k B$ has a natural
structure of a ${\bf D}({\cal C}om)$ - algebra. 

This result is analogous to Corollary 2.2.9 (b). The precise
meaning of the analogy will be explained in Chapter 4.

\vskip .2cm

\noindent {\sl Proof:} Since ${\cal C}om (n) = k$ for every $n$,
we find that ${\bf D}({\cal C}om)(n)$ is
the $dg$ - vector space whose $(-i)$ - th graded piece is the
direct sum of  spaces Det$(T)$ corresponding to $n$ - trees $T$
with $n-2-i$ interior edges.
We define the required maps
$\lambda_n: {\bf D}({\cal C}om) (n) \rightarrow {\cal P}(n) 
\otimes {\bf D}({\cal P})(n)$ by prescribing their restrtictions on
the summands ${\rm Det}(T)$.
We define
$$\lambda_n|_{{\rm Det} (T)}: {\rm Det}(T) \rightarrow 
 {\cal P}(n) \otimes {\cal P}^*(T) \otimes
\det(T) \i {\bf D}({\cal P})(n)$$
to be Id$_{{\rm Det} (T)} \otimes \gamma_T^\dagger$,
where the element $\gamma_T^\dagger \in {\cal P}(n) \otimes {\cal P}^*(T)$
is the transpose of the map $\gamma_T: {\cal P}(T) \rightarrow
{\cal P}(n)$ given by the composition in $\cal P$, see (1.2.4).

This defines the maps $\lambda_n$. The proof that these maps form 
a morphism of operads, is straightforward and left to the reader.

\hfill\vfill\eject

\centerline {\bf 3.3. The generating map of the dual  $dg$ - operad.}

\vskip 1cm

\noindent {\bf (3.3.1)} Let $\cal P$ be an admissible $dg$ - operad and $Q = 
{\bf D}({\cal P})$ be its dual (3.2.13). Let $r$ be the number
 of simple summands of the semisimple algebra ${\cal P}(1) = {\cal Q}(1)^{op}$
and let $g_{\cal P}, g_{\cal Q}: {\bf C}^r \rightarrow {\bf C}^r$ be the
generating maps of ${\cal P, Q}$, see (3.1.8). Let also $G_{\cal P},
G_{\cal Q}$ be the refined generating maps, see (3.1.16).

The following is the main result of this section.

\proclaim (3.3.2) Theorem. a) We have the following
identity of formal maps  ${\bf C}^r \rightarrow {\bf C}^r$:
$$g_{\cal Q}( - g_{\cal P} (-x)) = x, \quad \quad x = (x_1,...,x_r).$$
b) The same identity holds for the refined generating maps
$G_{\cal P}$ and $G_{\cal Q}$.

\noindent {\bf (3.3.3) Example.} Let ${\cal P} = {\cal C}om$ be the
commutative operad. As we saw in (3.1.12), $g_{\cal P}(x) = e^x-1$.
By (3.2.15), ${\cal Q} = {\bf D}({\cal P})$ is quasi - isomorphic to the
operad ${\cal L}ie$ so $g_{\cal Q} = g_{{\cal L}ie}(x) = 
- \log (1-x)$. One sees that these two series satisfy Theorem 3.3.2.

\vskip .3cm

\noindent {\bf (3.3.4)} Let $K = {\cal P}(1)$. As a first step
towards the proof of our theorem we   describe the $r$ - fold
$dg$ -  collection associated to the $K^{op} - dg$ - collection 
${\bf D}({\cal P}) = \{ {\bf D}({\cal P})(n)\}$.

\vskip .2cm

By an $r$ - colored tree we understand a tree $T$ together with a function
$c$ 
(''coloring") from the set of all edges of $T$ to
$\{1, ..., r\}$. For such a tree $T$ and a vertex $v\in T$ let
${\rm In}_i(v)$  denote the set of input edges at $v$ of color $i$. 
By an $(a_1, ..., a_r)$ - tree we mean an $r$ - colored tree $T$ with $a_i$
inputs of color $i$, which are labelled by the numbers $1, 2,...., a_i$. 

\vskip .3cm

\noindent {\bf (3.3.5)} Let $\{{\cal P}^{i}(a_1, ..., a_r)\}$ be the
$r$ - fold $dg$ - collection associated to $\cal P$. For any $i$ we
use the $\Sigma_{a_1} \times ... \times \Sigma_{a_r}$ - action
on each ${\cal P}^{i}(a_1, ..., a_r)$ to define, similarly to
(1.1.6), a functor
$$(I_1, ..., I_r) \longmapsto {\cal P}^{i}(I_1, ..., I_r)$$
on  the category of $r$ - tuples of finite sets and bijections. Let also
$\{{\cal Q}^{i}(a_1, ..., a_r)\}$ be the $r$ - fold $dg$ - collection
associated to the dual $dg$ - operad ${\cal Q} = {\bf D}({\cal P})$. 
Recall that  $|T|$ denotes the number of internal
edges of  a tree $T$.

\proclaim (3.3.6) Proposition. Each 
${\cal Q}^{i}(a_1, ..., a_r)$ is isomorphic, 
as a graded vector space with
$\Sigma_{a_1} \times ... \times \Sigma_{a_r}$ - action, 
 to
$$\bigoplus_{ (a_1, ..., a_r) - {\rm trees} \,\,\,T 
\atop c({\rm Out}(T)) = i}
\bigotimes_{v\in T} {\cal P}^{c({\rm Out}(v))} ({\rm In}_1(v), ..., 
{\rm In}_r(v))^* \left[ |T| - \sum a_i - 2 \right] \otimes \det (T).$$

The proof is straightforward.

\vskip .3cm

\noindent {\bf (3.3.7)} To prove Theorem 3.3.2, we invoke
a purely algebraic
formula describing the inversion of a formal (or analytic) map
$g: {\bf C}^r \rightarrow {\bf C}^r$.
This formula is due to J. Towber, see [W], Th. 3.10 or [MTWW], Th. 2.13.
(We are grateful to D. Wright for pointing out the references to
this formula, which was rediscovered by  us.)

Suppose we have $r$ formal power series
$$g^{(i)}(x_1, ..., x_r) \,\, = \,\, x_i \,\, + \,\,
\sum_{\matrix {a_1,...,a_r \geq 0\cr a_1 + ... + a_r \geq 2}}
p^{(i)}(a_1,...,a_r) \,\,\, {x_1^{a_1}\over a_1!} \cdot ... \cdot 
{x_r^{a_r}\over a_r!}$$
 where $p^{(i)}(a_1,...,a_r)$ are some complex coefficients and let
$$y_i = g^{(i)}(x_1, ..., x_r)\leqno {\bf (3.3.8)}$$
be a formal change of variables given by these series.
Since $g^{(i)}(x) = x_i +$ terms of order $\geq 2$, this change of
variables is formally invertible, i.e., we can express $x_i$ as power
series in $y_1,...,y_n$. The question is  to find the coefficients
of these power series provided
 the coefficients $p^{(i)}(a_1,...,a_r)$
are known.

\vskip .2cm

\proclaim (3.3.9) Theorem. The inverse of the formal map (3.3.8)
is given by
$$x_i =  h_i (y_1, ..., y_r)  = y_i + \sum_{\matrix {a_1,...,a_r
 \geq 0\cr a_1 + ... + a_r \geq 2}}
q^{(i)}(a_1,...,a_r)\,\,\,{y_1^{a_1}\over a_1!} \cdot ... \cdot 
{y_r^{a_r}\over a_r!}\leqno {\bf (3.3.10)}$$
where
$$ q^{(i)}(a_1,...,a_r) = \sum_{m=0}^{a_1+ ... + a_r -1}
(-1)^{a_1 + ... + a_r - m}\biggl\{ \leqno {\bf (3.3.11)}$$
$$   \sum_{\matrix{ (a_1, ...,a_r) - {\rm trees}\,\,\, T\cr
|T| = m, c({\rm Out}(T)) = i}}
 \prod_{v\in T} p^{c \bigl( {\rm Out}(v)\bigl)}
\bigl(|{\rm In}_1(v)|, \,\, ..., \,\, |{\rm In}_r(v)|
\bigl)\biggl\}.$$

Several other inversion formulas for analytic maps can be found in 
[Sy] [Gd] [Jo] [BCW] [MTWW].

\vskip .3cm

\noindent {\bf (3.3.12)} We prefer to give here a simple direct proof of
Theorem 3.3.9, partly for the convenience of the reader, partly
because an intermediate lemma in the proof will be used later.
The proof proceeds
by direct substitution of the proposed answer into the
claimed equation. Suppose first that the series $h_i$
 (i.e., the coefficients
$q^{(i)}(a_1, ..., a_r)$ ) in (3.3.10) are arbitrary. 
Let us form the composition
$$ f^{(i)}(y) \,\, = \,\, g^{(i)}\bigl(h^{(1)}(y), ..., h^{(r)}(y)
\bigl),
\quad y = (y_1,...,y_r)$$
and write
$$f^{(i)}(y) = y_i + 
\sum_{\matrix {a_1,...,a_r \geq 0\cr a_1 + ... + a_r \geq 2}}
u^{(i)}(a_1,...,a_r)\,\,\,{x_1^{a_1}\over a_1!} \cdot ... \cdot 
{x_r^{a_r}\over a_r!}.$$
It is immediate to get  a general formula for the coefficients
$u^{(i)}(a_1,...,a_r)$.
Call a rooted tree $T$ {\it short} if $T$ has no consecutive
internal edges (see Fig. 10)

\vbox to 5cm{}

The lowest vertex of a tree $T$ (the one which is adjacent to the
output edge Out$(T)$) will be denoted by Lw$(T)$.

\proclaim (3.3.13) Lemma. We have
$$u^{(i)}(a_1,...,a_r) = \hbox to 10cm{}$$
$$\sum_{\matrix{{\rm short}\,\, (a_1,...,a_r) - {\rm trees}
\,\,\,T\cr c\bigl({\rm Out}(T)\bigl) = i}}\biggl\{ 
p^{(c({\rm Out}(T)))} \bigl(|{\rm In}_1 ({\rm Lw}(T))|, \,\,
..., \,\, |{\rm In}_r ({\rm Lw}(T))|\bigl) \times$$
$$ \times
\prod_{v\in T, v\neq {\rm Lw}(T)} q^{(c({\rm Out}(v)))}
\bigl(|{\rm In}_1 (v)|, \,\,
..., \,\, |{\rm In}_r (v)|\bigl)\biggl\}$$
where $T$ runs over short $(a_1,...,a_r)$ -  trees 
such that the vertex ${\rm Lw}(T)$ may have just
one input edge but all other vertices have $\geq 2$ input edges.

The lemma is proved by explicit substitution of  power series. 

\vskip .2cm

\noindent Let us now substitute into Lemma 3.3.13 the 
particular values of the
coefficients \overfullrule 0cm $q^{(i)}(a_1,...,a_r)$
 from (3.3.11). We get a formula
for $u^{(i)}(a_1,...,a_r)$ which involves a summation over all 
$(a_1,...,a_r)$  - labelled trees $T$, not necessarily short,
such that  the lowest vertex Lw$(T)$ is allowed to have one
input (whose color coincides with the color of Out$(T)$), 
but no other vertex is allowed to have one input.
That means that any summand in the arising formula for 
$u^{(i)}(a_1,...,a_r)$ will enter twice: it will
correspond once  to a tree whose lowest vertex has $\geq 2$ inputs and
 it will correspond  for the
second time to the same tree but with a  one - input
vertex appended at the bottom. These two summands will enter with 
opposite signs, hence  they will cancel each other.
 Therefore all the coefficients
$u^{(i)}(a_1,...,a_r)$ are equal to zero. So $f_i(y) = y_i$
and Theorem 3.3.9 is proven.

\vskip .2cm

This completes the proof of part (a) of Theorem 3.3.2. The proof of
part (b) is entirely similar and consists of repeating the
argument in the new context when the coefficients of power
series are not just numbers but elements of the ring ${\cal R}^{\otimes r}$.
We omit the details.

\vskip .3cm

\noindent {\bf (3.3.14) Example.} Let $b_n$ be the number of non - 
isomorphic binary $n$ - trees. Let $\cal P$ be  the $k$ - linear
operad  with ${\cal P}(1) = {\cal P}(2) = k$,
${\cal P}(n) = 0, n\geq 3$. Applying the definition, we see that the
complex ${\bf D}({\cal P})(n)$ consists of one vector space
$V(n)$, of dimension $b_n$, placed in degree 0. The generating maps
of $\cal P$ and ${\bf D}({\cal P})$ are therefore $g(x) = x + x^2/2$
and $h(x) = x + \sum_{n=2}^\infty b_n x^n/n!$, 
respectively. The identity $g(-h(-x))
= x$ yields  $h(x) = 1 - \sqrt{1-2x}$ whence $b_n = (2n-3)!! =
1\cdot 3\cdot 5 \cdot ... \cdot (2n-3)$.

\hfill\vfill\eject

\centerline {\bf 3.4. The configuration operad and duality.}

\vskip 1cm

\noindent {\bf (3.4.1)} Let $X$ be a topological space. A complex of
sheaves of $k$ - vector spaces on $X$ can be regarded as a sheaf with
values in the Abelian category $dgVect$ of $dg$ - vector spaces. Such
 objects will be referred to as $dg$ - {\it sheaves}.

For a $dg$ - sheaf ${\cal F}^\bullet$ on $X$ we denote by
 $R\Gamma(X, {\cal F}^\bullet)$ the $dg$ - vector space of
 global sections of the canonical Godement resolution of ${\cal F}^\bullet$,
see [Go] [KS]. The cohomology spaces of $R\Gamma(X, {\cal F}^\bullet)$ are
$H^i(X, {\cal F}^\bullet)$, the usual topological hypercohomology
with coefficients in ${\cal F}^\bullet$. 

Given two topological spaces $X_1, X_2$ and $dg$ - sheaves 
${\cal F}_i^\bullet$
on $X_i$, there is a natural morphism
$$R\Gamma(X_1, {\cal F}^\bullet_1) \otimes R\Gamma(X_2, {\cal F}_2^\bullet)
\rightarrow R\Gamma(X_1 \times X_2, {\cal F}_1^\bullet \otimes 
{\cal F}_2^\bullet) \leqno {\bf (3.4.2)}$$
which is a quasi - isomorphism. 

\vskip .3cm

\noindent {\bf (3.4.3)} Let $\cal Q$ be a topological operad. The notion of
a $dg$ - sheaf on $\cal Q$ is completely analogous to the notion of a
sheaf on $\cal Q$ in (1.5.3). Similarly to (1.5.9) we have:

\proclaim (3.4.4) Proposition. If $\cal Q$ is a topological operad and
${\cal F}^\bullet$ is a $dg$ - sheaf on $\cal Q$ then the collection
$$R\Gamma({\cal Q}, {\cal F}^\bullet) = \{ R\Gamma({\cal Q}(n),
{\cal F}^\bullet(n))^*, \,\,\, n\geq 1\}$$
forms a $dg$ - operad, so that the cohomology spaces
$H^\bullet({\cal Q}(n), {\cal F}^\bullet(n))^*$ form an operad
in the category $gVect^-$. 

\vskip .2cm

\noindent {\bf (3.4.5) Example: logarithmic forms of the configuration 
operad.} 
Let $k = {\bf C}$ be the field of
complex numbers and let $\cal M$ be the configuration operad (1.4.4).
For any $n\geq 1$ let $j_n: M_{0, n+1} \hookrightarrow {\cal M}(n)$
be the embedding of the open stratum and $\underline{\bf C}_{M_{0, n+1}}$
be the constant sheaf on $M_{0, n+1}$. Let $D(n) = {\cal M}(n) - 
M_{0, n+1}$. This is a divisor with normal crossings
in ${\cal M}(n)$, consisting of $2^{n-1} -1$ smooth components $\overline
{{\cal M}(T)}$, where $T$ is an $n$ - tree with exactly one
internal edge.  Let 
$\Omega^\bullet_{{\cal M}(n)} (\log D(n))$ be the corresponding 
logarithmic de Rham complex [D 1]. Consider the collection of shifted 
and twisted complexes
$$ \Omega^\bullet_{\cal M}(\log D) = \left\{ \Omega^\bullet_{{\cal M}(n)}
(\log D(n)) [n-2] \otimes {\rm Det} ({\bf C}^n), \quad n\geq 1 \right\}.
\leqno {\bf (3.4.6)}$$
This collection has a natural structure of
 a $dg$ - sheaf  on $\cal M$,
 whose structure maps (1.5.3)(ii) are
given by the {\it Poincar\'e residue} maps ([D 1], n. 3.1.5.2).
Let us briefly recall this notion in the generality we need.
Let $X$ be a smooth variety and $Y\i X$ be a divisor with
normal crossings which we assume to consist of smooth components $Y_1, ..., 
Y_M$.
Let $Z \i Y$ a codimension $m$ subvariety given by intersection
of some $m$ of these components: $Z = \bigcap_{i\in I} Y_i$, $|I| = m$,
and let $\gamma: Z \hookrightarrow X$ be the embedding. The Poincar\'e
residue is the map
$${\rm res}: \,\, \gamma^* \Omega^\bullet_{X}(\log Y) [m]
\otimes {\rm Det} ({\bf C}^I) \longrightarrow
\Omega^\bullet_Z(\log (Z\cap Y).\leqno {\bf (3.4.7)}$$
The appearance of Det$({\bf C}^I)$ stems from the fact that the
operations of taking the residues along individual hypersurfaces $Y_i,
i\in I$ anticommute with each other due to the residue theorem. 

Returning to our situation,
let $m_1, ..., m_l \geq 1$ be given. The structure map
of the operad $\cal M$
$$\gamma_{m_1, ..., m_l}: {\cal M}(l) \times {\cal M}(m_1)
\times ... \times {\cal M}(m_l) \rightarrow {\cal M}(m_1 + ... + m_l),$$
described in (1.4.4), is a closed embedding whose image is the
intersection of several components of the divisor $D(m_1 + ... + m_l)$.
If we denote the set of these components by $I$ then we have a canonical
identification
$${\rm Det}({\bf C}^I) \cong \left({\rm Det} ({\bf C}^l)\otimes 
\bigotimes {\rm Det}({\bf C}^{m_i})\right)^* \otimes {\rm Det}
({\bf C}^{m_1 + ... + m_l}).$$ 
Thus the  maps (3.4.6) indeed supply the necessary structure.

\vskip .3cm

\noindent {\bf (3.4.8) The gravity operad ${\cal G}$.}
We denote the $dg$ - operad $R\Gamma({\cal M}, \Omega^\bullet_{\cal M}
(\log D))^*$ simply by $\tilde {\cal G}$ and by $\cal G$ we denote
the cohomology operad of $\tilde{\cal G}$ i.e., the
operad in the category of graded vector spaces $gVect^-$ given by
${\cal G}(n) = H^\bullet (\tilde {\cal G}(n))$. Following a
suggestion of E. Getzler, we call $\cal G$ the {\it gravity operad}.
 Since $\Omega^\bullet_{{\cal
M}(n)}(\log D(n))$ is quasi - isomorphic to the direct image
$Rj_{n*}\underline{\bf C}_{M_{0, n+1}}$, we have
$$H^i({\cal M}(n), \Omega^\bullet (\log D(n))) = H^i(M_{0, n+1}, {\bf C}).$$
Thus the $i$ - th component of the graded vector space ${\cal G}(n)$
is $H_{n-2-i}(M_{0, n+1}, {\bf C})\otimes {\rm Det} ({\bf C}^n)$.

\proclaim (3.4.9) Proposition.
There is a quasi - isomorphism of $dg$ - operads $\tilde {\cal G} 
\rightarrow {\cal G}$ (where $\cal G$ is considered with zero
differential).

\noindent {\sl Proof:} Let $L^i(n)$ be the space of global
logarithmic $i$ - forms on ${\cal M}(n)$. It is known
[ESV] that $L^i(n)$ consists of closed forms and is naturally
identified with $H^i(M_{0, n+1}, {\bf C})$. Thus the embedding of the
graded vector space $L^\bullet (n) = \bigoplus L^i(n)$ with
zero differential into $R\Gamma({\cal M}(n), \Omega^\bullet_
{{\cal M}(n)} (\log D(n)))$ is a quasi - isomorphism. Denote
this embedding by $\phi_n$.
Clearly the Poincar\'e residue homomorphisms preserve 
global logarithmic forms so they make  the collection of 
graded vector spaces $L^\bullet (n)^* [-n+2]\otimes {\rm Det}({\bf C}^n)$
 into an operad
which is isomorphic to $\cal G$. Taking duals of the embeddings
$\phi_n$ above we get the required quasi - isomorphism
$\tilde {\cal G} \rightarrow {\cal G}$. 

\vskip .3cm

\noindent {\bf (3.4.10)} Let $\cal P$ be an admissible $dg$ - operad.
As in (1.5.2) we associate to $\cal P$ a $dg$ - sheaf ${\cal F}_{\cal P}$
on the configuration operad $\cal M$. By (3.4.4) the complexes
$R\Gamma({\cal M}(n), {\cal F}_{\cal P}(n))^*$ form another
$dg$ - operad. It is described as follows.

\proclaim (3.4.11) Theorem. Suppose that $k={\bf C}$. Then the operad 
$R\Gamma({\cal M}, {\cal F}_{\cal P})^*$ is quasi - isomorphic to
${\bf D}({\cal P} \otimes {\cal G})$ where 
{\bf D} is the duality for
$dg$ - operads introduced in (3.2.13). 

\noindent {\sl Proof:} 
 Let $X$ be a finite (compact)
CW - complex and $S = \{X_\alpha\}$
be its  Whitney stratification by locally closed CW - subcomplexes.
 We denote by $j_\alpha: X_\alpha \hookrightarrow X$ the embeddings
of the strata. Let also $X_{\geq m}$ be the union of strata of 
dimensions $\geq m$. This is an open subset in $X$ and we denote
by $j_{\geq m}: X_{\geq m} \hookrightarrow X$
the embedding.
Let ${\cal F}$ be any ($dg$ -)
sheaf on $X$. Let ${\cal F}_\alpha = j_\alpha^*
{\cal F}$ be the restrictions of ${\cal F}$ to the strata. Note that 
${\cal F}$ has  a decreasing filtration
$${\cal F}\,\,\, = \,\,\, (j_{\geq 0})_! \,\, j_{\geq 0}^* 
\,\, {\cal F} \,\,\, \supset
(j_{\geq 1})_! \,\, j_{\geq 1}^* \,\, {\cal F} \,\,\, \supset 
\,\, ...\leqno {\bf (3.4.12)}$$
with quotients
$$\bigoplus_{{\rm dim} X_\alpha = m} j_{\alpha !} j_\alpha^* {\cal F}.$$
Here $j_{\alpha !}$ are the direct images with proper support 
(extensions by zero).
We can regard this filtration as a Postnikov system  realizing
${\cal F}$ as a convolution [Ka 4] of the following complex of objects of
$D^b(Sh_X)$, the derived category of sheaves on $X$:
$$\bigoplus_{{\rm dim} X_\alpha = 0} j_{\alpha !}j_\alpha^* {\cal F}
\rightarrow \bigoplus_{{\rm dim} X_\alpha = 1} j_{\alpha !}j_{\alpha}^*
{\cal F}
[1] \rightarrow 
\bigoplus_{{\rm dim} X_\alpha = 2} j_{\alpha !}j_{\alpha}^*
{\cal F}[2]
 \rightarrow ...\leqno {\bf (3.4.13)}$$
By replacing $j_{\alpha!} j_\alpha^*{\cal F}$ with appropriate 
injective resolutions,
we can realize (3.4.13) as an actual double complex of sheaves on $X$.
By applying the functor $R\Gamma(X, -)$ to (3.4.13) we get that $R\Gamma(X,
{\cal F})$ is the total complex of the double complex
$$\bigoplus_{{\rm dim} X_\alpha = 0} R\Gamma_c( X_\alpha,
j_\alpha^*{\cal F}) \rightarrow \bigoplus_{{\rm dim} X_\alpha = 1}
R\Gamma_c( X_\alpha,
j_\alpha^*{\cal F}) [1] \rightarrow ...\leqno {\bf (3.4.14)}$$
where $R\Gamma_c$ is the derived functor of global sections with
compact support. The horizontal grading in this complex is so
normalized that the sum with $\dim X_\alpha = 0$ has horizontal
degree 0. 

\vskip .2cm

Assume now that $\cal F$ is $S$ - combinatorial (1.5.1)
 and given by 
($dg$-) vector spaces $F_\alpha$. Then we have the equalities
 $$R\Gamma_c(X_\alpha, j_\alpha^*{\cal F}) = R\Gamma_c(X_\alpha, {\bf C})
\otimes F_\alpha.\leqno {\bf (3.4.15)}$$ 

We specialize now to the case 
 when $X = {\cal M}(n)$, the stratification $S$
consists of ${\cal M}(T)$ and ${\cal F} = {\cal F}_{\cal P}(n)$
is the sheaf corresponding to an admissible $dg$ - operad ${\cal P}$.
For an $n$ - tree $\cal T$ let $D(T) \i \overline {{\cal M}(T)}$ be the
complement to ${\cal M}(T)$. This is a divisor with normal crossing.
As usual, we denote by $|T|$ the number of internal edges in $T$ 
so dim$_{\bf C} \,\,\, {\cal M}(T) = n-2 - |T|$. 
Note that we have the following quasi - isomorphisms
$$R\Gamma_c( {\cal M}(T), {\bf C}) \cong  R\Gamma({\cal M}(T), 
{\bf C})^* [-2(n-2-|T|)] 
 \cong $$
$$\cong R\Gamma(\overline {{\cal M}(T)}, \Omega^\bullet_
{\overline {{\cal M}(T)}}(\log D(T))^* 
[-2(n-2-|T|)].\leqno {\bf (3.4.16)}$$
Here the first quasi - isomorphism is the Poincar\'e duality.
 The second 
quasi - isomorphism is
the consequence of the fact  [D 1] that the logarithmic de Rham complex
is quasi - isomorphic to the full direct image of the constant sheaf. 
Note that in virtue of the product decompositions  (1.4.8) and
the definition of the $dg$ - operad $\tilde{ \cal G}$ 
we have a quasi - isomorphism 
$$R\Gamma \left(\overline {{\cal M}(T)}, \Omega^\bullet_
{\overline {{\cal M}(T)}}(\log D(T)\right)^* \cong 
\tilde {\cal G}(T) [n-2 - |T|] \otimes {\rm Det}(T),$$
where Det$(T)$ was introduced in (3.2.0). 
Taking into account Proposition 3.4.9, we see that
 the double complex (3.4.14), i..e., $R\Gamma( {\cal M}(n),
{\cal F}_{\cal P}(n))$, can be replaced by the double complex
$$\bigoplus_{|T| = 0} {\cal P}(T) \otimes
{\cal G}(T) \otimes {\rm Det} (T) \rightarrow 
\bigoplus_{|T| = 1} {\cal P}(T) \otimes
{\cal G}(T) \otimes {\rm Det} (T) \rightarrow ...$$
So the dual complex will be ${\bf D}({\cal P} \otimes {\cal G})(n)$
as claimed.

\vskip .3cm

\noindent {\bf (3.4.17) Example.} Taking ${\cal P} = {\cal C}om$,
the commutative operad, we get that every ${\cal F}_{\cal P}(n)$ is
the constant sheaf {\bf C} on ${\cal M}(n)$. Thus the operad formed by
 the total homology spaces $H_\bullet ({\cal M}(n), {\bf C})$ is
quasi - isomorphic to ${\bf D}({\cal G})$.

\hfill\vfill\eject

\centerline {\bf 3.5. The building co - operad and duality.}

\vskip 1cm

\noindent {\bf (3.5.1)} By a {\it cooperad} in a symmetric 
monoidal category
$(\cal A, \otimes)$ we mean an operad in the opposite category
 ${\cal A}^{op}$.
Explicitly, a cooperad ${\cal B}$ is a collection of objects
${\cal B}(n) \in {\cal A}, n\geq 1$ with $\Sigma_n$ - action on each
${\cal B}_n$ and morphisms
$$\rho_{m_1, ..., m_l}: {\cal B}(m_1 + ... + m_l) \rightarrow 
{\cal B}(l) \otimes {\cal B}(m_1) \otimes ... \otimes {\cal B}(m_l)
\leqno {\bf (3.5.2)}$$
satisfying the conditions dual to the corresponding conditions for an
operad. 

\vskip .2cm

If a collection ${\cal C} = \{ {\cal C}(n)\}$ forms a cooperad
in
 the category of ($dg$ -) vector spaces
then the dual $dg$ - vector spaces ${\cal C}(n)^*$ form a 
$dg$ - operad ${\cal C}^*$.
 A $dg$ - cooperad ${\cal C}$  is called {\it admissible} if 
${\cal C}^*$ is an admissible operad in the sense of (3.1.5). 

Cooperads in the category of topological spaces will be called
{\it topological cooperads}. 

\vskip .3cm

\noindent {\bf (3.5.2) The building ${\cal W}(n)$.} Let $T$ be a $n$ -
tree. By a {\it metric} on $T$ we mean an assignment of a positive
number (length) to each internal edge of $T$. The external edges
 can be thought  of as 
as having length 1. 
Let ${\cal W}(n)$ be the set of isometry classes of all $n$ - 
trees with metrics. This set has a natural topology.
In plain words, when the length of some edge is going to 0, we
say that the limit tree is obtained by contracting this edge into
a point. 
For any  $n$ -  tree $T$  (without metric) 
let  ${\cal W}^0(T)$ denote
the subset
in ${\cal W}(T)$ consisting of all  trees  with metric 
isomorphic to $T$.
This subset is clearly a cell (a product of
several copies of ${\bf R}_+$, corresponding to internal edges of $T$). 
Thus ${\cal W}(n)$ is a union of these non - compact cells. For example,
the space ${\cal W}(3)$ is the union of three half - lines glued along
a common end (Fig. 11):

\vbox to 4cm{}

The minimal cell in ${\cal W}(n)$ has  dimension 0 and  corresponds
 to the $n$ -  tree  $T_n$ with
 a single vertex (and no internal edges). Maximal cells of ${\cal W}
(n)$ correspond to binary trees. Thus, the cells in ${\cal W}(n)$
are parametrized by the same set as the strata in the moduli space
${\cal M}(n)$, see (1.4.5). But the closure relations
among the cells are {\it dual} to those among the strata in
${\cal M}(n)$,
 where the maximal stratum corresponds to the
tree with one vertex and 0 -dimensional strata correspond to
binary trees.

\vskip .2cm

We call ${\cal W}(n)$ the {\it building} of $n$ - trees since it
is in many respects similar to the Bruhat - Tits building.

\vskip .3cm

\noindent {\bf (3.5.3) The building cooperad.} Suppose given  natural
numbers $m_1, ..., m_l$. We define a map
$$\rho_{m_1, ..., m_l}: {\cal W}(m_1 + ... + m_l) \rightarrow 
{\cal W}(l) \times
{\cal W}(m_1) \times ... \times {\cal W}(m_l) $$
as follows. Let $T(m_1, ..., m_l)$ be the $(m_1 + ... + m_l)$
-tree drawn in Fig. 3 above.
 Let $(T', \mu')$ be any other metrized $(m_1 + ... + m_l)$
-  tree. Let $T''$
 be the maximal tree which is obtained from both $T'$ and $T(m_1,
..., m_l)$ by contracting  edges. The tree $T''$ is naturally 
divided into blocks $T''_1, ..., T''_l, T''_\infty$ where
 $T''_i$ has $m_i$ inputs and $T''_\infty$ has $l$ inputs, see
Fig.12.

\vbox to 5cm{}

  Each block $T''_\nu, \nu = 1,..., \infty$, is naturally
equipped with a metric $\mu''_\nu$ on its internal edges. We
redefine the lengths of edges which become loose edges in $T''_\nu$
be leting their lengths be equal to 1. The map $\rho_{m_1,...,m_l}$
takes
$$(T', \mu') \,\,\, \longmapsto \,\,\, \bigl((T''_\infty, \mu''_\infty),
 (T''_1, \mu''_1), ..., 
(T''_n, \mu''_n) \bigl).$$

\proclaim (3.5.4) Proposition. The collection of
maps $\rho_{m_1,...,m_l}$
and  the natural actions of $\Sigma_n$ on ${\cal W}(n),
n\geq 1$, define on the collection of ${\cal W}(n)$ the  structure
of a topological cooperad which will be denoted by ${\cal W}$.

We call $\cal W$ the {\it building cooperad.}
Note that the structure maps $\rho_{m_1,...,m_l}$ for
${\cal W}$ are surjective, whereas the structure maps for the 
operad ${\cal M}$ are, dually, injective.

\vskip .2cm

Just as any operad gives rise to a sheaf on the topological operad ${\cal M}$
(Theorem 1.5.11),
any cooperad gives a sheaf on the cooperad ${\cal W}$. Let us give the
corresponding definitions.

\proclaim  (3.5.5) Definition. Let $\cal B$ is a topological cooperad and let
$$\rho_{m_1,...,m_l}: {\cal B}(m_1+...+m_l) \rightarrow 
{\cal B}(l) \times {\cal B}(m_1) \times ...
\times {\cal B}(m_l) $$
be its structure maps (3.5.2). A sheaf (resp. a $dg$ - sheaf) 
${\cal F}$ on $\cal B$ consists of the following data:
\item{(i)} A collection
of $\Sigma_n$ - equivariant sheaves (resp. $dg$ -  sheaves)
 ${\cal F} (n)$ on ${\cal B}(n)$, one for each $n \geq 1$;
\item{(ii)}
  Homomorphisms
of sheaves on ${\cal B}(m_1 + ... + m_l)$
$$\mu_{m_1,...,m_l}: \rho_{m_1,...,m_l}^*( {\cal F}(l) \otimes
{\cal F}({m_1}) 
\otimes ... \otimes 
{\cal F}(m_l)) \rightarrow {\cal F} (m_1+
... + m_l).$$
These data should satisfy the  compatibility conditions dual to those
given in (1.5.4) - (1.5.6).

If, for example, all the spaces ${\cal B}(n)$ consist of
a single point
then a $dg$ - sheaf on $\cal B$ is the same as a $dg$ - operad.
More generally, we have the following.

\proclaim (3.5.6) Proposition. If $\cal B$ is a topological cooperad and
${\cal F}$ is a $dg$ - sheaf on $\cal B$ then the collection of
complexes $R\Gamma({\cal B}(n), {\cal F}(n))$ forms a $dg$ - operad.  
Similarly, the complexes
$R\Gamma_c({\cal B}(n), {\cal F}(n))$ (the derived functors of sections
with compact support) form a $dg$ - operad.

Let  $R\Gamma ({\cal B}, {\cal F})$
and $R\Gamma_c ({\cal B}, {\cal F})$ denote the $dg$ - 
operads thus obtained.

\proclaim  (3.5.7) Definition. 
 A sheaf ${\cal F}$ on a topological cooperad $\cal B$ is called
an iso - sheaf if the morphisms
$$ {\cal F}(l) \otimes {\cal F} ({m_1}) \otimes ... \otimes 
{\cal F}({m_l})  \rightarrow (\rho_{m_1,...,m_l})_*
{\cal F}(m_1 + ... + m_l)$$
obtained from $\mu_{m_1, ..., m_l}$ by adjunction, are isomorphisms
and higher direct images 
\hfill\break $R^j(\rho_{m_1,...,m_l})_*
{\cal F}( m_1 + ... + m_l)$ vanish for $j\geq 1$.

\vskip .2cm

\noindent {\bf (3.5.8)} Let $\cal Q$ be any $k$ - linear cooperad.
Then ${\cal Q}(1)$ is a $k$ - coalgebra and each ${\cal Q}(n)$
is a $({\cal Q}(1)^{\otimes n}, {\cal Q}(1))$ - bi - comodule.
As in (1.1.5), we extend the collection of $\Sigma_n$ - modules
${\cal Q}(n)$ to a functor $I \mapsto {\cal Q}(I)$ on 
finite sets and bijections. For any tree $T$ we define, similarly to
(1.2.13), the space
$${\cal Q}(T) = \bigodot_{v\in T}{}_{_{{\cal Q}(1)}} {\cal Q}({\rm In}(v))$$
where $\bigodot$ is the cotensor product
of comodules over a coalgebra. If ${\cal Q}(1) = k$
then ${\cal Q}(T) = \bigotimes_{v\in T} {\cal Q}({\rm In}(v))$
is the usual tensor product over $k$. If all ${\cal Q}(n)$
are finite  - dimensional over $k$ then we can use the operad ${\cal Q}^*$
to define ${\cal Q}(T) = ({\cal Q}^*(T))^*$ where ${\cal Q}^*(T)$
was defined by (1.2.13).

\vskip .2cm

If $T \leq T'$ then
 the cooperad structure on $\cal Q$
defines a linear map $\rho_{T,T'}: {\cal Q} (T) \rightarrow 
{\cal Q} (T')$. The condition
$T\leq T'$ means that the cell ${\cal W}^0(T)$ is contained in the
closure of ${\cal W}^0(T')$. Therefore 
the  maps $\rho_{T, T'}$ give rise
to a  combinatorial (1.5.1)
sheaf ${\cal G}_{\cal Q}(n)$ on ${\cal W}(n)$ 
with fiber ${\cal Q} (T)$ on  ${\cal W}^0(T)$. 

Here is the ``dual" of Theorem 1.5.11..

\proclaim (3.5.9) Theorem. Let $\cal Q$ be a $k$ - linear operad. Then:
\hfill\break
a) The sheaves ${\cal G}_{\cal Q}(n)$ on the spaces ${\cal W}(n)$
form a sheaf ${\cal G}_{\cal Q}$ on the cooperad $\cal W$. \hfill\break
b) If ${\cal Q}(1) = k$ then ${\cal G}_{\cal Q}$ is an iso - sheaf.
\hfill\break
c) Any iso - sheaf ${\cal G}$ on $\cal W$ such that each ${\cal G}(n)$
is constant on each cell ${\cal W}^0(T)$, has the form ${\cal G}_{\cal Q}$
for some $k$ - linear cooperad $\cal Q$ with ${\cal Q}(1) = k$.

\noindent The proof is straightforward and left to the reader.

\vskip .2cm

 In a similar way, for 
any $dg$ - operad $\cal Q$ we construct a $dg$ - sheaf ${\cal G}_{\cal Q}
^\bullet$ on $\cal W$.

\vskip .3cm

\noindent {\bf (3.5.10)}Let $X$ be a CW - complex and $S$ be
 its stratification into strata 
which are topological manifolds. 
Recall [KS] that the Verdier duality gives a contravariant
functor ${\cal F}^\bullet \mapsto {\bf V}({\cal F}^\bullet)$ from the
derived category of $dg$ - sheaves on $X$ with $S$ - constructible
cohomology into itself. 

The following result gives two 
different sheaf - theoretic interpretations of the duality
{\bf D} on  $dg$ - operads introduced in \S 3.2. 

\proclaim (3.5.11) Theorem. Let $\cal Q$ be an admissible $dg$ - cooperad
(3.5.1) so that ${\cal Q}^*$ is an admissible $dg$ - operad. Then:
\hfill\break
a) There is a natural quasi - isomorphism of $dg$ - operads
$${\bf D}({\cal Q}^*) \quad \cong \quad R\Gamma_c({\cal W}, {\cal G}^\bullet
_{\cal Q}) \otimes \Lambda$$
where $\Lambda$ is the determinant $dg$ - operad (1.3.21).
\hfill\break
b) For any $n$, the $dg$ - sheaves  ${\cal G}_{\cal Q}^\bullet (n)$
and ${\cal G}_{{\bf D}({\cal Q}^*)^*}$ on ${\cal W}(n)$
are Verdier dual to each other. 

\noindent {\sl Proof:} Let $(X, S)$ be any space stratified into
cells. Giving an $S$ - combinatoial   ($dg$-) sheaf ${\cal F}$
 on $X$ is the
same (1.5.11) as giving ($dg$ -) vector spaces $F_\sigma \bigl(=
R\Gamma(\sigma, {\cal F}\bigl)$ together with
 generalization maps
$g_{\sigma\tau}: F_\sigma \rightarrow F_\tau$ for $\sigma \i \bar\tau$
satisfying  transitivity conditions. For any cell $\sigma$ let
 ${\rm OR}(\sigma) = H^{{\rm dim}(\sigma)}_c (\sigma, k)$ be
its (1-dimensional) orientation space. Note that for the cell
${\cal W}^0(T)$ in the building ${\cal W}(n)$ corresponding to an
$n$ - tree $T$, the space OR$({\cal W}^0(T))$ can be naturally
identified with the space $\det (T)$, see (3.2.0). Theorem
3.5.11 is a consequence of the following well known combinatorial  
construction of the Verdier duality and hypercohomology
functors. 

\proclaim (3.5.12) Proposition. a) Let $(X,S)$ be as above and
$\cal F$ be an $S$ - combinatorial $dg$ - sheaf  on $X$ given
by $dg$ - vector spaces $F_\sigma$ and maps $g_{\sigma\tau}$.
Then:
\hfill\break
a) The $dg$ - vector space $R\Gamma_c(X, {\cal F})$ 
is quasi - isomorphic  naturally to ( the total $dg$ - vector 
space arising from) the complex 
$$\bigoplus_{{\rm dim}(\sigma) = 0} F_\sigma \otimes {\rm OR}(\sigma)
\longrightarrow 
\bigoplus_{{\rm dim}(\sigma) = 1} F_\sigma \otimes {\rm OR}(\sigma)
\longrightarrow ...$$
where the sum over $\sigma$ with ${\rm dim}(\sigma) = m$ is placed
in degree $m$. \hfill\break
b) The Verdier dual $dg$ - sheaf ${\bf V}({\cal F})$ is 
represented by the collection of
by $dg$ - vector spaces associated to complexes 
$${\bf V}({\cal F})_\sigma = \biggl\{ ... \rightarrow \hbox to 10cm{} $$
$$\hskip -.5cm  
\bigoplus_{\matrix{\tau\supset \sigma\cr {\rm dim} (\tau) = {\rm dim}
(\sigma) + 2}} { F}_\tau^* \otimes {\rm OR}(\tau)\rightarrow 
\bigoplus_{\matrix{\tau\supset \sigma\cr {\rm dim} (\tau) = {\rm dim}
(\sigma) + 1}} { F}_\tau^* \otimes {\rm OR}(\tau) \rightarrow 
{ F}_\sigma^*\otimes {\rm OR}
(\sigma) \biggl\},$$
where the sum over $\tau$ with dim$(\tau) = {\rm dim}(\sigma) + m$
is placed in degree $(-m)$.

This completes the proof of Theorem 3.5.11.

\vskip .2cm

The isomorphism ${\bf D}({\cal C}om) \cong {\cal L}ie$, see (3.2.15),
yields the following.

\proclaim (3.5.13) Corollary. The cohomology of ${\cal W}(n)$
with compact support (and constant coefficients) are as follows:
$$H^i_c({\cal W}(n), k) = \cases {$0$, & $i\neq n-2$\cr ${\cal L}ie (n)$,&
$i=n-2$.\cr}$$

This shows that the space ${\cal L}ie(n)$ is analogous to the Steinberg
representation of the group $GL_n({\bf F}_q)$. Note that the
standard Steinberg representation of this group, see  [Lu]
has dimension $q^{n(n-1)/2} = q^1\cdot q^2 \cdot ... \cdot q^{n-1}$.
The modified Steinberg (discrete series) representation introduced by
Lusztig [Lu] has dimension $(q-1)(q^2-1) ... (q^{n-1} - 1)$.
Both numbers can be seen as $q$ - analogs of $(n-1)! = {\rm dim}\,\,\,
{\cal L}ie (n)$. The role of the Lie operad as the dualizing module 
in our theory (2.2.9) is analogous to the role of the Steinberg
representation in the Deligne - Lusztig theory   [DL 1-3].

\hfill\vfill\eject

\centerline {\bf 4. KOSZUL OPERADS.}

\vskip 1.5cm

\centerline {\bf 4.1. Koszul operads and Koszul complexes.}

\vskip 1cm

\noindent {\bf (4.1.1)} Let $K$ be a semisimple $k$ - algebra and
 ${\cal P}  = {\cal P}(K, E, R)$ be a quadratic operad. 
Let ${\cal P}^! = {\cal P}(K^{op}, E^\vee, R^\bot)$ be
 the dual quadratic operad
and ${\bf D}({\cal P})$ be the  dual $dg-$ operad.
 Observe that for every
$n$ the degree 0 part, ${\bf D}(P)(n)^0$, of the $dg$ - vector 
space
${\bf D}({\cal P})(n)$ is equal to $F(E^\vee)(n)$, the $n$ - th
space of the free operad generated by the single space $E$, i.e.,
$$ {\bf D}({\cal P})(n)^0 = 
\bigoplus_{{\rm binary.}\atop n - {\rm trees}\,\, T}
 E^\vee(T) \otimes {\rm Det}(T) =
F(E^\vee)(n).$$
We define a morphism of $dg$ - operads $\gamma_P: {\bf D}({\cal P})
 \rightarrow {\cal P}^!$ (here ${\cal P}^!$ is equipped with the
trivial $dg$ - structure) to be given by compositions
$${\bf D}({\cal P})(n) \rightarrow {\bf D}({\cal P})(n) =
F(E^\vee)(n)  \longrightarrow F(E^\vee)(n)/(R^\bot) = {\cal P}^!(n).$$

\proclaim (4.1.2) Lemma. 
For every $n$
the morphism $\gamma_{\cal P}$ induces an isomorphism
$H^0({\bf D}({\cal P})(n)) \rightarrow {\cal P}^!(n)$.

\noindent {\sl Proof:} Observe that the penultimate term
$${\bf D}({\cal P})(n)^{-1} = \bigoplus_T {\cal P}^*(T) \otimes
 {\rm Det}(T)$$
is the sum over the $n$ -  trees $T$ such that all but one vertices
of $T$ are binary and just one vertex is ternary. Note also that
 ${\cal P}(3)^\vee$
is dual to the space of relations of the quadratic operad ${\cal P}^!$.
Thus, the image of the last differential in ${\bf D}({\cal P})(n)$ 
is precisely the space of consequences of the relations in ${\cal P}^!$,
 whence
the statement of the lemma. 

\proclaim (4.1.3) Definition. A quadratic operad $\cal P$ is called Koszul
if the morphism $\gamma_{\cal P}: {\bf D}({\cal P}) \rightarrow {\cal P}^!$
 is a
quasi - isomorphism, i.e., each complex ${\bf D}({\cal P})(n)$
is exact everywhere but the right end. 

\proclaim (4.1.4) Proposition. a) A quadratic operad $\cal P$ is Koszul 
if and only if so is 
${\cal P}^!$.
\hfill\break
b) Let $r$ be the number of simple summands of the algebra
 ${\cal P}(1) = {\cal P}^!(1)^{op}$ and   let $g_{\cal P},
g_{{\cal P}^!}: {\bf C}^r \rightarrow {\bf C}^r$ be the generating maps
of $\cal P$ and ${\cal P}^!$ respectively. 
If $\cal P$ is Koszul,  we have the formal
power series identity 
$$g_{{\cal P}^!}( - g_{\cal P} (-x)) = x, \quad x = (x_1,...,x_r).$$

\noindent {\sl Proof:} a) Form the composition
$${\bf D}({\cal P}^!) \buildrel {\bf D}(\gamma_{\cal P})\over\longrightarrow
{\bf D}({\bf D}({\cal P})) \buildrel f_{\cal P} \over\longrightarrow 
{\cal P},$$
where $f_{\cal P}$ is the quasi - isomorphism constructed in Theorem 3.2.16. 
It is immediate that this composition coincides with
$\gamma_{_{{\cal P}^!}}$. If $\cal P$ is Koszul then $\gamma_{\cal P}$ and
hence ${\bf D}(\gamma_{\cal P})$ are quasi -  isomorphisms,
whence the statement.

\vskip .2cm

b) Follows from Theorem 3.3.2 and the observation that quasi - isomorphic
$dg$ - operads have the same generating maps.

\vskip .3cm

\noindent {\bf (4.1.5)} Let $K$ be a semisimple $k$ - algebra and
$P = \{P(n)\}$ be a $K$ - collection (1.2.11). 
We extend the collection $P$
to a functor $I \mapsto P(I)$ on the category of finite
sets and their bijections as in (1.1.5). For 
any
surjection $f: I \rightarrow J$   of finite sets  we put
$$P(f) = \bigotimes_{j\in J} P(f^{-1}(j)).$$

If, in addition, $P$ is  an operad then for any composable pair
of surjections
$I \buildrel f \over\longrightarrow J
\buildrel g \over\longrightarrow H$ 
the compositions in $P$ give rise to a map
 $$\mu_{g,f}: P(g) 
\otimes P(f) \rightarrow P(gf).$$

Let $P, Q$ be two $K$ - collections. We define a new $K$ - collection
$P(Q)$, called the {\it composition} of $P$ and $Q$,  as follows:
$$P(Q)(n) = \bigoplus_{m=1}^n P(Q)(n)^m, \quad {\rm where}\quad
P(Q)(n)^m = \bigoplus_{f: [n] \rightarrow [m]} \left[
P(m) \otimes_A Q(f)\right]_{\Sigma_m}\leqno {\bf (4.1.6)}$$
and the last sum is taken over all surjections $[n] \rightarrow [m]$.

\proclaim (4.1.7)Proposition. The generating map (3.1.8)
 of  the collection $P(Q)$ is given by the composition
$$g_{P(Q)}(x) = g_P(g_Q(x)), \quad x = (x_1, ..., x_r).$$

\noindent {\sl Proof:} This follows from Lemma 3.3.13.

\vskip .3cm

\noindent {\bf (4.1.8) Koszul complexes.}
Recall [Pr] that the Koszul complex of a quadratic $K$ - algebra $A$ 
is the vector space $A \otimes_K (A^!)^\vee$ (where $(A^!)^\vee = 
{\rm Hom}_{K^{op}} (A^!, K^{op})$) equipped with the differential
defined in a certain natural way. 

\vskip .2cm

Let now ${\cal P} = {\cal P}(K, E, R)$ be a quadratic operad and ${\cal P}^!$
be the dual operad. We define the $n$ - th Koszul complex of $\cal P$
to be the $n$ -th space  of the 
composition ${\cal P}(({\cal P}^!)^\vee)(n)$. Here $({\cal P}^!)^\vee$
is the $K$ - collection consisting of ${\cal P}^!(n)^\vee =
{\rm Hom}_{K^{op}}({\cal P}^!(n), K^{op})$ with the $\Sigma_n$ - action
being the transposed one twisted by sign. We put a grading on the space
${\cal P}(({\cal P}^!)^\vee)(n)$ by (4.1.6) and define the differential $d$
in the following way. 

For $X \in {\cal P}(2)$ and a surjection $g: [m+1] \rightarrow [m]$ 
let 
$$\mu_{g,X}: {\cal P}(m) \rightarrow {\cal P}(m+1)$$
be the operator of composition with $X$ along $g$.

\vskip .2cm

For $\Xi \in {\cal P}^!(2)  = {\cal P}(2)^\vee$ and surjections
 $g: [m+1] \rightarrow [m]$,
$h: [n] \rightarrow [m+1]$  let
$$\mu_{g,h, \Xi}: {\cal P}^!(h) \rightarrow {\cal P}^!(hg)$$
 be the operator of composition with $\Xi$ induced by $\mu_{g,h}$. 

Let us write the identity element
$${\rm Id}_{{\cal P}(2)} \in {\rm Hom}_K({\cal P}(2),
{\cal P}(2)) = {\cal P}(2) \otimes_K {\cal P}^!(2)$$
in the form $\sum X_i \otimes \Xi_i$ where 
$X_1, ..., X_d \in {\cal P}(2)$ and $\Xi_1, ..., \Xi_d \in
{\cal P}^!(2)$. For a surjection $f: [n] \rightarrow [m]$ we define a map
$$d_f: {\cal P}(m) \otimes ({\cal P}^!)^\vee(f) \longrightarrow 
\bigoplus_{h: [n]\rightarrow [m+1]} {\cal P}(m+1) \otimes
 ({\cal P}^!)^\vee(h) \quad {\rm by}$$
$$d_f = \sum_{\matrix{g: [m+1] \rightarrow [m]\cr gh = f}} 
\sum_{i=1}^d \mu_{h, X_i} \otimes \mu_{g,h,\Xi_i}^*.
\leqno {\bf (4.1.9)}$$
Clearly $d_f$ is $\Sigma_n$ - equivariant, hence
factors through  the spaces
of $\Sigma_n$ - coinvariants. Therefore the formula $d =
 \sum_{f: n\rightarrow
m} d_f$ defines a linear map
$$d: {\cal P}(({\cal P}^!)^\vee)(n)^m \rightarrow 
{\cal P}(({\cal P}^!)^\vee)(n)^{m+1}.$$

\vskip .2cm

\proclaim (4.1.10) Proposition. The morphisms $d$
satisfy, for various $m$, 
the condition $d^2 = 0$ thus making each ${\cal P}(({\cal P}^!)^\vee)(n)$
 into 
a complex.

\vskip .2cm

\noindent {\bf (4.1.11)} To prove Proposition 4.1.10, 
we introduce the following notation. 
 Let $E = {\cal P}(2)$ and $F(E)$
be the corresponding free operad. By definition (2.1.1), $F(E)(n) = 
\bigoplus E(T)$ (sum over binary $n$ - trees).
Let  $R \i F(E)(3)$ be the space of
relations of $\cal P$. Call an $n$ -  tree $T$ a {\it 1-ternary tree
} if $T$ contains exactly one ternary vertex $v = v(T)$ , all other vertices
being binary. For such a tree $T$  there are exactly three
binary trees
$T', T'', T'''$  such that $T$ can be obtained from each of them
by contracting an edge. 
We  denote by
 $$R_T \i  E(T') \oplus E(T'') \oplus  E(T''') \i F(E)(n)$$
 the result of substituting $R$ at the place $v$.
Thus,
$${\cal P}(n) = {F(E)(n) \over \sum_{\matrix { 1- {\rm ternary}\cr
n - {\rm trees} \,\,\,T}} R_T}, \quad {\cal P}^!(n)^\vee = \bigcap_
{\matrix {1- {\rm ternary}\cr
n - {\rm trees}\,\,\,T}} R_T.$$
Looking at any summand of 
 ${\cal P}(({\cal P}^!)^\vee)(n)^m$ we find that this is a subquotient
of $F(E)(n)$. More precisely, for a surjection $f: [n]\rightarrow [m]$
let $T(F)$ be the  reduced tree obtained by depicting the map $f$
(and ignoring  vertices with one input (Fig.13).

\vbox to 4cm{}

This tree has one vertex (root) at the bottom and several vertices at the
top. 
We have
$${\cal P}(m) \otimes {\cal P}^!(f) 
 =  {\bigcap_{ 1-{\rm ternary \,\,\,trees}\,\, T:
\atop  T\geq T(f), \,\,v(T) \rightarrow {\rm top}} 
R_T \over
\left(\sum_{ 1-{\rm ternary \,\,\,trees}\,\, T:
\atop  T\geq T(f), \,\,v(T) \rightarrow {\rm bottom}}R_T\right) 
 \cap \left(
 \bigcap_{ 1-{\rm ternary \,\,\,trees}\,\, T:
\atop  T\geq T(f), \,\,v(T) \rightarrow {\rm top}} R_T\right)}
 .\leqno {\bf (4.1.12)}$$

Applying $d$ amounts to deleting some $R_T$ from the intersection
in the numerator of (4.1.12) and
 adding some other $R_T$ to the sum in the denominator,
  according to surjections $g: [m+1] \rightarrow [m]$. Hence applying $d^2$
 to (4.1.12)  moves some  term $R_T$ from the numerator to the denominator.
Therefore $d^2 = 0$. 

\vskip .3cm

Now we state the main result of this section.

\proclaim (4.1.13) Theorem. Let $\cal P$ be a quadratic operad.
The following conditions are equivalent:
\item{(i)} $\cal P$ is Koszul.
\item{(ii)} The Koszul complexes ${\cal P}(({\cal P}^!)^\vee)(n)$ 
are exact for all
$n\geq 2$.

\noindent {\bf (4.1.14)} We   prove Theorem 4.1.13 by reducing it
to a similar result about quadratic associative algebras or, 
rather, quadratic categories, a result which was proven in the
required generality in [BGS]. 

\vskip .2cm

For any operad $\cal P$ we construct, following Mac Lane, a 
certain PROP [Ad]. By definition, this is a category 
Cat$({\cal P})$ whose objects  are symbols
$[n], n=0,1,2,...$. 
The morphisms are defined by
$${\rm Hom}_{{\rm Cat}({\cal P})} ([n], [m])  = \bigoplus_{f: [n]\rightarrow
[m]} {\cal P}(f), \quad {\rm for} \quad n\geq m.$$
where $f$ runs over all surjections $[n] \rightarrow [m]$.
For $n<m$ we set ${\rm Hom}_{{\rm Cat}{\cal P})} ([n], [m]) = 0$.
The composition in the category Cat$({\cal P})$ is induced by maps
$\mu_{g,f}$, see (4.1.5). It is well known [Ad] that Cat$({\cal P})$
has a natural structure of  a symmetric monoidal category
 given  on objects
by $[m] \otimes [n] = [m+n]$. 
We   write  $ {\cal P}(n,m)=
{\rm Hom}_{{\rm Cat}({\cal P})} ([n], [m])$.
 Thus the space ${\cal P}(n)$ of our operad can be
written as ${\cal P}(n,1)$. 

\vskip .3cm

If $\cal P$ is a quadratic operad then Cat$({\cal P})$ is a quadratic
category (${\bf Z}_+$ - algebra) in the sense of [BGS, \S 3] and
 Cat$({\cal P}^!)$ is the dual quadratic
category. Let us compare the cobar - duals and the Koszul duals
for ${\cal P}$ and Cat$({\cal P})$. The cobar - dual category 
of Cat$({\cal P})$ is the
category ${\bf D}({\rm Cat}({\cal P}))$ with the same objects $[m]$
as Cat$({\cal P})$. 
The complex
${\rm Hom}_{{\bf D}({\rm Cat}({\cal P}))} ([n], [m])$ is given by
$${\cal P}(n,m)^* \rightarrow \bigoplus_{n>r>m} {\cal P}(n,r)^*
 \otimes_{{\cal P}(r,r)} {\cal P}(r,m)^* 
\longrightarrow$$
$$\longrightarrow \bigoplus_{n> r_1 > r_2 > m} {\cal P}(n, r_1)^* \otimes_
{{\cal P}(r_1 ,r_1)} {\cal P}(r_1, r_2)^* \otimes_{{\cal P}(r_2,
 r_2)}{\cal P}(r_2, m)^* \rightarrow
...$$
The category Cat$({\cal P})$ is Koszul if and only if each of these complexes
is exact off the rightmost term.
Observe that
$$ {\rm Hom}_{{\bf D}({\rm Cat}({\cal P}))} ([n], [m]) = \left(
{\rm Hom}_{{\bf D}({\rm Cat}({\cal P}))} ([n], [1])\right)^{\otimes m}$$
and 
$${\rm Hom}_{{\bf D}({\rm Cat}({\cal P}))} ([n], [1]) = {\bf D}({\cal P})(n)$$
is the $n$ -th complex of the cobar - operad ${\bf D}({\cal P})$.
This gives the following.

\proclaim (4.1.15) Lemma. 
A quadratic operad ${\cal P}$ is Koszul if and
only if ${\rm Cat}({\cal P})$ is a Koszul quadratic category. 

Let us now recall the interpretation of Koszulness for categories
in terms of Koszul complexes. Let ${\cal Q} = {\cal P}^!$.
For the category Cat$({\cal P})$ we have, according to [BGS, n.4.4],
 the complexes in the form
$$\overfullrule 0cm  
\biggl\{ {\cal Q}(n,m)^* \rightarrow {\cal Q}(n, m+1)^* 
\otimes_{K[\Sigma_{m+1}]}{\cal P}(m+1, m) 
\rightarrow {\cal Q}(n, m+2)^* \otimes_{K[\Sigma_{m+2}]}
{\cal P}(m+2, m) \rightarrow
...$$
$$ \rightarrow {\cal Q}(n, n-1)^*\otimes_{K[\Sigma_{n-1}]}{\cal P}(n-1, m) 
\rightarrow {\cal P}(n,m)\biggl\}.$$
Let us denote the above complex by ${\bf K}^\bullet(n,m)$.
Again, ${\bf K}^\bullet (n,m)$ is the $m$ - th tensor power of
the complex
${\bf K}^\bullet(n,1)$ and ${\bf K}^\bullet(n,1)$ is
 the $n$ - th Koszul complex
of the quadratic operad ${\cal P}$. Thus  Koszulness of ${\cal P}$ is
equivalent to the exactness of all the Koszul
complexes for Cat$({\cal P})$ hence is equivalent to exactness of 
the complexes ${\bf K}^\bullet(n,1)$ only,
i.e., of Koszul complexes of ${\cal P}$. This concludes the proof.

\vskip .3cm

\hfill\vfill\eject

\centerline {\bf 4.2. Homology of algebras over a quadratic operad
and Koszulness.}

\vskip 1cm

\noindent {\bf (4.2.1)} Given a quadratic operad ${\cal P} =
{\cal P}(K, E, R)$ and a $\cal P$ - algebra $A$, we define a chain
complex $(C^{\cal P}_\bullet (A), d)$ as follows. We put
$$C_n(A) = C^{\cal P}_n (A) = \left(A^{\otimes n} \otimes_{K^{\otimes n}}
{\cal P}^!(n)^\vee
\right)_{\Sigma_n}.$$
To define $d_n: C_n (A) \rightarrow C_{n-1}(A)$, we first define a map 
$$\overfullrule 0cm \bar d_n: \left( A^{\otimes n} \otimes \left(
\bigoplus_{{\rm binary}\atop n- {\rm trees}\,\,\,T} E(T)\right)\otimes
{\rm Det} (T)\right)_{\Sigma_n}
\rightarrow \left( A^{\otimes (n-1)} \otimes \left(
\bigoplus_{{\rm binary}\atop (n-1) - {\rm trees}\,\, S} E(S)
\otimes {\rm Det} (S)\right)\right)_
{\Sigma_{n-1}}.$$
If $T$ is an $n$ - tree, we call a vertex $v\in T$ {\it extremal} if all
inputs of $v$ are inputs of $T$. For such a $v$ we define the set $[n]/v$
obtained by replacing the subset In$(v) \i [n]$ by a single
 element. Let $T/v$ be the
$[n]/v$ - tree  obtained by erasing $v$ and replacing it by
an external edge.

Let $T$ be a binary $n$ - tree and $v \in T$ an extremal vertex.
The $\cal P$ - action on $A$ defines a map
$$\bar d_{T, v}: A^{\otimes n} \otimes E(T) \otimes {\rm Det} (T)
\longrightarrow A^{\otimes [n]/v} \otimes E(T/v) \otimes
{\rm Det} (T/v).$$

We define  $\bar d_n$ to be given by the matrix 
with entries $\bar d_{T,v}$.
(Note that the ambiguity in numbering the set $[n]/v$ by $1, 2, ..., n-1$
will disappear after we take coinvariants of $\Sigma_{n-1}$.) 
Observe that $C_n(A) \i \left( A^{\otimes n} \otimes
\bigl( \bigoplus E(T) \otimes {\rm Det} (T)\bigl)\right)_{\Sigma_n}$. 

\proclaim (4.2.2) Proposition. For any $n$ the map $\bar d_n$
takes the subspace $C_n(A)$ into $C_{n-1}(A)$. The maps
$d_n$ defined as restrictions of $\bar d_n$ to $C_n(A)$,
satisfy  $d_{n-1} \circ d_n = 0$. 

\noindent {\sl Proof:} Let  $X$ be an element of $A^{\otimes n} 
\otimes R_T$ for some
1-ternary tree $T$ (here $R_T$ stands for the space
of relations at the ternary vertex of $T$).
 Let $v\in T$ be the ternary vertex and
$T_1, T_2, T_3$ be the binary trees obtained by splitting this
vertex. Suppose first that the vertex $v\in T$ is extremal.
Let $e_i$ be the edge of $T_i$ which is contracted into $v$.
Let also $v_i \in T_i$ be the source of the edge $e_i$.
Then $v_i$ is an extremal vertex. Therefore after erasing
$v_i$ the
edge $e_i$ will become external and cannot be contracted.

\vskip .2cm

If $v$ is not extremal then all the extremal vertices of $T_i$
come from extremal binary vertices of $T$. If $w\in T$ is such a vertex
then $T/w$ is 1-ternary and $\sum_{i=1}^3 \bar d_{T_i,W}(X)$
belongs to $R_{T/w}$. Clearly all 1-ternary $(n-1)$ - trees
can be represented as $T/w$. So if $X\in A^{\otimes n} \otimes R_T$
for all 1-ternary $n$ - trees $T$ then $\bar d_n(X) \in
A^{\otimes (n-1)} \otimes R_S$ for all 1-ternary $(n-1)$ - trees
$S$. This shows that $\bar d_n(C_n(A)) \i C_{n-1}(A)$.

\vskip .2cm

Let us prove that $\bar d_{n-1} \circ \bar d_n = 0$ on $C_n(A)$. Let
$T$ be a 1-ternary tree $T$ whose 1-ternary vertex $v$ is extremal. 
Suppose that $X \in A^{\otimes n} \otimes R_T \i \bigoplus_{i=1}^3 
A^{\otimes n} \otimes E(T_i)$ where $T_i$, are as before. 
Denote by $w_i$ the end of the edge $e_i$. Then $w_i$ will
become extremal after erasing $v_i$. Thus
$$\sum_{i=1}^3 \bar d_{T_i/e_i, w_i} \left( \bar d_{T_i, v_i}(X)
\right)$$
is a sum of quantities like $r(a,b,c)$ with $r\in R \i F(E)(3)$
and $a,b,c \in A$. Since $A$ is a $P$ - algebra and $R$ is the
space of relations for $P$, every such quantity is equal to 0.
This implies that $\bar d_{n-1}(\bar d_n(X)) = 0$
if $X \in \bigcap_T A^{\otimes n} \otimes R_T$.
The proposition follows.

\proclaim (4.2.3) Definition. The complex $C_\bullet(A)  = \bigoplus
C_n(A)$ with the differential $d$ defined above
is called the chain complex of the $\cal P$ - algebra $A$. Its homology
will be denoted $H_n(A)$ or $H_n^{\cal P}(A)$.

\noindent {\bf (4.2.4) Examples.} a) Let ${\cal P} = {\cal A}s$ be the 
associative operad. 
The dual operad ${\cal P}^!$ is isomorphic to $\cal P$ itself.
Let $A$ be an associative algebra.
Since ${\cal P}^!(n)$ is the regular representation of
 $\Sigma_n$, we find that the
complex $C_\bullet (A)$ has the form
$$ ... \rightarrow A\otimes A \otimes A \longrightarrow
A\otimes A \longrightarrow A.$$
A straightforward calculation of the differential shows that it
is the standard Hochschild complex of $A$. Thus,
 $H_i^{\cal P}(A)$ is the
Hochschild homology of $A$.

\vskip .2cm

b) Let ${\cal P} = {\cal L}ie$ be the Lie operad. Then ${\cal P}^! = 
{\cal C} om$
is the commutative operad; in particular, ${\cal P}^!(n)$
is the trivial 
1-dimensional $\Sigma_n$ - module and ${\cal P}^!(n)^\vee$
is the sign representation. Let ${\cal G}$ be a Lie algebra. By the
above, its chain complex, as a $\cal P$ - algebra, has the form
$$...\rightarrow \Lambda ^3 {\cal G} \longrightarrow \Lambda^2 {\cal G}
\longrightarrow {\cal G}.$$
Again, a straightforward calculation of the differential shows that this
is the standard Chevalley - Eilenberg chain complex of ${\cal G}$, so that
$H_i^{\cal P}({\cal G}) = H_i({\cal G}, k)$ is the Lie algebra homology of
${\cal G}$ with constant coefficients.

\vskip .2cm

c) Similarly, if ${\cal P} = {\cal C}om$ is the commutative operad, and
$A$ is a commutative algebra, one finds that $H_i^{\cal P}(A)$ is the
Harrison homology of $A$, see e.g., [Lo]. 

\proclaim (4.2.5) Theorem. Let $\cal P$ be a quadratic operad. 
Then $\cal P$ is
Koszul if and only if for any free $\cal P$ - algebra $F = F_{\cal P}(V)$
we have $H_i^{\cal P}(F) = 0$ for $i>0$.

\noindent {\sl Proof:} Let $F_n$  denote the free $\cal P$
 - algebra on generators
$x_1, ..., x_n$. This algebra has an obvious $({\bf Z}_+)^n$ - grading
such that deg $(x_i) = (0,..., 1, ..., 0)$ (the unit on the $i$ - th place)
and ${\rm deg} \, (\mu(a,b)) = {\rm deg} (a) + {\rm deg} (b)$ for any
$a,b \in F_n, \mu \in {\cal P}(2)$. 

The chain complex $C_\bullet^{\cal P}(F_n)$ also inherits this grading. The
following result which is immediate from the definitions implies the
''if" part of Theorem 4.2.5.

\proclaim (4.2.6) Proposition. The $n$ - th Koszul complex 
$K^{\cal P}_\bullet (n)$
of $\cal P$ is isomorphic to the 
multihomogeneous part of $C_\bullet^{\cal P}
(F_n)$
of multi - degree  $(1,...,1)$. 

Before proving the ``only if" part of (4.2.5), let us mention the
following corollary from what we have already done.

\proclaim (4.2.7) Corollary. The operads ${\cal A}s, {\cal C}om, {\cal L}ie$
are Koszul. 

\noindent {\sl Proof:} It is well known that a free associative algebra 
without unit has
 higher Hochschild homology trivial and a free Lie algebra has
higher homology with constant coefficients trivial. This proves that
${\cal A}s$ and ${\cal L}ie$ are Koszul. The case of ${\cal C}om$
follows by duality (4.1.4). 

\vskip .2cm

Another situation of applicability of Theorem 4.2.5 is provided by the
work of E. Getzler and J.D.S. Jones [Ge J]. They proved,
 in the context of their work
on iterated integrals, the vanishing
of the homology of free algebras over the operads associated to
homology of configuration space. 

\vskip .3cm

\noindent {\bf (4.2.8) End of the proof of Theorem 4.2.5.} Let $A$ be any
$\cal P$ - algebra. Along with the chain complex $C_\bullet^{\cal P}(A)$
we introduce the ``big" chain complex $BC_\bullet^{\cal P} (A)$. By
definition,
$$BC_n^{\cal P}(A) = \bigoplus_{i+j=n\atop j\leq 0} \left(
\bigoplus_{ i - {\rm trees} \,\,T \atop |T| = i-2+j}
 A^{\otimes i} \otimes {\cal P}
(T) \otimes {\rm Det} (T)\right)_{\Sigma_i}. \leqno{\bf (4.2.9)}$$
The differential $d: BC_n^{\cal P}(A) \rightarrow BC_{n-1}^{\cal P}(A)$
is given by the following two types of matrix elements:
$$d_{T,e}: A^{\otimes i} \otimes {\cal P}(T) \otimes {\rm Det}(T)
\rightarrow A^{\otimes i} \otimes {\cal P}(T/e)\otimes {\rm Det}(T/e)$$
defined for any $i$ - tree $T$ and any internal edge $e\in {\rm Ed}(T)$,
and $$d_{T,v}: A^{\otimes i} \otimes {\cal P}(T) \otimes {\rm Det}
(T) \rightarrow
A^{\otimes [i]/v} \otimes {\cal P}(T/v) \otimes {\rm Det} (T/v)$$
defined for any $i$ - tree $T$ and any extremal vertex $v\in T$. 

\vskip .2cm

The operator $d_{T,e}$ is $1\otimes \mu_{T,e} \otimes l_e^*$ where 
$\mu_{T,e}: {\cal P}(T) \rightarrow {\cal P}(T/e)$ is induced by
the composition in $\cal P$ and the map
$l_e^*: {\rm Det} (T) \rightarrow {\rm Det} (T/e)$ is dual to the
exterior multiplication by $e$. 
The operator $d_{T,v}$ is induced, in a similar way, by the $\cal P$ - action
on $A$.

\vskip .2cm

It is immediate to verify that $d^2=0$. The embedding
${\cal P}^!(n)^\vee \i \bigoplus_{ {\rm binary}\atop {\rm trees}\,\,T} R(T)$
gives rise to an embedding of complexes
$$j: C_\bullet^{\cal P}(A) \rightarrow BC_\bullet^{\cal P}(A).
\leqno {\bf (4.2.10)}$$
The complex $BC_\bullet^{\cal P}(A)$ has an increasing filtration $F$
with the $m$ - th term
 $F_m BC_n^{\cal P}(A)$ being the sum of the summands in (4.2.9)
with $i\leq m$. The associated graded complex has the form
$${\rm gr}_m^F BC^{\cal P}_\bullet (A) = A^{\otimes m} \otimes
{\bf D}({\cal P})(m)^\vee$$
where {\bf D} is the duality for $dg$ - operads (3.2.13). 
This implies the following result.

\proclaim (4.2.11) Proposition. If the operad $\cal P$ is
Koszul then the embedding (4.2.11) is a quasi - isomorphism for any
$\cal P$ - algebra $A$.

 To complete the proof of Theorem 4.2.5,
it is enough, in view of Proposition 4.2.11, to establish the following.

\proclaim (4.2.11) Proposition. Let $\cal P$ be any $k$ - linear operad,
$K = {\cal P}(1)$ and
and $A = F_{\cal P}(V)$ be the free $\cal P$ - algebra generated by a
finite - dimensional $K$ - module $V$. Then
$$H_i(BC_\bullet^{\cal P}(A)) = \cases{ 0,& $i>0$\cr
$V$,& $i=0$.\cr}$$

\noindent {\sl Proof:} To each tree $T$ and a function $\nu: {\rm In}(T)
\rightarrow \{ 1,2,3, ...\}$ we associate the vector space
$$C(T,\nu) = {\rm Det} (T) \otimes \bigotimes_{v\in T} {\cal P}({\rm In}(v))
\otimes\bigotimes_{i\in {\rm In}(T)} A_{\nu(i)}$$
where $A_n, n= 1,2,3,...$,  denotes the degree $n$ component of $A$.
By construction the graded vector space $BC_\bullet^{\cal P} (A)$
is the direct sum of spaces $C(T, \nu)$ over all
(isomorphism classes of) pairs $(T,\nu)$. 

\vskip .2cm

Given a pair $(T, \nu)$ we call an input edge $e\in {\rm In}(T)$
{\it non- degenerate} if $\nu(e) > 1$. Given a non- degenerate edge
$e \in {\rm In}(T)$ we define the {\it degeneration} of $(T,\nu)$
along $e$ as the pair $(\hat T, \hat\nu)$ where $\hat T$ is obtained
by attaching at $e$ a star (1.1.1) with $\nu(e)$ inputs.
The function $\hat\nu$ is set equal 1 on the new inputs
and remain unchanged on other inputs.  A pair $(T,\nu)$ is called
{\it non- degenerate} if for any extremal vertex $v\in T$ there is 
an edge $e
\in {\rm In} (v)$ such that $\nu (e) > 1$. A non - degenerate pair
cannot be obtained by degeneration.

\vskip .2cm

 Let Nd$(T, \nu)$ denote the
number of non-degenerate input edges of $T$. We introduce an increasing
filtration $G$ on $BC_\bullet^{\cal P}(A)$ by putting
$$G_i BC_\bullet^{\cal P}(A) = \bigoplus_{\# ({\rm vertices \,\,\,\, of}
\,\,\,T) + {\rm Nd} (T, \nu) \leq i} C(T,\nu).$$
The differential preserves the filtration. Furthermore, the complex
${\rm gr}_i^G BC_\bullet (A)$ splits into a direct sum of
complexes $E_\bullet(T,\nu)$
labelled by non- degenerate pairs $(T, \nu)$. These complexes have
the form
$$E_\bullet (T, \nu) = \bigoplus_{(\hat T, \hat \nu)} C(\hat T, \hat\nu)$$
where the sum runs over all possible pairs $(\hat T, \hat\nu)$
obtained from $(T, \nu)$ by (iterated) degeneration. Such pairs
$(\hat T, \hat\nu)$ are parametrized by all possible subsets of the set
$S$ of non - degenerate inputs of $T$. It is clear that the corresponding 
spaces
$C(\hat T, \hat\nu)$ are all the same. Observe that subsets of $S$
correspond to faces of a simplex with $|S|$ vertices (denote this simplex
by $\Delta$). Moreover, we find that $E_\bullet(T, \nu)$ is the
tensor product of the fixed vector space $C(T, \nu)$ and the augmented
chain complex of $\Delta$. This complex is acyclic unless $S = \emptyset$.
The only non- degenerate pair $(T, \nu)$ with $S = \emptyset$
consists of $T = \{\rightarrow\}$ (the tree with no vertices)
and $\nu = 1$. This gives $H_0(BC_\bullet (A)) = V$, and the theorem 
follows.

\vskip .3cm

\noindent {\bf (4.2.12) Homotopy $\cal P$ - algebras.}
Let ${\cal P} = {\cal P}(K, E, R)$
 be a Koszul quadratic operad and ${\cal P}^!$ its quadratic dual.
Then the canonical morphism of operads (4.1.1)
$\gamma_{\cal P}: {\bf D}({\cal P}^!) \rightarrow {\cal P}$ is
a quasi - isomorphism. Hence, any $\cal P$ - algebra can be viewed as
a ${\bf D}({\cal P}^!)$ - algebra. This motivates the following

\proclaim (4.2.13) Definition. A $dg$ - algebra over 
${\bf D}({\cal P}^!)$ is called a homotopy $\cal P$ - algebra.

\proclaim (4.2.14) Proposition. Let $A = \bigoplus A_n$ be a graded $K$
- bimodule with ${\rm dim} (A_n) < \infty$ for all $n$. A
Giving astructure of
a homotopy $\cal P$ - algebra on $A$ is the same as 
giving a non - homogeneous differential
$d$ on the free algebra $F_{{\cal P}^!}(A^*[1])$ satisfying the conditions:
\item{(i)} $d^2 = 0$;
\item{(ii)} $d$ is a derivation with respect to any binary operation
in ${\cal P}^!$, i.e., we have
$$d(\mu(a,b)) = \mu(d(a), b) + (-1)^{{\rm deg}(a)} \mu (a, d(b)),\quad
\mu \in {\cal P}(2), a,b \in F_{{\cal P}^!}(A^*[1]).$$

\noindent {\sl Proof:} Let $d$ be a differential in the free algebra
$F = F_{{\cal P}^!}(A^*[1])$ satisfying (i) and (ii). 
Recall that the free algebra has a natural grading $F = \bigoplus_{n\geq 1}
F_n$ where (1.3.5)
$$F_n = \left( {\cal P}^!(n) \otimes_{K^{\otimes n}} (A^*)^{\otimes n}
\right)_{\rm sgn} \leqno {\bf (4.2.15)}$$
(the subscript ``sgn" stands for the ``anti - invariants"
of the symmetric group action). We decompose the differential $d$ into
homogeneous components $d = d_1 + d_2 + ...$ where the component
$d_i$ shifts degree by $i-1$, i.e., $d_i: F_j \rightarrow
F_{j+i-1}, \,\,\forall j$. The equality $d^2 = 0$ yields, in
particular $d_1^2 = 0$. 

\vskip .2cm

Observe that the operad ${\cal P}^!$ being quadratic, the
algebra $F$ is generated by its degree 1 component
$F_1 = {\cal P}^1 (1) \otimes_K A^* = A^*$. Hence, the differential $d$
is completely determined, due to the property (ii), by its restriction
to $F_1$.  Separating the degrees, we see that this restriction
is given by a collection of maps $d_m: A^* = F_1 \rightarrow F_m$,
$m=1,2,...$. In view of (4.2.15), we may view $d_m$ as a
$\Sigma_m$ - invariant element
$$d_m \in  {\rm Det} (k^m) \otimes
{\cal P}^!(m) \otimes_{K^{\otimes m}}(A^*)^{\otimes m}
\otimes_K A.\leqno {\bf (4.2.16)}$$
Next, we replace, using (1.1.6), the integer $m$
(4.2.16) by any $m$ - element set. Further, for any $n$ - tree $T$,
we form the tensor product of the elements $d_m$ over all vertices of $T$
to get an element $\bigotimes_{v\in T} d_{{\rm In}(v)}$. Rearranging
the factors in the tensor product, we get
$$\bigotimes_{v\in T} d_{{\rm In}(v)} \quad \in \quad
{\rm Det}(T) \otimes 
\bigotimes_{v\in T} \left( {\cal P}^!({\rm In}(v)) \otimes_K
\biggl(\bigotimes_{e\in {\rm In}(v)} A^*\biggl) \otimes_K A\right)
\leqno{\bf (4.2.17)}$$
where the factor $A$ on the right corresponds to the output edge of $v$
and Det$(T)$ was introduced in (3.2.0).
Each internal edge of $T$ occurs in the above tensor product 
twice, once as an input at some vertes, contributing to the
factor $A^*$ and once as an output at some vertex, contributing
to the factor $A$. Contracting the pairs of factors $A^*$ and $A$
corresponding to each internal edge and using the notation
(1.2.13), we obtain from (4.2.17) a well defined element
$$d(T) \quad \in \quad 
{\rm Det}
 (k^n) \otimes {\cal P}^!(T) \otimes_K (A^*)^{\otimes n} \otimes A.$$
This element can be regarded a morphism
$$d(T): {\cal P}^!(T)^\vee =
{\cal P}^!(T)^* \otimes {\rm Det}(k^n)
 \longrightarrow {\rm Hom}(A^{\otimes n}, A).
\leqno {\bf (4.2.18)}$$
The morphisms (4.2.18), assembled for various $n$ - trees together,
define, for each $n$, a morphism
$${\bf D}({\cal P}^!)(n) \rightarrow {\rm Hom} (A^{\otimes n}, A) =
{\cal E}_A(n).$$
One can check by a direct calculation that these morphisms give rise to
a morphissm of operads ${\bf D}({\cal P}^!) \rightarrow {\cal E}_A$
if and only if the original differential $d$ on the free algebra $F$
satisfies the property (i) of Proposition 4.2.13. That makes $A$
a ${\bf D}({\cal P}^!)$ - algebra. The opposite implication is proved
by reversing the above argument.

\vskip .3cm

\noindent {\bf (4.2.19) Examples.} a) 
Let ${\cal P} = {\cal L}ie$ be the Lie
operad so ${\cal P}^! = {\cal C}om$. Since ${\cal C}om (n) = k$ for any
$n$, a structure of a ${\bf D}({\cal C}om)$ - algebra on a $dg$ - 
vector space $A$ is determined by specification of $n$ - ary
antisymmetric operations $[x_1, ..., x_n]$ for any $n\geq 2$ (these
operations correspond to the basis vectors of ${\cal C}om (n)$).
Proposition 4.2.13 in this case gives the equivalence of two
definitions [SS] of a homotopy Lie algebra: first as a vector space with
brackets $[x_1,...,x_n]$ satisfying the generalized Jacobi
identity and, the second, as a vector space $A$ with a differential
on the exterior algebra $\bigwedge^\bullet (A^*)$, satisfying
the Leibnitz rule. 

 b) If we take ${\cal P} = {\cal C}om$, we get a notion of a 
homotopy ${\cal C}om$ - algebra which is a $dg$ - algebra over
${\bf D}({\cal L}ie)$. Such algebras are particular cases of
algebras ``associative and commutative up to all higher homotopies"
(May algebras) [HS 1] [KM 1,2]. More precisely, a homotopy
${\cal C}om$ - algebra is strictly commutative and
equipped with a system of natural homotopies which
ensure, in particular, the associativity of the cohomology
algebra. This structure is not the same as a
homotopy between
$ab$ and $ba$ in an associative $dg$ - algebra,
 the data often referred to as "homotopy 
commutativity".

\hfill\vfill\eject

\centerline {\bf References.}

\vskip 1cm

\item{[Ad]} J. F. Adams, Infinite loop spaces. Princeton Univ. Press,
1978.

\item{[Ar 1]} V.I.  Arnold, On the representation of a continuous function of
three variables by superpositions of functions of two variables,
{\it Mat. Sbornik}, {\bf 48} (1959) 3 - 74 (in Russian).

\item{[Ar 2]} V.I.  Arnold, Topological invariants of algebraic
functions II, {\it Funct. Anal. Appl.} {\bf 4}
(1970) 91 - 98. 

\item{[BCW]} H. Bass, E.H. Connell, D. Wright, The jacobian conjecture:
reduction of degree and formal expansion of the inverse, {\it Bull. AMS}
{\bf 7} (1982), 287 - 330.

\item{[Be]} A. Beilinson, Coherent sheaves on $P^n$ and
problems of linear algebra, {\it Funct. Anal. Appl.} {\bf 12} (1978),
214 - 216. 

\item{[BGG]} I.N.  Bernstein, I. M. Gelfand, S.I. Gelfand, Algebraic
vector bundles on $P^n$ and problems of linear algebra,
{\it Funct. Anal. Appl.} {\bf 12} (1978), 212 - 214.

\item{[BGS]} A.  Beilinson, V.  Ginzburg, V.
 Schechtman, Koszul duality,
{\it J. of Geometry and \hfill\break
Physics}, {\bf 5} (1988) 317 - 350.

\item{[BGSo]} A. Beilinson, V. Ginzburg, W. Soergel, Koszul duality patterns
in representation theory, {\it J. AMS}, to appear. 

\item{[BG 1]} A.  Beilinson, V.  Ginzburg, Infinitesimal structure of
moduli spaces of G - bundles, {\it Int. Math. Research Notices}
(appendix to {\it Duke Math. J.}) 1992, p. 63 - 74.

\item{[BG 2]}  A.  Beilinson, V.  Ginzburg, Resolution of diagonals,
homotopy algebra and moduli spaces, preprint 1993.

\item{[BV]} J.M.  Boardman, R. Vogt, Homotopy invariant algebraic structures
on topological spaces, Lect. Notes in Math., {\bf 347}, Springer -
Verlag, 1973.

\item{[C]} F.R.  Cohen, The homology of $C_{n+1}$ - spaces, in:
Lect. Notes in Math., {\bf 533}, p. 207 - 351, Springer - Verlag, 1976.

\item{[D 1]} P. Deligne, Th\'eorie de Hodge II, {\it Publ. IHES}, 
{\bf 40} (1971) 5-58.

\item{[D 2]}  P. Deligne, Resum\'e des premi\'ers expos\'es de
A. Grothendieck, in: SGA 7, Expos\'e I, Lecture Notes in Math., 
{\bf 288}, p. 1-24, Springer - Verlag, 1972.

\item{[DL 1]} P. Deligne, G. Lusztig, Representations of reductive
groups over finite fields, {\it Ann. Math.} {\bf 103} (1976) 103 -
161.

\item{[DL 2]} P. Deligne, G. Lusztig, Duality for representations of
a reductive group over a finite field I, {\it J. of Alg.} {\bf 74}
(1982), 284 - 291.

\item{[DL 3]} P. Deligne, G. Lusztig, Duality for representations of
a reductive group over a finite field II, {\it J. of Alg.}
{\bf 81} (1983), 540 - 549. 

\item{[Dr]} V.G.  Drinfeld, letter to V.V.  Schechtman, September
1988 (unpublished).

\item{[ESV]} H.Esnault, V.V.Schechtman, E.Vieweg,
Cohomology of local systems on the complement of hyperplanes,
{\it Invent. Math.}, {\bf 109} (1992), 557 - 561.

\item{[F]} Z. Fiedorowicz, The symmetric bar - construction,
preprint 1991.

\item{[FM]} W. Fulton, R.D.  MacPherson, Compactification
 of configuration 
spaces, {\it Ann. Math.}, to appear.

\item{[G]} E. Getzler, Batalin - Vilkovisky algebras and
two - dimensional topological field theories, preprint 1992.

\item{[Ge J]} E. Getzler, J.D.S. Jones, $n$ - algebras, in preparation.

\item{[Go]} R. Godement, Topologie algebrique and th\`eorie des faisceaux,
Hermann, Paris, 1962. 

\item{[Gd]} I. Good, Generalization to several variables of
Lagrange's expansion with application to stochastic processes,
{\it Proc. Camb. Phil. Soc.}, {\bf 56} (1960) 367 - 380.

\item{[GM]} M. Goreski, R.D.  MacPherson, Stratified Morse theory,
Springer - Verlag 1990.

\item{[HS 1]} V.A. Hinich, V.V. Schectman, On the homotopy limit of
homotopy algebras, in: Lecture Notes in Math., {\bf 1289}, p. 240 - 264,
Springer - Verlag 1987.

\item{[HS 2]} V.A. Hinich, V.V. Schectman, Homotopy Lie algebras,
preprint 1993.

\item{[Jo]} S. A. Joni, Lagrange's inversion in higher dimension and
umbral operators, {\it Linear and Multilinear Alg.} {\bf 6}
(1978) 111 - 121.

\item{[Ka 1]} M.M.  Kapranov, On the derived categories and the K - functor
of coherent sheaves on intersections of quadrics, {\it Russian Math.
Izv.} {\bf 32} (1989) 191- 204. 

\item{[Ka 2]} M.M.  Kapranov, Veronese curves and Grothendieck - Knudsen
moduli space $\overline{M_{0,n}}$, {\it J. Algebr. Geometry},
{\bf 2} (1993), to appear.

\item{[Ka 3]} M.M. Kapranov, Chow quotients of Grassmannians I, preprint
1992.

\item{[Ka 4]} M.M.Kapranov, On the derived categories of coherent sheaves
on some homogeneous spaces, {\it Invent. Math.}, {\bf 92} (1988),
479 - 508.

\item{[KS]} M. Kashiwara, P. Schapira, Sheaves on Manifolds
(Grundlehren der Math. Wiss {\bf 292}), Springer -
Verlag 1991.

\item{[Kl]} A.A.  Klyachko, Lie elements in the tensor algebra,
{\it Siberian Math. J.}, {\bf 15} (1974), 914 - 920. 

\item{[Kn]} F.F.  Knudsen, The projectivity of the moduli space of
stable curves II. The stacks $\overline {M_{g,n}}$, {\it Math.
Scand.,} 52 (1983), 161 - 189. 

\item{[Kol]} A.N. Kolmogoroff, On  the representation of
 continuous functions of
several variables as superpositions of functions of smaller
number of variables, {\it Sov. Math. Dokl.} {\bf 108}
(1956) n. 2, p. 179 - 182 (in Russian), reprinted in his
Selected Papers (V. Tikhomirov Ed.) vol.I, p. 378 - 382, 
Kluwer Publ., 1991.

\item{[Kon 1]} M.L.  Kontsevich, Intersection theory on moduli spaces
and matrix Airy function, {\it Comm. Math. Phys.},  1992.

\item{[Kon 2]} M.L.  Kontsevich, Formal (non- commutative)
differential geometry, preprint 1992.

\item{[Kon 3]} M.L.  Kontsevich, Graphs, homotopy Lie algebras and
low - dimensional topology, preprint 1992.  

\item{[KM 1]} I. Kriz, J.P. May, Operads, algebras and modules I.
Preprint 1993.

\item{[KM 2]} I. Kriz, J.P. May, Differential graded algebras up to
homotopy and their derived categories, preprint 1992.

\item{[L]} M. Lazard, Lois de groupes et analyseurs, {\it Ann. ENS},
{\bf 62} (1955) 299 - 400.

\item{[Li Z]} B. Lian, G. Zuckerman, New perspectives on the BRST - 
Algebraic structures of string theory, preprint  1992. 

\item{[Lo]} J.L. Loday, Cyclic homology (Grund. Math. Wiss. {\bf 301}), 
Springer - Verlag 1992.

\item{[L\"o]} C. L\"ofwall, On the sublgebra generated by one -
dimensional elements in the Yoneda Ext - algebra, in: Lect. Notes in
Math., {\bf 1183} p. 291 - 338, Springer - Verlag 1986.

\item{[Lu]} G. Lusztig, The discrete series of $GL_n$ over a finite
field (Ann. of Math. Studies {\bf 81}) Princeton Univ. Press 1974.

\item{[MTWW]} J.H. MacKay, J. Towber, S. S.-S. Wang, D. Wright,
Reversion of a system of power series with application to
combinatorics, Preprint 1992. 

\item{[Mac]} S. Mac Lane, Natural associativity and commutativity,
{\it Rice Univ. Studies} {\bf 49} (1963), 28 - 46, reprinted in his
Selected papers, p. 415 - 433, Springer - Verlag 1979.

 \item{[Macd]} I. Macdonald, Symmetric functions and Hall
polynomials, Oxford 1979.

\item{[Man]} Y.I.  Manin, Some remarks on Koszul algebras and 
quantum groups, {\it Ann. Inst. Fourier} {\bf 37}
(1987) 191 - 205.

\item{[May]} J. P. May, Geometry of iterated loop spaces, Lect. Notes
in Math. {\bf 271}, Springer - Verlag 1972.

\item{[Mo]} J.C. Moore, Differential homological algebra,
in: {\it Proc. Congr\`es Int. Math. Nice 1970}, t.1, p. 335 - 339,
Gauthier - Villars, Paris 1971.

\item{[Pe]} R.C. Penner, The 
decorated Teichm\"uller space of  puctured surfaces,
{\it Comm. Math. Phys.}, {\bf 113} (1987), 299 - 339.

\item{[Q 1]} D. Quillen, Rational homotopy theory, {\it Ann. Math.} 
{\bf 90} (1969) 205 - 295.

\item{[Q 2]} D. Quillen, On the cohomology of commutative rings,
{\it Proc. Symp. Pure Math.} {\bf 17} (1970)
65 - 87.

\item{[Pr]} S. Priddy, Koszul resolutions, {\it Trans. AMS} {\bf 152}
(1970) 39 - 60. 

\item{[SV]} V.V. Schechtman, A.N. Varchenko, Arrangements of
hyperplanes and Lie algebra cohomology, {\it Inv. Math.} {\bf 106}
(1991) 139 - 194.

\item{[SS]} M. Schlesinger, J.D. Stasheff,
The Lie algebra structure of tangent cohomology
and deformation theory, {\it J. Pure and Appl . Alg.}, {\bf 38} (1985),
313 - 322.

\item{[St]} J.D. Stasheff, Homotopy associativity of H - spaces,
{\it Trans. AMS} {\bf 108} (1963) 275 - 312.

\item{[Sy]} J.J. Sylvester, On the change of systems of independent
variables, {\it Quart. J. of Math.} {\bf 1} (1857), 42 - 56. 

\item{[W]} D. Wright, The tree formula for reversion of power series,
{\it J. Pure Appl  Alg.} {\bf 57} (1989), 191 - 211.

\vskip 2cm

{\sl Authors'  addresses:

V.G.: Department of Mathematics, University of Chicago, Chicago IL
60637, email: 
 beigin@cpd.landau.free.msk.su, ginzburg@zaphod.uchicago edu

M.K.: Department of Mathematics, Northwestern University,
Evanston Il 60208, email: kapranov@chow.math.nwu.edu

\bye